\documentclass[sn-mathphys,Numbered]{sn-jnl}

\usepackage{graphicx}%
\usepackage{amsmath,amssymb,amsfonts,bm}%
\usepackage{multirow}
\usepackage{amsthm}%
\usepackage{mathrsfs}%
\usepackage[title]{appendix}%
\usepackage{xcolor}%
\usepackage{textcomp}%
\usepackage{manyfoot}%
\usepackage{booktabs}%
\usepackage{algpseudocode}%
\usepackage{listings}%
\usepackage{bigints}
\usepackage{outlines}
\usepackage{geometry}
\usepackage{subfigure}
\usepackage{siunitx}
\usepackage{float}
\usepackage{comment}

\usepackage[onehalfspacing]{setspace}

\usepackage[utf8]{inputenc}
\usepackage{graphicx}
\usepackage[ruled,vlined]{algorithm2e}
\usepackage{amsmath}
\usepackage[section]{placeins}
\usepackage{enumitem}
\usepackage{color,soul}
\usepackage{xfrac}

\newcommand{\bb}{\boldsymbol}
\usepackage{nicefrac}
\usepackage{mathtools}

\begin{document}

\title[Article Title]{A Jacobian-free Newton-Krylov method for cell-centred finite volume solid mechanics}

\author*[1]{\fnm{Philip} \sur{Cardiff}}\email{philip.cardiff@ucd.ie}
\author[1]{\fnm{Dylan} \sur{Armfield}}
\author[2]{\fnm{\v{Z}eljko} \sur{Tukovi\'{c}}}
\author[1,2]{\fnm{Ivan} \sur{Batisti\'{c}}}

\affil*[1]{\orgdiv{School of Mechanical and Materials Engineering}, \orgname{University College Dublin}, \orgaddress{\country{Ireland}}}
\affil[2]{\orgdiv{Faculty of Mechanical Engineering and Naval Architecture}, \orgname{University of Zagreb}, \orgaddress{\country{Croatia}}}

\abstract
{
This study proposes a Jacobian-free Newton-Krylov approach for finite-volume solid mechanics.
Traditional Newton-based approaches require explicit Jacobian matrix formation and storage, which can be computationally expensive and memory-intensive.
In contrast, Jacobian-free Newton-Krylov methods approximate the Jacobian’s action using finite differences, combined with Krylov subspace solvers such as the generalised minimal residual method (GMRES), enabling seamless integration into existing segregated finite-volume frameworks without major code refactoring.
This work proposes and benchmarks the performance of a compact-stencil Jacobian-free Newton-Krylov method against a conventional segregated approach on a suite of test cases that span varying geometric dimensions, nonlinearities, dynamic responses, and material behaviours.
Key metrics, including computational cost, memory efficiency, and robustness, are evaluated, along with the influence of preconditioning strategies and stabilisation scaling.
Results show that the proposed Jacobian-free Newton-Krylov method outperforms the segregated approach in all linear and nonlinear elastic cases, achieving order-of-magnitude speedups in many instances; however, divergence is observed in elastoplastic cases, highlighting areas for further development.
It is found that preconditioning choice affects performance: a LU direct solver is fastest for small to moderately sized cases, while a multigrid method is more effective for larger problems.
The findings demonstrate that Jacobian-free Newton-Krylov methods are promising for advancing finite-volume solid mechanics simulations, particularly for existing segregated frameworks where minimal modifications enable their adoption.
The described implementations are available in the solids4foam toolbox for OpenFOAM, inviting the community to explore, extend, and compare these procedures.
}


\maketitle

\section{Introduction}\label{sec:intro}
%
Finite volume formulations for solid mechanics are heavily influenced by their fluid mechanics counterparts, favouring fully explicit \citep{Trangenstein1991, Kluth2010, Lee2013, Haider2017} or segregated implicit \citep{Demirdzic1988, Fryer1991, Demirdzic1995, Jasak2000, Tukovic2013, Cardiff2017, Tukovic2018, Batistic2022} methods.
Segregated approaches, where the governing equations are temporarily decomposed into scalar component equations, offer memory efficiency and implementation simplicity, but the outer-coupling Picard iterations often exhibit slow convergence.
Explicit formulations are straightforward to implement and offer superior robustness, but are only efficient for high-speed dynamics, where the physics requires small time increments.
In contrast, the finite element community commonly employs Newton-Raphson-type solution algorithms, which require repeated assembly of the Jacobian matrix and the solution of the resulting block-coupled, non-diagonally dominant linear system.
A disadvantage of traditional Newton-based approaches is that they typically require explicit formation and storage of the Jacobian matrix, which can be computationally expensive and memory-intensive.
A further disadvantage from a finite-volume perspective is that extending existing code frameworks from segregated algorithms to coupled Newton-Raphson-type approaches is challenging due to the required assembly, storage, and solution of the resulting block-coupled system.
In addition, the derivation of the true Jacobian matrix is non-trivial.
Consequently, similar block-coupled finite-volume solution methods are rare in the literature \citep{Das2011, Cardiff2016, Castrillo2024}.
The motivation for the current work is to assess the robustness and efficiency of block-coupled Newton-Raphson approaches that can be easily incorporated into existing segregated solution frameworks.
To this end, the current article proposes and examines the efficacy of \emph{Jacobian-free} Newton-Krylov methods, where the quadratic convergence of Newton methods can be achieved without deriving, assembling, or storing the exact Jacobian.

Jacobian-free Newton-Krylov methods circumvent the need for the Jacobian matrix by combining the Newton-Raphson method with Krylov subspace iterative linear solvers, such as the generalised minimal residual method (GMRES), which do not explicitly require the Jacobian matrix.
Instead, only the action of the Jacobian matrix on a solution-type vector is required.
The key step in Jacobian-free Newton-Krylov methods is the approximation of products between the Jacobian matrix and a vector using the finite difference method; that is
\begin{eqnarray}
	\bb{J} \bb{v} \approx \frac{\bb{R}(\bb{u} + \epsilon \bb{v}) - \bb{R}(\bb{u})}{\epsilon}
\end{eqnarray}
where $\bb{J}$ is the Jacobian matrix, $\bb{R}$ is the residual function of the governing equation, $\mathbf{u}$ is the current solution vector (e.g. nodal displacements), $\mathbf{v}$ is a vector (e.g., from a Krylov subspace), and $\epsilon$ is a small scalar perturbation.
With an appropriate choice of $\epsilon$ (balancing truncation and round-off errors), the characteristic quadratic convergence of Newton methods can be achieved without the Jacobian, hence the modifier \emph{Jacobian-free}.
This approach promises significant memory savings over Jacobian-based methods, especially for large-scale problems, and potentially also reduced execution time, with an appropriate choice of solution components.

A crucial aspect of ensuring the efficiency and robustness of the Jacobian-free Newton-Krylov method is the selection of an appropriate preconditioner for the Krylov iterations.
This preconditioner is often derived from the exact Jacobian matrix in traditional Newton methods.
However, the Jacobian-free approach does not allow direct access to the full Jacobian matrix, necessitating an alternative strategy to approximate its action.
To this end, and to extend existing segregated frameworks, this work proposes using a compact-stencil \emph{approximate} Jacobian as the preconditioner. This approximate Jacobian corresponds to the matrix typically employed in segregated solid mechanics approaches; similar approaches are successful in fluid mechanics applications \citep{Mchugh1994, Qin2000, Geuzaine2001, Pernice2001, Knoll2004, Nejat2008, Vaassen2008, Lucas2010, Nejat2011, Nishikawa2020}; however, it is unclear if such an approach is suitable for solid mechanics - a question this work aims to answer.
By leveraging this compact-stencil approximate Jacobian, it is aimed to effectively precondition the Krylov iterations, enhancing convergence while maintaining the memory and computational savings that define the Jacobian-free and segregated methods.
Similarly, if such an approach is efficient, it would naturally fit into existing segregated frameworks, as existing matrix storage and assembly can be reused.


As noted by \citet{Knoll2004}, Jacobian-free Newton–Krylov methods first appeared in the 1980s and early 1990s for the solution of ordinary and partial differential equations \citep{Gear1983, Chan1984, Brown1986, Brown1990}.
\citet{Mchugh1994} and co-workers demonstrated the potential of Jacobian-free Newton-Krylov methods for the solution of steady, incompressible Navier-Stokes problems using a staggered finite-volume formulation.
They found the GMRES linear solver to be faster than the conjugate gradient solver; however, the true Jacobian matrix was still evaluated by finite differencing during preconditioner construction.
\citet{Qin2000} later extended Jacobian-free Newton–Krylov methods to unsteady compressible Reynolds-averaged Navier–Stokes equations.
Once again, GMRES proved more robust than a conjugate gradient linear solver.
For preconditioning, a fully matrix-free approach was proposed, based on the approximate factorisation procedure of \citet{Badcock1996}.
\citet{Geuzaine2001} examined the use of lower and higher-order discretisations for the preconditioning matrix for steady, compressible high Reynolds number flows.
They found that the lower-order compact stencil preconditioner performed best when combined with the GMRES linear solver and the ILU($k$) preconditioner.
As is common with other Newton approaches, a pseudo-transient algorithm and mesh sequencing strategy were used to improve the convergence.
With the aim of exploiting mature segregated solution approaches, \citet{Pernice2001} recast the SIMPLE segregated algorithm as a method that operates on residuals, allowing it to serve as a preconditioner for a Jacobian-free Newton-Krylov approach.
It was found that the SIMPLE-preconditioned Jacobian-free Newton-Krylov substantially accelerated convergence for the incompressible Navier–Stokes problems examined.
GMRES was adopted as the linear solver, and the importance of choosing an appropriate restart parameter was noted.
They also found that a multigrid approach was more efficient than ILU($k$) for preconditioning, with the implementation based on the PETSc \citep{PETSc} package.

Over the subsequent two decades, Jacobian-free Newton-Krylov methods have seen increasing application in the field of finite-volume computational fluid dynamics, though they remain far from widespread.
The use of low-order approximate Jacobians for preconditioning has been established as an efficient and robust approach, e.g. \citep{Nejat2008, Vaassen2008, Nejat2011, Nishikawa2020}.
A particularly attractive application of the low-order preconditioning approach is for higher-order discretisations; that is, spatial discretisation of order greater than two.
\citet{Nejat2008} demonstrated such a higher-order approach for steady, inviscid compressible flows, with up to fourth-order accuracy.
GMRES was the adopted linear solver with the ILU($k$) preconditioner and fill-in values ($k$) ranging from 2 to 4. 
\citet{Nejat2011} later applied the second-order variant of the approach to non-Newtonian flows and found Newton-GMRES with ILU(1) to be the most efficient combination.
Like \citet{Vaassen2008} and many other authors, Nejat and co-workers \citep{Nejat2008, Nejat2011} proposed a special start-up phase to encourage the solution to efficiently reach the Newton method's domain of quadratic convergence.
The use of a pseudo-transient approach combined with an adaptive time step (e.g., based on the convergence of the Newton iterations) is common.
A further promotion of Jacobian-free Newton-Krylov methods in fluid dynamics was provided by \citet{Lucas2010}, who demonstrated order-of-magnitude speedups over nonlinear multigrid methods.
A particular benefit was the insensitivity of the Jacobian-free Newton-Krylov methods to changes in mesh aspect ratio, density, time step and Reynolds number.
Another relevant point is the effect of numerical stabilisation (damping) on convergence: \citet{Nishikawa2017} demonstrated that a Jacobian-free Newton-Krylov approach was less sensitive to damping magnitude and remained stable at lower damping levels than a segregated approach.

In recent years, Jacobian-free Newton-Krylov methods have gained increasing traction in finite volume discretisations across various engineering and physics applications.
For turbomachinery simulations, \citet{Zhang2024} developed a Jacobian-free Newton-Krylov method to improve convergence difficulties associated with simplified implicit iterative algorithms, employing a startup strategy based on an approximate Jacobian system coupled with an adaptive Courant-Friedrichs-Lewy number scheme before transitioning to full Newton iterations.
Their approach leverages GMRES with a graph-colouring-based finite-difference preconditioner.
Similarly, \citet{Zhang2024twofluid} introduced a staggered-grid JFNK solver for the two-fluid six-equation model, specifically designed to address discontinuities and ill-conditioning, using GMRES for the Krylov subspace iteration.
In the context of Reynolds-Averaged Navier-Stokes equations, \citet{Sukas2025} utilised PETSc-based \citep{PETSc} Jacobian-free Newton-Krylov solvers with both Jacobian-free finite difference and direct Jacobian-vector products, highlighting the versatility of the approach.
\citet{Johannes2021} explored a hybrid approach where a high-order Discontinuous Galerkin method was preconditioned using a finite volume-based multigrid solver, demonstrating effectiveness in unsteady three-dimensional compressible flow applications.
Beyond traditional fluid dynamics, Jacobian-free Newton-Krylov methods have been successfully extended to reactive transport in porous media, as demonstrated by \citet{Amir2021}.
Jacobian-free Newton-Krylov has also been applied in magnetohydrodynamics, where 
\citet{Nguyen2022} proposed higher-order spatial and temporal discretisation coupled with additive Schwarz preconditioning based on block incomplete LU factorisation, and \citet{Chacon2025} introduced a physics-based preconditioner leveraging Schur complement operator splitting, coupled with Newton-flexible-GMRES to improve convergence.

Beyond individual applications, several multiphysics simulation frameworks incorporate Jacobian-free Newton-Krylov as a core numerical strategy.
The Eilmer code \citep{Eilmer2023}, an open-source finite-volume solver for hypersonic flow simulations, relies on Jacobian-free Newton-Krylov methods for robust convergence in high-speed aerothermodynamics.
In the finite element community, MOOSE \citep{Giudicelli2024, Zhu2025, Wu2024} has emerged as a leading multiphysics framework that employs Jacobian-free Newton-Krylov via PETSc \citep{PETSc}, enabling highly scalable simulations across disciplines.
Built on MOOSE, the BISON nuclear fuel performance code \citep{Williamson2021} uses Jacobian-free Newton-Krylov to handle tightly coupled thermomechanical and radiation transport effects in large-scale nuclear reactor simulations.
Other FE-based solvers, such as deal.II \citep{Africa2024, Munch2024} also integrate Jacobian-free Newton-Krylov with PETSc \citep{PETSc} for efficient nonlinear solvers.
The adoption of JFNK in large-scale multiphysics codes underscores its scalability: matrix-free residual evaluations scale well on distributed-memory systems, and Krylov solvers can leverage multigrid and domain-decomposition preconditioners to keep computational costs manageable.

Jacobian-free Newton-Krylov methods have also been extended to finite-difference and nodal-integral methods.
\citet{Ahmed2025} explored preconditioned Newton-Krylov approaches for Navier-Stokes equations using a nodal integral method with GMRES-based physics-aware preconditioning, demonstrating improved robustness for complex flow problems.
In contrast, \citet{Liu2025} employed an explicit Jacobian-based Newton-Krylov solver for nuclear multiphysics coupling, arguing that the true Jacobian provides a better preconditioner, thereby enhancing stability compared to a purely Jacobian-free approach.
Similarly, JFNK has been extended to immersed boundary finite element methods for cardiac mechanics \citep{Ma2024}, where BiCGSTAB was preferred over GMRES to ensure scalability of computational time and memory requirements.

Although they are increasingly used in finite-volume computational fluid dynamics, Jacobian-free Newton-Krylov methods have yet to be applied to finite-volume solid mechanics.
This article addresses this point by being the first (to the authors' knowledge) to propose a Jacobian-free Newton-Krylov approach for finite volume solid mechanics.
By benchmarking the Jacobian-free Newton-Krylov solver on various structural mechanics problems – including static and dynamic cases, linear elasticity and nonlinear material behaviour – we assess its performance relative to the traditional segregated approach.
Particular focus is given to the use of a compact-stencil approximate Jacobian preconditioner, inspired by existing segregated solid mechanics procedures, and ease of implementation into an existing segregated finite volume framework -- the OpenFOAM-based solids4foam toolkit \citep{Cardiff2018}.
This study explicitly addresses whether the known benefits of the Jacobian-free Newton-Krylov method (fast convergence, stability, memory savings) can be realised for finite volume solid mechanics, and what challenges might arise (for example, we report and discuss robustness issues in elastoplastic cases, pointing to areas needing further research).
The proposed Jacobian-free Newton-Krylov approach is benchmarked against a segregated solution procedure, which remains the standard in finite-volume solid (and fluid) mechanics and forms the backbone of widely used codes such as OpenFOAM.
While the limitations of segregated approaches in terms of efficiency and robustness are well known, they can still be competitive - or even preferable - in situations where Newton-type methods may be overly expensive or insufficiently robust. For example, \citet{Whelan2025} recently demonstrated that a segregated procedure was more robust than a Newton method for large-strain elastoplastic damage problems. 
Against this backdrop, a primary motivation of the current study is to assess whether Jacobian-free Newton–Krylov methods can be adopted within such segregated frameworks with minimal disruption, thereby leveraging existing infrastructure while potentially improving robustness and efficiency

The remainder of the paper is structured as follows:
Section 2 summarises a typical solid-mechanics mathematical model and its cell-centred finite-volume discretisation.
Section 3 presents the solution algorithms, starting with the classic segregated method and then the proposed Jacobian-free Newton-Krylov method.
The performance of the proposed Jacobian-free Newton-Krylov approach is compared with the segregated approach across several benchmark cases in Section 4, where the effects of several factors are examined. 
Finally, the article concludes with a summary of the main conclusions.

\section{Mathematical Model and Numerical Methods}\label{sec:math_model}

%

\subsection{Governing Equations} \label{sec:governing_eqn}

In this work, interest is restricted to Lagrangian formulations of the conservation of linear momentum.
Assuming small strains, the linear geometry formulation is expressed in strong integral form as:
\begin{eqnarray} \label{eqn:momentum_lingeom}
    \int_{\Omega} \rho \frac{\partial^2 \bb{u} }{\partial t^2} \, d\Omega
    =
    \oint_{\Gamma} \bb{n} \cdot \bb{\sigma}_s \,  d\Gamma
    + \int_{\Omega}  \bb{f}_b \, d\Omega
\end{eqnarray}
where $\Omega$ is the volume of an arbitrary body bounded by a surface $\Gamma$ with outwards pointing normal $\bb{n}$.
The density is $\rho$, $\bb{u}$ is the displacement vector, $\bb{\sigma}_s$ is the engineering (small strain) stress tensor, and $\bb{f}_b$ is a body force per unit volume, e.g., $\rho \bb{g}$, where $\bb{g}$ is gravity.
Total and updated Lagrangian formulations, suitable for finite strains, are shown in Appendix \ref{app:TL_UL}.

The governing equations are complemented by boundary conditions, which are classified into three types: prescribed displacement, prescribed traction, and symmetry.
The definition of the engineering stress ($\bb{\sigma}_s$) and true stress ($\bb{\sigma}$) in Equation \ref{eqn:momentum_lingeom} (and Equations \ref{eqn:momentum_TL} and \ref{eqn:momentum_UL} in Appendix \ref{app:TL_UL})  is given by a chosen mechanical law.
Five mechanical laws are considered in this work, as briefly outlined in Appendix \ref{app:mechLaws}: linear elasticity (Hooke's law), three forms of hyperelasticity (St.\,Venant-Kirchhoff, neo-Hookean, and Guccione), and neo-Hookean $J_2$ hyperelastoplasticity.

\subsection{Newton-Type Solution Methods}
To facilitate the comparison between the classic segregated solution algorithm and the proposed Jacobian-free Newton-Krylov algorithm, the governing linear momentum conservation (Equations \ref{eqn:momentum_lingeom}, \ref{eqn:momentum_TL} and \ref{eqn:momentum_UL}) is expressed in the general form:
\begin{eqnarray} \label{eqn:residual}
	    \bb{R}(\bb{u}) = \bb{0}
\end{eqnarray}
where $\bb{R}$ represents the \emph{residual} (imbalance) of the equation, which is a function of the primary unknown field.
For example, in the linear geometry case, the residual is given as
\begin{eqnarray}
    \bb{R}(\bb{u})
    \;=\;
    \oint_{\Gamma} \bb{n} \cdot \bb{\sigma}_s(\bb{u}) \,  d\Gamma
    + \int_{\Omega}  \rho \bb{g} \, d\Omega
    -  \int_{\Omega} \rho \frac{\partial^2 \bb{u} }{\partial t^2} \, d\Omega
    \;=\; \bb{0}
\end{eqnarray}
where the dependence of the stress tensor on the solution vector is made explicitly clear: $\bb{\sigma}_s(\bb{u})$.

In Newton-type methods, a Taylor expansion about a current point $\bb{u}_k$ can be used to solve Equation \ref{eqn:residual} \cite{Knoll2004}:
\begin{eqnarray}
	\bb{R}(\bb{u}_{k+1}) = \bb{R}(\bb{u}_{k}) \;+\;  \bb{R}'(\bb{u}_{k}) (\bb{u}_{k+1} - \bb{u}_{k}) \;+\; \text{H.O.T.} = \bb{0}
\end{eqnarray}
Neglecting the higher-order terms ($\text{H.O.T.}$) yields the strict Newton method in terms of an iteration over a sequence of linear systems: 
\begin{eqnarray} \label{eq:NewtonRaphson}
	\bb{J}(\bb{u}_k) \delta \bb{u} &=& -\bb{R}(\bb{u}_k), \notag \\
	\bb{u}_{k+1} &=& \bb{u}_k + s \, \delta \bb{u}, \notag \\
	\quad
	k &=& 0,1,...
\end{eqnarray}
where $\bb{J} \equiv \bb{R}' \equiv \partial \bb{R}/\partial \bb{u}$ is the Jacobian matrix.
Starting the Newton procedure requires specifying $\bb{u}_0$.
The scalar $s > 0$ can be chosen to improve convergence, for example, via a line search or an under-relaxation/damping procedure, and is set to unity in the classic Newton-Raphson approach.
Iterations are performed over this system until the residual $\bb{R}(\bb{u}_k)$ and solution correction $\delta \bb{u}$ are sufficiently small, with appropriate normalisation.

For problems with $N$ scalar equations and $N$ scalar unknowns, the residual $\bb{R}$ and solution $\bb{u}$ vectors have dimensions of $N \times 1$. 
The components of the $N \times N$ Jacobian are
\begin{eqnarray} \label{eq:J}
	{J}_{ij} = \frac{\partial {R}_i (\bb{u})}{\partial u_j}
\end{eqnarray}

The current work focuses on vector problems, where the governing momentum equation is formulated in terms of the unknown displacement solution vector.
In this case, Equation \ref{eq:J} refers to the individual scalar components of the residual, solution, and Jacobian.
That is, for 3-D analyses, the residual takes the form
\begin{eqnarray}
	\bb{R}(\bb{u}) = \left\{ R_1^x, R_1^y, R_1^z, R_2^x, R_2^y, R_2^z, ..., R_n^z \right\}
\end{eqnarray}
and the solution takes the form
\begin{eqnarray}
	\bb{u} = \left\{ u_1^x, u_1^y, u_1^z, u_2^x, u_2^y, u_2^z, ..., u_n^z \right\}
\end{eqnarray}
In practice, it is often more practical and efficient to form and store the residual, solution and Jacobian in a \emph{blocked} manner, where the residual and solution can be considered as vectors of vectors.
Similarly, the Jacobian can be formed in terms of sub-matrix block coefficients.

In the strict Newton procedure, the residuals converge at a quadratic rate when the current solution is close to the true solution; that is, the iteration error decreases proportionally to the square of the error at the previous iteration.
Once the method gets sufficiently close to the true solution, the number of correct digits in the approximation roughly doubles with each iteration. 
However, quadratic convergence is only possible when using the exact Jacobian.
In contrast, a quasi-Newton method uses an approximation to the Jacobian, sacrificing strict quadratic convergence in exchange for a more computationally efficient overall procedure.
From this perspective, the segregated solution algorithm commonly employed in finite volume solid mechanics can be viewed as a quasi-Newton method, where an approximate Jacobian replaces the exact Jacobian: 
\begin{eqnarray} \label{eq:Seg}
    \bb{\tilde{J}}(\bb{u}_k) \;\delta \bb{u} = -\bb{R}(\bb{u}_k)
\end{eqnarray}
In this case, the approximate Jacobian $\bb{\tilde{J}}$ comes from the inertia term and compact stencil discretisation of a simple diffusion (Laplacian) term.
A benefit of this approach is that the inter-component coupling is removed from the Jacobian, allowing the solution of three smaller scalar systems rather than one larger vector system in 3-D (or two smaller systems in 2-D).

A fully explicit procedure can also be viewed from this perspective by selecting a \emph{diagonal} approximate Jacobian $\bb{\tilde{D}}$ (only the inertia term), making the solution of the linear system trivial:
\begin{eqnarray} \label{eq:exp}
    \bb{\tilde{D}}(\bb{u}_k) \;\delta \bb{u} = -\bb{R}(\bb{u}_k)
\end{eqnarray}

\subsection{Cell-Centred Finite Volume Discretisation}
\label{sec:discretisation}
\begin{figure}[htbp]
	\centering
   		\includegraphics[width=\textwidth]{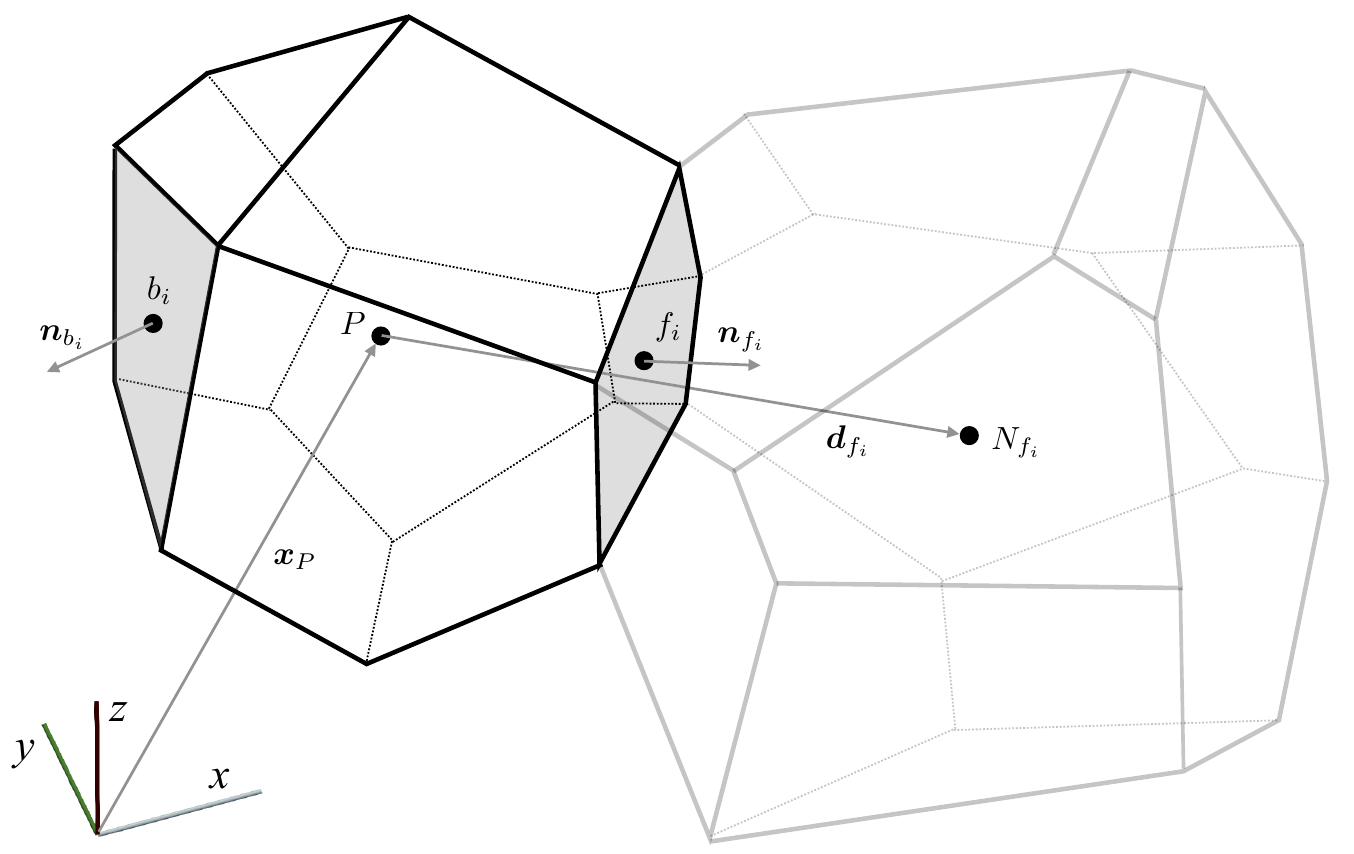} 
	\caption{Representative convex polyhedral cell $P$ and neighbouring cell $N_{f_i}$, which share a face $f_i$.}
	\label{fig:cell}
\end{figure}

In this work, a nominally second-order cell-centred finite volume discretisation is employed.
The solution domain is discretised in both space and time.
The total simulation period is divided into a finite number of time increments, denoted as $\Delta t$, and the discretised governing momentum equation is solved iteratively in a time-marching fashion.
The spatial domain is partitioned into a finite set $\mathcal{P}$ of contiguous convex polyhedral cells, where each cell is denoted by $P \in \mathcal{P}$.
The number of cells in the mesh is indicated by $|\mathcal{P}|$.
A representative cell $P$ is shown in Figure~\ref{fig:cell}.
The set of all faces of cell $P$ is denoted by $\mathcal{F}_P$. This set is further subdivided into two disjoint subsets:
\begin{itemize}
    \item \textbf{Internal Faces} ($\mathcal{F}^{\text{int}}_P$): Faces that are shared with neighbouring cells.
    \item \textbf{Boundary Faces} ($\mathcal{F}^{\text{bnd}}_P$): Faces that lie on the boundary of the spatial domain.
    The boundary faces of cell $P$ are further classed into three disjoint sets, $\mathcal{F}^{\text{bnd}}_P \coloneqq \mathcal{F}^{\text{disp}}_P \cup \mathcal{F}^{\text{trac}}_P \cup \mathcal{F}^{\text{symm}}_P$, representing boundary faces where displacement ($\mathcal{F}^{\text{disp}}_P$), traction ($\mathcal{F}^{\text{trac}}_P$) and symmetry ($\mathcal{F}^{\text{symm}}_P$) conditions are prescribed.
\end{itemize}
Each internal face $f_i \in \mathcal{F}^{\text{int}}_P$ corresponds to a neighbouring cell $N_{f_i} \in \mathcal{N}_P$, where $\mathcal{N}_P$ is the set of all neighbouring cells of $P$. The outward unit normal vector associated with an internal face $f_i$ is denoted by $\mathbf{n}_{f_i}$, while the outward unit normal vector associated with a boundary face $b_i$ is denoted by $\mathbf{n}_{b_i}$.
The vector $\mathbf{d}_{f_i}$ connects the centroid of cell $P$ with the centroid of the neighbouring cell $N_{f_i}$, whereas the vector $\mathbf{d}_{b_i}$ connects the centroid of cell $P$ with the centroid of boundary face ${b_i}$.
For convenience, we will also define the set of all faces of cell $P$ excluding those on a traction boundary as $\mathcal{F}_P^{\text{non-trac}} \coloneqq \mathcal{F}_P^{\text{int}} \cup \mathcal{F}_P^{\text{disp}} \cup \mathcal{F}_P^{\text{symm}}$.

The conservation equation (Equations \ref{eqn:momentum_lingeom}, \ref{eqn:momentum_TL}, or \ref{eqn:momentum_UL}) is applied to each cell $P$ and discretised in terms of the displacement at the centroid of the cell $\boldsymbol{u}_P$ and the displacements $\boldsymbol{u}_{N_{f_i}}$ at the centroids of the neighbouring cells.
Proceeding with the discretisation, the volume and surface integrals in the governing equation are approximated by algebraic equations as described below.

\subsubsection{Volume Integrals}
To discretise the volume integrals, the integrand $\bb{\phi}$ is assumed to locally vary according to a truncated Taylor series expansion about the centroid of cell $P$:
\begin{eqnarray}
	\bb{\phi}(\bb{x})  \approx \bb{\phi}_P + (\bb{x} - \bb{x}_P) \cdot \nabla \bb{\phi}_P
\end{eqnarray}
where subscript $P$ indicates a value at the centroid of the cell $P$.
Consequently, volume integrals over a cell $P$ can be approximated to second-order accuracy as
\begin{eqnarray} \label{eq:volume_integral}
	\int_{\Omega_P} \bb{\phi} \, d \Omega_P
		&\approx& \int_{\Omega_P}  \left[ \bb{\phi}_P + (\bb{x} - \bb{x}_P) \cdot \nabla \bb{\phi}_P \right] d \Omega_P \notag \\
		&\approx& \bb{\phi}_P \Omega_P
\end{eqnarray}
where $\Omega_P$ is the volume of cell $P$ and $\int_{\Omega_P} (\bb{x} - \bb{x}_P) d\Omega_P \equiv 0$ by definition of the cell centroid.
This approximation corresponds to the midpoint rule and one-point quadrature.

Using Equation \ref{eq:volume_integral}, the inertia term (e.g. left-hand side term of Equation \ref{eqn:momentum_lingeom}) becomes
\begin{eqnarray} \label{eq:inertia}
	\int_{\Omega_P} \rho \frac{\partial \bb{u} }{\partial t}  d \Omega_P
	\;&\approx&\;
	\rho_P \left(\frac{\partial^2 \bb{u} }{\partial t^2}\right)_P  \Omega_P
\end{eqnarray}
Similarly, the body force term (e.g. the second term on the right-hand side of Equation \ref{eqn:momentum_lingeom}) becomes:
\begin{eqnarray}
	\int_{\Omega_P} \, \rho \, \bb{g} \,  d \Omega_P
	\;&\approx&\;
	\rho_P \, \bb{g}\,  \Omega_P
\end{eqnarray}
The discretisation of the acceleration in Equation \ref{eq:inertia} can be achieved using one of many finite difference schemes, e.g. first-order Euler, second-order backwards, or second-order Newmark-beta.
In the current work, the second-order backwards (BDF2) scheme is used:
\begin{eqnarray} \label{eq:inertia2}
	\left(\frac{\partial^2 \boldsymbol{u}_P}{\partial t^2}\right)_P
	&\approx& \frac{3\boldsymbol{v}_P^{[t+1]} - 4\boldsymbol{v}_P^{[t]} + \boldsymbol{v}_P^{[t-1]}}{2\Delta t} \notag \\
	&\approx&
	\frac{3\left( 
		\dfrac{3\boldsymbol{u}_P^{[t+1]} - 4\boldsymbol{u}_P^{[t]} + \boldsymbol{u}_P^{[t-1]}}{2\Delta t} 
		\right) 
	- 4\boldsymbol{v}_P^{[t]} + \boldsymbol{v}_P^{[t-1]}}{2\Delta t}
\end{eqnarray}
where $\Delta t$ is the time increment -- assumed constant here.
Superscript $[t]$ indicates the time level, with $\bb{u}_P^{[t+1]}$ corresponding to the unknown displacement at the current time step.
 The velocity vector $\bb{v} = \partial \bb{u}/\partial t$ at the current time step is also updated using the BDF2 scheme as
 \begin{eqnarray}
	\boldsymbol{v}_P^{[t+1]}	&\approx&
		\dfrac{3\boldsymbol{u}_P^{[t+1]} - 4\boldsymbol{u}_P^{[t]} + \boldsymbol{u}_P^{[t-1]}}{2\Delta t} 
\end{eqnarray}
Consequently, the displacement and velocity at the two previous time steps must be stored, or alternatively, the displacement at the previous four time steps.

\subsubsection{Surface Integrals}
The surface integral term can be discretised using one-point quadrature at each face $f_i$ as
\begin{eqnarray} \label{eq:divStressDiscret}
	\oint_{\Gamma_P} \bb{n} \cdot \bb{\sigma}  \; d\Gamma_P
	&=& \sum_{f_i \in \mathcal{F}_P} \int_{\Gamma_{f_i}} \bb{n} \cdot \bb{\sigma}  \,  d \Gamma_{f_i} \notag \\
	&\approx&
	\sum_{f_i \in \mathcal{F}_P^{\text{int}}} \bb{\Gamma}_{f_i} \cdot \bb{\sigma}_{f_i}
	+ \sum_{d_i \in \mathcal{F}_P^{\text{disp}}} \bb{\Gamma}_{d_i} \cdot \bb{\sigma}_{P}
	+ \sum_{s_i \in \mathcal{F}_P^{\text{symm}}} \bb{\Gamma}_{s_i} \cdot \bb{\sigma}_{s_i}
	+ \sum_{t_i \in \mathcal{F}_P^{\text{trac}}} |\bb{\Gamma}_{t_i}| \bar{\bb{T}}_{t_i}
\end{eqnarray}
where $\Gamma_P$ indicates the surface of cell $P$, $\bb{\Gamma}_{f_i}$ indicates the area vector of face $f_i$, and vector $\bar{\bb{T}}_{t_i}$ represents the prescribed traction on the traction boundary face $t_i$.

It is noted that the displacement is assumed to vary linearly within each cell; hence, the displacement gradient and stress are constant within each cell.
In the current work, a unique definition of stress $\bb{\sigma}_{f_i}$ at each cell face $f_i$ is given as a weighted-averaged of values in the two cells ($\bb{\sigma}_P $, $\bb{\sigma}_{N_{f_i}}$) straddling the face \citep{Jasak1996}:
\begin{eqnarray} \label{eq:stressInterp}
	\bb{\sigma}_{f_i} &=& w_{f_i} \bb{\sigma}_P + (1 - w_{f_i}) \bb{\sigma}_{N_{f_i}}
\end{eqnarray}
where the interpolation weight is defined as $w_{f_i} = (\bb{n}_{f_i} \cdot [\bb{x}_{N_{f_i}} - \bb{x}_{f_i}])/(\bb{n}_{f_i} \cdot [\bb{x}_{N_{f_i}} - \bb{x}_{P} ])$; however, achieving second-order accuracy of the displacement field is independent of the value of the weights, and, for example, $w_{f_i} = \nicefrac{1}{2}$ would also be sufficient.

The stress $\bb{\sigma}_{s_i}$ at a symmetry boundary face is calculated as
\begin{eqnarray} \label{eq:symm}
	\bb{\sigma}_{s_i}
		&=& \frac{1}{2} \left (\bb{\sigma}_P + \bb{R}_{s_i} \cdot \bb{\sigma}_{P} \right) \notag \\
		&=& \left (\textbf{I} - \bb{n}_{s_i} \otimes \bb{n}_{s_i} \right) \cdot \bb{\sigma}_P
\end{eqnarray}
where $\bb{R}_{s_i} \cdot \bb{\sigma}_{P}$ represents the mirror reflection of $\bb{\sigma}_P$ across the symmetry boundary face $s_i$.
The reflection tensor is $\bb{R}_{s_i} = \textbf{I} - 2 \bb{n}_{s_i} \otimes \bb{n}_{s_i}$ \citep{Demirdzic2022}, with $\bb{n}_{s_i}$ indicating the unit normal of the symmetry boundary face $s_i$.
From Equation \ref{eq:symm}, it is clear that shear stresses are zero on a symmetry plane boundary face.
Note from Equation \ref{eq:divStressDiscret} that the stress on displacement boundary faces is assumed to be equal to the stress $\bb{\sigma}_P$ at the centroid of cell $P$.

The cell-centred stress $\bb{\sigma}_P$ is calculated as a function of the displacement gradient according to the chosen mechanical law, for example, as shown in Appendix \ref{app:mechLaws}.
The presented discretisation is second-order accurate in space for displacement if the cell-centred displacement gradients (and the stress) are at least first-order accurate, even if the cell faces are not flat.
To achieve this, the cell-centred displacement gradients are determined using a weighted first-neighbours least squares method \citep{Jasak1996},
\begin{eqnarray} \label{eq:leastSquaresGrad}
	\left(\bb{\nabla}\bb{u}\right)_P
		&=&\sum_{f_i \in \mathcal{F}^{\text{int}}_P} w_{f_i} |\bb{\Gamma}_{f_i}| \frac{\bb{G}^{-1}_P \cdot \bb{d}_{f_i}}{\bb{d}_{f_i} \cdot \bb{d}_{f_i}}  \otimes \left(\bb{u}_{N_{f_i}} - \bb{u}_P \right) \notag \\
		&&+ \sum_{s_i \in \mathcal{F}^{\text{symm}}_P} w_{f_i} |\bb{\Gamma}_{s_i}| \frac{\bb{G}^{-1}_P \cdot \bb{d}_{s_i}}{\bb{d}_{s_i} \cdot \bb{d}_{s_i}} \otimes \left(\bb{R}_{s_i} \cdot \bb{u}_P - \bb{u}_P \right)
\end{eqnarray}
which is exact for linear functions.
The vector $\bb{u}_{b_i}$ indicates the displacement at the centroid of boundary face ${b_i}$, while vector $\bb{d}_{d_i}$ connects the centroid of cell $P$ to the centroid of displacement boundary face $d_i$.
The quantity $\bb{R}_{s_i} \cdot \bb{u}_P$ represents the mirror reflection of $\bb{u}_P$ across the symmetry boundary face $s_i$.
The vector $\bb{d}_{s_i}$ connects the centroid $\bb{x}_P$ of cell $P$ with its mirror reflection $\bb{R}_{s_i}  \cdot \bb{x}_P$ through boundary face $s_i$.
Traction and displacement boundary faces are excluded in Equation \ref{eq:leastSquaresGrad}; this is in contrast to the default approach in OpenFOAM \citep{Jasak2011}.
The effect of the prescribed boundary displacements is introduced through the compact component of the Rhie-Chow stabilisation term, described in Section \ref{sec:RhieChow}.
The $\bb{G}_P$ tensor for cell $P$ is calculated as
\begin{eqnarray}
	 \bb{G}_P &=&
	 \sum_{{f_i} \in \mathcal{F}^{\text{int}}_P} (1 - w_{f_i}) |\bb{\Gamma}_{f_i}|  \frac{\bb{d}_{f_i} \otimes \bb{d}_{f_i}}{\bb{d}_{f_i} \cdot \bb{d}_{f_i}} \notag \\
	 && +  \sum_{{s_i} \in \mathcal{F}^{\text{symm}}_P} (1 - w_{s_i}) |\bb{\Gamma}_{s_i}|  \frac{\bb{d}_{s_i} \otimes \bb{d}_{s_i}}{\bb{d}_{s_i} \cdot \bb{d}_{s_i}}
\end{eqnarray}
As the $w |\bb{\Gamma}| \ \bb{G}^{-1}_P \cdot \bb{d}/(\bb{d}\cdot \bb{d})$ vectors in Equation \ref{eq:leastSquaresGrad} are purely a function of the mesh, they can be computed once (or each time the mesh moves) and stored.
Equation \ref{eq:leastSquaresGrad} approximates the cell-centre gradients to at least a first-order accuracy, increasing to second-order accuracy on certain smooth grids \citep{Syrakos2023};
first-order accurate gradients are sufficient to preserve second-order accuracy of the cell-centre displacements.

If required, the displacement $\bb{u}_{t_i}$ on a traction boundary face $t_i$ can be calculated by extrapolation from the centre of cell $P$ as
\begin{eqnarray}
	\bb{u}_{t_i} = \bb{u}_P + \bb{d}_{t_i} \cdot \left(\bb{\nabla} \bb{u} \right)_P
\end{eqnarray}
where $\bb{d}_{t_i}$ represents the vector from the centroid of cell $P$ to the centroid of the traction boundary face $t_i$.
Similarly, if required, the displacement $\bb{u}_{s_i}$ at a symmetry plane face $s_i$ is calculated using the same approach as Equation \ref{eq:symm}:
\begin{eqnarray} 
	\bb{u}_{s_i}
		&=&  \frac{1}{2} \left( \bb{u}_P + \bb{R}_{s_i} \cdot \bb{u}_P \right) \notag \\
		&=& \left (\textbf{I} - \bb{n}_{s_i} \otimes \bb{n}_{s_i} \right) \cdot \bb{u}_P
\end{eqnarray}

\subsubsection{Rhie-Chow Stabilisation} \label{sec:RhieChow}
To quell zero-energy solution modes (i.e. checkerboarding oscillations), a Rhie-Chow-type stabilisation term \cite{Rhie1983} is added to the residual (Equation \ref{eqn:residual}).
The Rhie-Chow stabilisation term  $\mathcal{D}_P^{\text {Rhie-Chow }}$ for a cell $P$ takes the following form:
\begin{eqnarray} \label{eq:RhieChow}
	\mathcal{D}_P^{\text {Rhie-Chow}}
	&=& \sum_{f_i \in \mathcal{F}^{\text{int}}_P} \alpha \bar{K}_{f_i} \left[
		\left|\bb{\Delta}_{f_i} \right| \frac{ \bb{u}_{N_{f_i}} - \bb{u}_P}{\left|\bb{d}_{f_i}\right|}	- \bb{\Delta}_{f_i} \cdot \left(\bb{\nabla} \bb{u} \right)_{f_i}
		\right]    \left|\bb{\Gamma}_{f_i}\right| \notag \\
	&&+ \sum_{d_i \in \mathcal{F}^{\text{disp}}_P} \alpha \bar{K}_{d_i} \left[
		\left|\bb{\Delta}_{d_i} \right| \frac{ \bar{\bb{u}}_{d_i} - \bb{u}_P}{\left|\bb{d}_{d_i}\right|}	- \bb{\Delta}_{d_i} \cdot \left(\bb{\nabla} \bb{u} \right)_{P}
		\right]    \left|\bb{\Gamma}_{d_i}\right| \notag \\
	&&+ \sum_{s_i \in \mathcal{F}^{\text{symm}}_P} \alpha \bar{K}_{s_i} \left[
		\left|\bb{\Delta}_{s_i} \right| \frac{ \bb{R}_{s_i} \cdot \bb{u}_{P} - \bb{u}_P}{\left|\bb{d}_{s_i}\right|} - \bb{\Delta}_{s_i} \cdot \left(\bb{\nabla} \bb{u} \right)_{s_i}
		\right]    \left|\bb{\Gamma}_{s_i}\right|
\end{eqnarray}
where $\alpha > 0$ is a user-defined parameter for globally scaling the amount of stabilisation.
Parameter $\bar{K}$ is a stiffness-type parameter that gives the stabilisation term an appropriate scale and dimension.
Here, $\bar{K} = \frac{4}{3}\mu + \kappa = 2\mu + \lambda$ following previous work \cite{Jasak2000, Cardiff2017, Cardiff2018}, where $\mu$ is the shear modulus (first Lam\'{e} parameter), $\kappa$ is the bulk modulus, and $\lambda$ is the second Lam\'{e} parameter.
Vector $\bar{\bb{u}}_{d_i}$ represents the prescribed displacement at the centroid of displacement boundary face $d_i$.
The quantities $\bb{\Delta} = \nicefrac{\bb{d}}{\bb{d} \cdot \bb{n}}$ are termed the \emph{over-relaxed orthogonal} vectors \cite{Jasak1996}.
The displacement gradient $\left(\bb{\nabla} \bb{u} \right)_{f_i}$ at the internal face $f_i$ is calculated by interpolation from adjacent cell centres (like in Equation \ref{eq:stressInterp}).
Similarly, the displacement gradient $\left(\bb{\nabla} \bb{u} \right)_{s_i}$ at a symmetry boundary face $s_i$ is averaged from cell $P$ and its mirror reflection across face $s_i$, as in Equation \ref{eq:symm}.
Note that the displacement gradient at the displacement boundary $d_i$ is assumed to be equal to the displacement gradient at the centroid of cell $P$.
In addition, no stabilisation term is applied on a traction boundary face $t_i$.
The form of Rhie-Chow stabilisation given in Equation \ref{eq:RhieChow} can also be equivalently expressed in the form of a \emph{jump} term, for example, as in \citet{Nishikawa2010}: this is shown in Appendix \ref{app:RhieChowJump}.

\section{Solution Algorithms}\label{sec:sol_alg}


\subsection{Segregated Solution Algorithm} 
\label{sec:seg_alg}
The classic segregated solution algorithm can be viewed as a quasi-Newton method, where an approximate Jacobian is derived from the inertia term and a compact-stencil discretisation of a diffusion term (from the stabilisation term):
\begin{eqnarray} \label{eq:diffusion}
	\tilde{\bb{J}} &=& \frac{\partial}{\partial \bb{u}} \left[ \oint_{\Gamma_P} \alpha \bar{K} \, \bb{n} \cdot \bb{\nabla} \bb{u} \; d\Gamma_P
	 \; -\;  \int_{\Omega_P} \rho \frac{\partial^2 \bb{u} }{\partial t^2} \, d\Omega_P \right]
\end{eqnarray}
The inertia term is discretised as described in Equations \ref{eq:inertia} and \ref{eq:inertia2}, while the diffusion term is discretised in the same manner as the compact-stencil component of the Rhie-Chow stabilisation term:
\begin{eqnarray}
	\oint_{\Gamma_P} \alpha \bar{K} \, \bb{n} \cdot \bb{\nabla} \bb{u} \; d\Gamma_P &\approx&
		\sum_{f_i \in \mathcal{F}^{\text{int}}_P}  \alpha \bar{K}
		\left|\bb{\Delta}_{f_i} \right| \frac{ \bb{u}_{N_{f_i}} - \bb{u}_P}{\left|\bb{d}_{f_i}\right|}    \left|\bb{\Gamma}_{f_i}\right| \notag \\
	&&+  \sum_{d_i \in \mathcal{F}^{\text{disp}}_P}  \alpha \bar{K}
		\left|\bb{\Delta}_{d_i} \right| \frac{ \bar{\bb{u}}_{d_i}  - \bb{u}_P}{\left|\bb{d}_{d_i}\right|} 
		\left|\bb{\Gamma}_{d_i}\right| \notag \\
	&&+ \sum_{s_i \in \mathcal{F}^{\text{symm}}_P}  \alpha \bar{K}
		\left|\bb{\Delta}_{s_i} \right| \frac{ \bb{R}_{s_i} \cdot \bb{u}_{P} - \bb{u}_P}{\left|\bb{d}_{s_i}\right|}
		\left|\bb{\Gamma}_{s_i}\right|
\end{eqnarray}
When a diffusion term is typically discretised using the cell-centre finite volume method, non-orthogonal corrections are included in a deferred correction manner to preserve the order of accuracy on distorted grids.
However, in the current quasi-Newton method, the exact value of the approximate Jacobian does not affect the final converged solution, but only the convergence behaviour.
Consequently, non-orthogonal corrections are not included in the approximate Jacobian here. However, grid distortion is appropriately accounted for in the residual calculation.
Nonetheless, it is expected that the convergence behaviour of the segregated approach may degrade as mesh non-orthogonality increases.

The linearised system (Equation \ref{eq:Seg}) is formed for each cell in the domain, resulting in a system of algebraic equations:
\begin{eqnarray} \label{eq:SegSys}
    \left[ \bb{\tilde{J}} \right]  \; \left[ \delta \bb{u} \right] = - \left[\bb{R}(\bb{u}_k)\right]
\end{eqnarray}
where $\left[ \bb{\tilde{J}} \right]$ is a symmetric, $M \times M$ stiffness matrix, where $M = 3|\mathcal{P}|$ in 3-D and $M = 2|\mathcal{P}|$ in 2-D.
If $\Delta t < \infty$ or $\mathcal{F}^{\text{disp}}_P \neq \emptyset$, matrix $\left[ \bb{\tilde{J}} \right]$ is strongly diagonally dominant; otherwise, it is weakly diagonally dominant.
The block ($3\times3$ for 3-D, $2\times2$ for 2-D) diagonal coefficient for cell $P$ (row $P$, column $P$) can be expressed as
\begin{eqnarray}
	 \left[ \tilde{\bb{J}}\right]_{PP} &=&
		- \sum_{f_i \in \mathcal{F}^{\text{int}}_P}  \alpha \bar{K}
		\frac{ \left|\bb{\Delta}_{f_i} \right| }{\left|\bb{d}_{f_i}\right|}    \left|\bb{\Gamma}_{f_i}\right| \textbf{I} 
	    \quad-  \sum_{d_i \in \mathcal{F}^{\text{disp}}_P}  \alpha \bar{K}
		 \frac{ \left|\bb{\Delta}_{d_i} \right| }{\left|\bb{d}_{d_i}\right|} 
		\left|\bb{\Gamma}_{d_i}\right| \textbf{I} \notag \\
	 &&\quad - \sum_{s_i \in \mathcal{F}^{\text{symm}}_P} \alpha \bar{K}
		 \frac{ \left|\bb{\Delta}_{s_i} \right|}{\left|\bb{d}_{s_i}\right|}
		\left|\bb{\Gamma}_{s_i}\right|  \bb{n}_{s_i} \otimes \bb{n}_{s_i} 
	\quad - \quad \frac{9}{4}  \frac{\rho_P \Omega_P}{\Delta t^2} \textbf{I}
\end{eqnarray}
while the off-diagonal coefficients (row $P$, column $Q$) can be expressed as
\begin{eqnarray}
	\left[\tilde{\bb{J}}\right] _{PQ} &=&
		\alpha \bar{K} \frac{ \left|\bb{\Delta}_{f_{PQ}} \right| }{\left|\bb{d}_{f_{PQ}}\right|} \left|\bb{\Gamma}_{f_{PQ}}\right| \textbf{I} 
\end{eqnarray}
where subscript $PQ$ indicates a quantity associated with the internal face $f$ shared between cells $P$ and $Q$.

By design, matrix $\left[ \bb{\tilde{J}} \right]$  contains no inter-component coupling, i.e., each block coefficient is diagonal; consequently, three equivalent smaller linear systems can be formed in 3-D (or two linear systems in 2-D) and solved for the Cartesian components of the displacement correction, e.g.
\begin{eqnarray} \label{eq:SegSysX}
     \left[ \bb{\tilde{J}}_x \right]  \;  \left[ \Delta \bb{u}_x \right] = - \left[ \mathcal{R}_x(\bb{u}_k) \right] \label{eq:segX} \\
     \left[ \bb{\tilde{J}}_y \right]  \;  \left[ \Delta \bb{u}_y \right] = - \left[ \mathcal{R}_y(\bb{u}_k) \right] \label{eq:segY} \\
     \left[ \bb{\tilde{J}}_z \right]  \;  \left[ \Delta \bb{u}_z \right] = - \left[ \mathcal{R}_z(\bb{u}_k) \right] \label{eq:segZ}
\end{eqnarray}
where $ \bullet_x$ represents the components in the $x$ direction, $ \bullet_y$ represents the components in the $y$ direction, and $ \bullet_z$ represents the components in the $z$ direction.
Matrices $\left[ \bb{\tilde{J}}_x \right]$, $ \left[ \bb{\tilde{J}}_y\right]$ and $\left[ \bb{\tilde{J}}_z\right]$ have the size $|\mathcal{P}| \times |\mathcal{P}|$.
An additional benefit from a memory and assembly perspective, is that matrices $\left[ \bb{\tilde{J}}_x \right]$, $ \left[ \bb{\tilde{J}}_y\right]$ and $\left[ \bb{\tilde{J}}_z\right]$ are identical, except for the effects from including boundary conditions ($\mathcal{F}^{\text{symm}}_P$ terms).
From an implementation perspective, this allows a single scalar matrix to be formed and stored, where the boundary condition contributions are inserted before solving a particular component.

The \emph{inner} linear sparse systems (Equations \ref{eq:segX}, \ref{eq:segY} and \ref{eq:segZ}) can be solved using any typical direct or iterative linear solver approach; however, an incomplete Cholesky pre-conditioned conjugate gradient method \cite{Jacobs1986} is often preferred as the diagonally dominant characteristic leads to good convergence characteristics.
Algebraic multigrid can also accelerate convergence.
%

In the literature, the segregated solution algorithm is typically formulated with the total displacement vector (or its difference between time steps) as the primary unknown; in contrast, in the quasi-Newton interpretation presented here, the primary unknown is the correction to the displacement vector, which converges to zero.
In the total displacement approach, the matrix remains the same, and the only difference is the formulation of the right-hand side, which is computed as the residual minus the matrix times the previous solution.
Equivalence of both approaches can be seen by adding $\left[ \bb{\tilde{J}} \right]  \; \left[ \bb{u} \right]_k$ to both sides of Equation \ref{eq:SegSys}:
\begin{eqnarray} \label{eq:SegSys}
    \left[ \bb{\tilde{J}} \right]  \; \left[ \delta \bb{u} \right] + \left[ \bb{\tilde{J}} \right]  \; \left[ \bb{u} \right]_k  = - \left[\bb{R}(\bb{u}_k)\right] + \left[ \bb{\tilde{J}} \right]  \; \left[ \bb{u} \right]_k \notag \\
    \left[ \bb{\tilde{J}} \right]  \; \left[ \bb{u} \right]_{k+1} = \left[ \bb{\tilde{J}} \right]  \; \left[ \bb{u} \right]_k - \left[\bb{R}(\bb{u}_k)\right] \
\end{eqnarray}
where $\left[ \bb{u} \right]_k$ is the solution field at the previous iteration, and noting that $\left[ \bb{u} \right]_{k+1} = \left[ \bb{u} \right]_k + \left[ \delta \bb{u} \right]$.
In practice, the $\left[ \bb{\tilde{J}} \right]  \; \left[ \bb{u} \right]_k$ term on the right-hand side can be evaluated directly in terms of the discretised temporal and diffusion terms, without the need for matrix multiplication.
In this work, the standard total-displacement segregated formulation is employed.

As a final comment, the $\alpha$ scaling factor in the stabilisation term (Equation \ref{eq:RhieChow}) need not take the same value as in the approximate Jacobian (Equation \ref{eq:diffusion}).
In the current work, $\alpha$ is taken as unity in the approximate Jacobian, while the variation of its value in the residual calculation is examined in Section \ref{sec:RhieChowResults}.


\subsection{Jacobian-free Newton-Krylov Algorithm}
\label{sec:JFNK_alg}

As noted in the introduction, the Jacobian-free Newton-Krylov method avoids the need to explicitly construct the Jacobian matrix by approximating its action on a solution vector using the finite difference method, repeated here:
\begin{eqnarray} \label{eq:JF}
	\bb{J} \bb{v} \approx \frac{\bb{R}(\bb{u} + \epsilon \bb{v}) - \bb{R}(\bb{u})}{\epsilon}
\end{eqnarray}

The derivation of this approximation can be shown for a $2 \times 2$ system as \cite{Knoll2004}:
\begin{eqnarray}
	\frac{\bb{R}(\mathbf{u} + \epsilon \mathbf{v}) - \mathbf{R}(\mathbf{u})}{\epsilon}
	&=&
	\begin{pmatrix}
	\dfrac{R_1 (u_1 + \epsilon v_1, u_2 + \epsilon v_2) - R_1 (u_1, u_2)}{\epsilon}\\
	\dfrac{R_2 (u_1 + \epsilon v_1, u_2 + \epsilon v_2) - R_2 (u_1, u_2)}{\epsilon}
	\end{pmatrix} \notag \\
	&\approx&
	\begin{pmatrix}
	\dfrac{R_1 (u_1,u_2) + \epsilon v_1 \dfrac{\partial R_1}{\partial u_1} + \epsilon v_2 \dfrac{\partial R_1}{\partial u_2} - R_1 (u_1, u_2)}{\epsilon}\\
	\dfrac{R_2 (u_1, u_2) + \epsilon v_1 \dfrac{\partial R_2}{\partial u_1} + \epsilon v_2 \dfrac{\partial R_2}{\partial u_2}  - R_2 (u_1, u_2)}{\epsilon}
	\end{pmatrix} \notag \\
	&\approx&
	\begin{pmatrix}
	v_1 \dfrac{\partial R_1}{\partial u_1} +  v_2 \dfrac{\partial R_1}{\partial u_2} \\
	v_1 \dfrac{\partial R_2}{\partial u_1} + v_2 \dfrac{\partial R_2}{\partial u_2}
	\end{pmatrix} \notag \\
	&\approx&
	\bb{J} \bb{v}
\end{eqnarray}
where a first-order truncated Taylor series expansion about $\bb{u}$ was used to approximate $\bb{R} (\bb{u} + \epsilon \bb{v})$.
As noted above, choosing an appropriate value for $\epsilon$ is non-trivial, and care must be taken to balance truncation error (reduced by decreasing $\epsilon$) and round-off error (increased by decreasing $\epsilon$).


The purpose of preconditioning the Jacobian-free Newton-Krylov method is to reduce the number of inner linear solver iterations.
The current work uses the GMRES linear solver for the inner system.
Using right preconditioning, the finite difference approximation of Equation \ref{eq:JF} becomes
\begin{eqnarray}
	\bb{J} \bb{P}^{-1} \bb{v}
	\approx
	\frac{\bb{R}(\bb{u} + \epsilon \bb{P}^{-1} \bb{v}) - \bb{R}(\bb{u})}{\epsilon}
\end{eqnarray}
where $\bb{P}$ is the preconditioning matrix or process.
In practice, only the action of $\bb{P}^{-1}$ on a vector is required, and $\bb{P}^{-1}$ need not be explicitly formed.
Concretely, the preconditioner needs to approximately solve the linear system $\bb{y} = \bb{P}^{-1} \bb{v}$.
Alternatively, left preconditioning may be applied. In this case, the finite difference approximation of Equation \ref{eq:JF} becomes
\begin{eqnarray}
	\bb{P}^{-1} \bb{J} \bb{v} \;\approx\; \frac{\bb{P}^{-1}\left[ \bb{R}(\bb{u}+\epsilon \bb{v}) - \bb{R}(\bb{u}) \right]}{\epsilon}
\end{eqnarray}
where the preconditioner acts directly on the residual evaluations. As with right preconditioning, the explicit inverse of $\bb{P}$ is not required; only the ability to approximately apply $\bb{P}^{-1}$ to a given vector is necessary.
In practice, the choice between left- and right-preconditioning in Jacobian-free Newton–Krylov methods is not always straightforward. Right preconditioning preserves the true nonlinear residual for convergence checks, while left preconditioning alters the residual norm through the preconditioner. In this work, we use left preconditioning, but no significant difference in convergence behaviour was observed between the two approaches.

In the current work, we propose to use the compact-stencil approximate Jacobian from the segregated algorithm $\tilde{\bb{J}}$ as the preconditioning matrix $\bb{P}$ for the preconditioned Jacobian-free Newton-Krylov method.
This preconditioning approach can be considered a ``physics-based" preconditioner in the classifications of \citet{Knoll2004}.
The approach is conceptually similar to approximating the Jacobian of a higher-order, large-stencil scheme by a lower-order, compact-stencil scheme.
A benefit of the proposed approach is that existing segregated frameworks can reuse their existing matrix assembly and storage implementations.
Concretely, the Jacobian-free Newton-Krylov method requires only a procedure for forming this preconditioning matrix and a procedure for explicitly evaluating the residual.
Both routines are easily implemented - and are likely already available - in an existing segregated framework.
The only additional required procedure is an interface to an existing Jacobian-free Newton-Krylov implementation.
In the current work, the PETSc package (version 3.22.2) \cite{PETSc} is used as the nonlinear solver, driven by a finite volume solver implemented in the solids4foam toolbox \citep{Cardiff2018, Tukovic2018} for the OpenFOAM toolbox \citep{Weller1998} (version OpenFOAM-v2312).
The codes are publicly available at \url{https://github.com/solids4foam/solids4foam} on the \texttt{feature-petsc-snes} branch, and the cases and plotting scripts are available at \url{https://github.com/solids4foam/solid-benchmarks}.

Several preconditioners are available in the literature, with incomplete Cholesky and ILU($k$) being popular; however, multigrid methods offer the greatest potential for large-scale problems.
As noted by \citet{Knoll2004}, algorithmic simplifications within a multigrid procedure, which may result in loss of convergence for multigrid as a solver, have a much weaker effect when multigrid is used as a preconditioner.
In this work, three preconditioners are considered:
\begin{enumerate}
	\item ILU($k$): incomplete lower-upper decomposition with fill-in $k$. 
	\item Multigrid: the HYPRE Boomerang \citep{hypre} multigrid implementation.
	\item LU: the MUMPS \citep{MUMPS:1, MUMPS:2} lower-upper decomposition direct solver.
\end{enumerate}

A challenge with Newton-type methods, including Jacobian-free versions, is poor convergence when far from the true solution, and divergence is often a real possibility.
Globalisation refers to steering an initial solution towards the quadratic convergence range of the Newton method.
Several strategies are possible, and it is common to combine approaches \cite{Knoll2004}.
In the current work, a line search procedure is used to select the $s$ parameter in the solution update step (the second line in Equations \ref{eq:NewtonRaphson}).
Line search methods assume the Newton update direction is correct and aim to find a scalar $s > 0$ that decreases the residual $\bb{R}(\bb{u}_k + s \delta \bb{u}) < \bb{R}(\bb{u}_k)$.
In addition to a line search approach, a \emph{transient continuation} globalisation approach is used in the current work, where the displacement $\bb{u}_P$ for cell $P$ at time $t + \Delta t$ is predicted at the start of a new time (or loading) step, based on a truncated second-order Taylor series expansion:
\begin{eqnarray} \label{eq:predictor}
	\bb{u}^{[t+\Delta t]}_P = \bb{u}^{[t]}_P + \Delta t \, \bb{v}^{[t]}_P + \frac{1}{2} \Delta t^2 \left( \frac{\partial \bb{v}}{\partial t} \right)^{[t]}_P
\end{eqnarray}
where $\bb{v}^{[t]}_P$ is the velocity of cell $P$ at time $t$, and $\left( \frac{\partial \bb{v}}{\partial t} \right)^{[t]}_P$ is the acceleration.
In this way, for highly nonlinear problems, the user can decrease the time step size $\Delta t$ as a globalisation approach to improve the Newton method's performance.
The predictor step in Equation \ref{eq:predictor} has been chosen to be consistent with the discretisation of the inertia term (Equation \ref{eq:inertia2}).

A final comment on the Jacobian-free Newton-Krylov solution algorithm is the potential importance of \emph{oversolving}.
Here, oversolving refers to solving the linear system to a tolerance that is too tight during the early Newton iterations, essentially wasting time when the solution is far from the true solution.
In addition, some authors 
\cite{Knoll2004} has shown that Newton convergence is worse when earlier iterations are solved to too tight a tolerance.
The concept of oversolving also applies to segregated solution procedures and has been well-known since the early work of Demird\v{z}i\'{c} and co-workers \cite{Demirdzic1995}, where the residuals are typically reduced by one order of magnitude in the inner linear system.
In the current work, the residual is reduced by a factor of 0.9 in the segregated approach and by three orders in the Jacobian-free Newton-Krylov approach.


\subsection{Checking Convergence}
The current work adopts three checks for determining whether convergence has been achieved within each time (or loading) step, closely following the default strategy in the PETSc toolbox:
\begin{itemize}
    \item Norm of the residual $|\bb{R}(\bb{u})|$:
    \begin{itemize}
        \item \emph{Absolute}: Convergence is declared when $|\bb{R}(\bb{u})|$ falls below an absolute tolerance $a_{\text{tol}}$, taken here as $1\times 10^{-50}$. This acts as a failsafe in situations where the \emph{relative} criterion (below) is ineffective because the initial solution already nearly satisfies the governing equation.
        \item \emph{Relative}: Convergence is declared when $|\bb{R}(\bb{u})|$ falls below $r_{\text{tol}} \times |\bb{R}(\bb{u}_0)|$, with $|\bb{R}(\bb{u}_0)|$ denoting the residual norm at the start of the time (or loading) step. In this work, $r_{\text{tol}} = 10^{-6}$.
    \end{itemize}
    \item Norm of the solution correction $|\delta \bb{u}|$: Convergence is declared when the norm of the change in the solution between successive iterations is less than $s_{\text{tol}} \times |\bb{u}|$.
\end{itemize}

Although applying the same tolerances for the segregated and Jacobian-free Newton-Krylov solution procedures may seem reasonable, care must be taken with $s_{\text{tol}}$. The Jacobian-free Newton-Krylov method typically makes a small number of large corrections to the solution, whereas the segregated procedure often makes a large number of small corrections. Consequently, if the same $s_{\text{tol}}$ value were used, the segregated procedure might be declared convergent prematurely. To avoid this, $s_{\text{tol}}$ is set to zero in the present work, meaning convergence is based solely on residual reduction. Nonetheless, setting $s_{\text{tol}} > 0$ could yield a more robust procedure that can handle potential stalling.

\section{Test Cases}\label{sec:test_cases}

This section compares the performance of the proposed Jacobian-free Newton-Krylov solution approach with that of the segregated approach across several benchmark cases.
The cases have been chosen to exhibit a variety of characteristics in terms of
\begin{itemize}
	\item Geometric dimension (2-D vs. 3-D),
	\item Geometric nonlinearity (small strain vs. large strain),
	\item Geometric complexity (basic geometric shapes vs. complex geometry),
	\item Statics vs. dynamics, and
	\item Material behaviour (elasticity, elastoplasticity, hyperelasticity).
\end{itemize}

Although the presented analyses aim to be extensive, several common features of modern solid mechanics procedures are left for future work, including nonlinear boundary conditions (e.g., contact, fracture) and mixed formulations (e.g., incompressibility). 

It is noted that the benchmark cases considered in this study are modest in size compared with the largest-scale finite element studies. This choice was deliberate, as the aim was to benchmark the Jacobian-free Newton–Krylov method against the established segregated finite-volume approach under controlled, reproducible conditions, while keeping memory and runtime requirements practical. Application to much larger problems is feasible and will be the subject of future work.

The remainder of this section is structured as follows:
The proposed discretisation's accuracy and order are assessed on several test cases of varying dimensions and phenomena.
Subsequently, the efficiency of the Jacobian-free Newton-Krylov approach is compared with the standard segregated procedure in terms of time, memory and robustness. Finally, the effects of the preconditioning strategy and stabilisation scaling are examined.

\subsection{Testing the Accuracy and Order of Accuracy}
\label{sec:accuracy}
This section focuses solely on assessing the accuracy and the order of accuracy of the discretisation.
Assessment of computational efficiency in terms of time and memory requirements is left to the subsequent sections.
In all cases where convergence is achieved, the differences between the Jacobian-free Newton-Raphson and segregated approaches are minimal; this is expected, since both use the same discretisation and the solution tolerances have been chosen to ensure the iteration errors are small.
Consequently, the results presented below (Section \ref{sec:accuracy}) have been solely generated using the Jacobian-free Newton-Raphson solution algorithm, unless stated otherwise.


 
\begin{figure}[htbp]
	\centering
	\subfigure[Magnitude of the manufactured displacement solution]
	{
		\label{fig:mms_solution}
   		\includegraphics[height=0.45\textwidth]{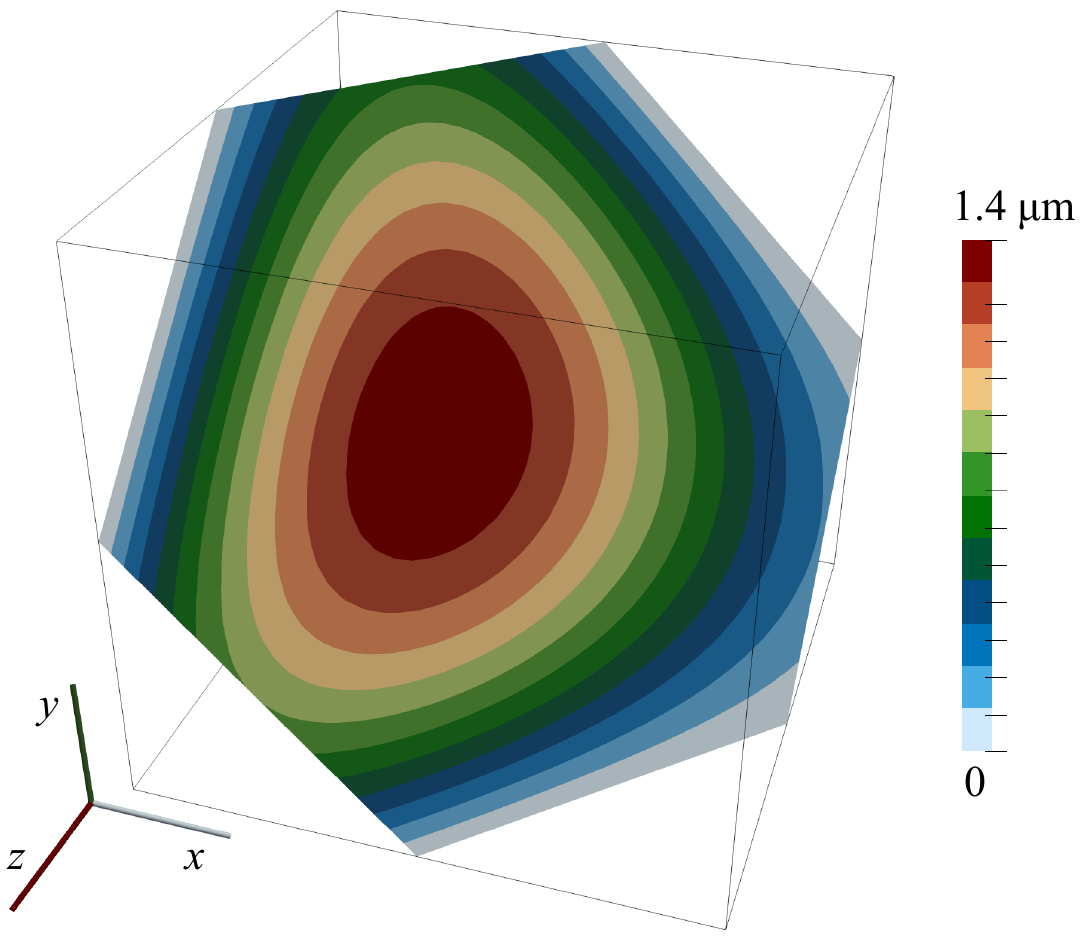} 
   	}
	\subfigure[Regular polyhedral mesh with $1\,000$ cells]
	{
		\label{fig:mms_mesh}
   		\includegraphics[height=0.45\textwidth]{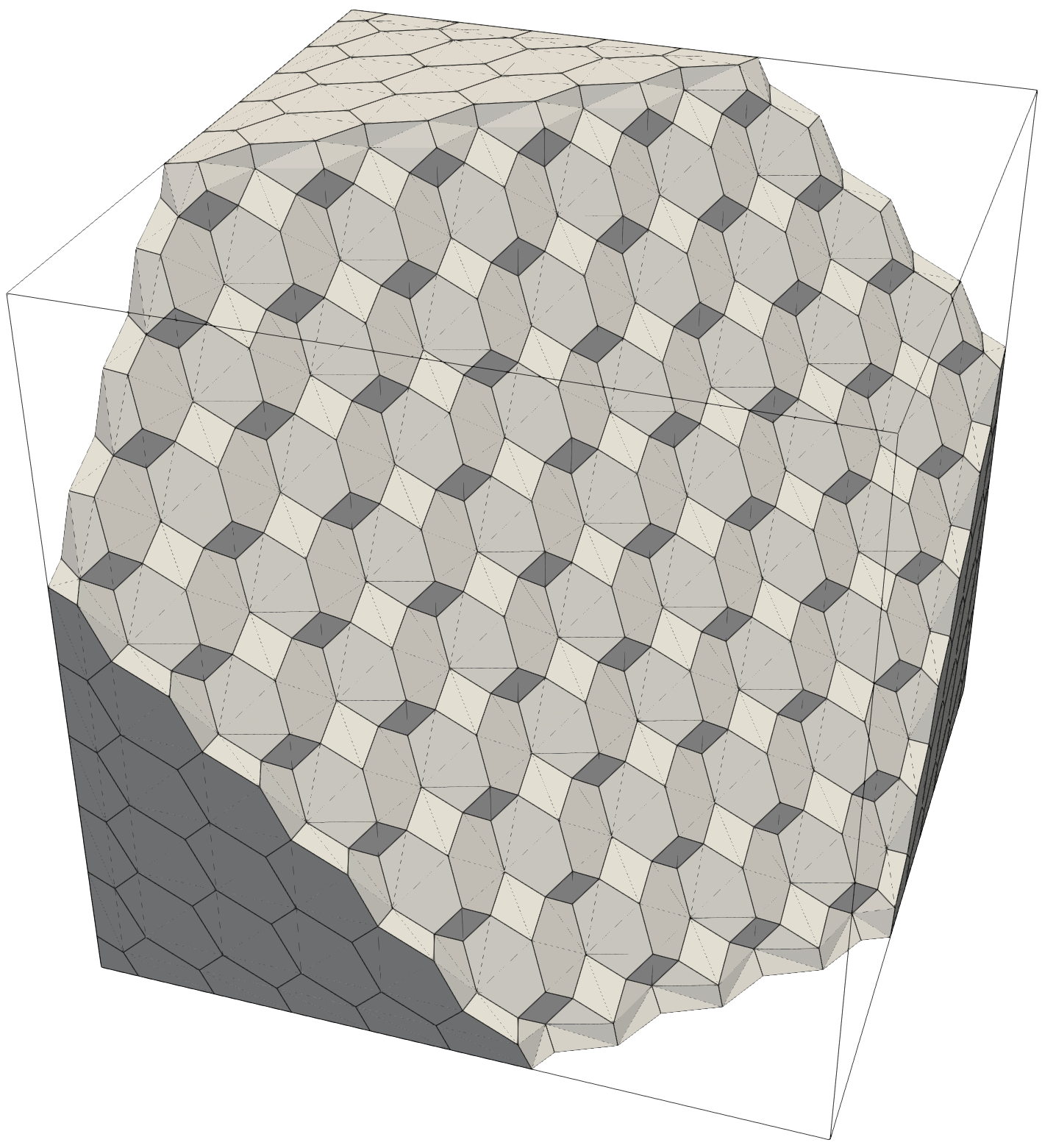}  
   	}
	\caption{A cut plane through the cube case geometry showing the magnitude of the manufactured displacement solution (left) and a polyhedral mesh (right). The cut plane passes through the centre of the cube and has the unit normal $\bb{n} = (\sfrac{1}{\sqrt{3}} \quad \sfrac{1}{\sqrt{3}} \quad \sfrac{1}{\sqrt{3}})$.}
	\label{fig:mms}
\end{figure}

\paragraph{Case 1: Order Verification via the Manufactured Solution Procedure}
The first test case consists of a $0.2 \times 0.2 \times 0.2$ m cube with linear elastic ($E = 200$ GPa, $\nu = 0.3$) properties.
A manufactured solution for displacement (Figure \ref{fig:mms_solution}) is employed of the form \citep{Mazzanti2024}
\begin{eqnarray}
	\bb{u} =
	\begin{pmatrix}
	a_x \sin(4\pi x) \sin(2\pi y) \sin(\pi z) \\
	a_y \sin(4 \pi x) \sin(2 \pi y) \sin(\pi z) \\
	a_z \sin(4 \pi x) \sin(2 \pi y) \sin(\pi z) 
	\end{pmatrix}
\end{eqnarray}
where $a_x = 2\,\mu$m, $a_y = 4\,\mu$m, and $a_z = 6\,\mu$m.
The Cartesian coordinates are $x$, $y$, and $z$.
The corresponding manufactured body force term ($\bb{f}_b$  in Equation \ref{eqn:momentum_lingeom}) is given in Appendix \ref{app:mms}.

The manufactured displacement solution is applied at the domain's boundaries (prescribed displacements), and inertial effects are neglected.
Four mesh types are examined:
(i) \emph{regular} tetrahedra, (ii) \emph{regular} polyhedra (shown in Figure \ref{fig:mms_mesh}), (iii) \emph{regular} hexahedra, and (iv) \emph{distorted} hexahedra.
The second coarsest meshes are shown in Appendix \ref{app:meshes}.
The tetrahedral meshes are created using Gmsh \citep{geuzaine2009gmsh}, while the polyhedral meshes are created by converting the tetrahedral meshes to their dual polyhedral representations using the OpenFOAM \texttt{polyDualMesh} utility.
The regular hexahedral meshes are created using the OpenFOAM \texttt{blockMesh} utility, and the distorted hexahedral meshes are created by perturbing the points of the regular hexahedra by a random vector of magnitude less than 0.3 times the local minimum edge length.
Starting from an initial mesh spacing of $0.04$ m, six meshes of each type are created by successively halving the spacing.
The cell numbers for the hexahedral and polyhedral meshes are 125, $1\,000$, $8\,000$, $64\,000$, $512\,000$, and $4\,096\,000$, while the tetrahedral mesh cell counts are 384, $4\,374$, $41\,154$, $355\,914$, $2\,958\,234$, and $24\,118\,074$.
For the same average cell width, the cell counts show that tetrahedral meshes have 3 to 6 times as many cells.

\begin{figure}[htbp]
	\centering
	\subfigure[Displacement magnitude discretisation errors]
	{
   		\includegraphics[width=0.45\textwidth]{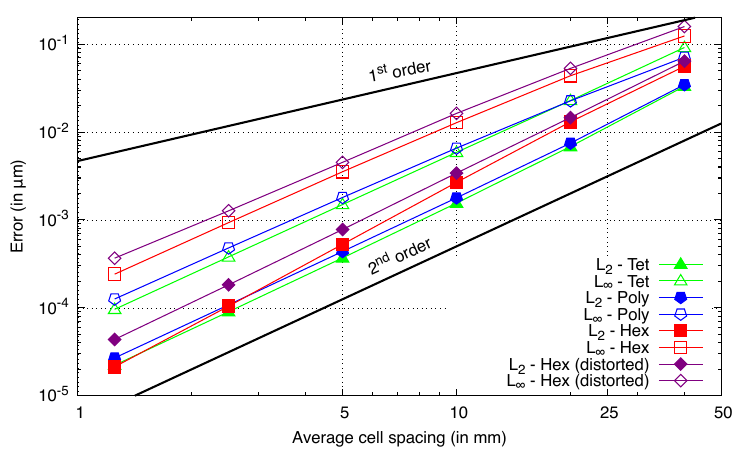} 
   	}
	\subfigure[Displacement discretisation error order of accuracy]
	{
   		\includegraphics[width=0.45\textwidth]{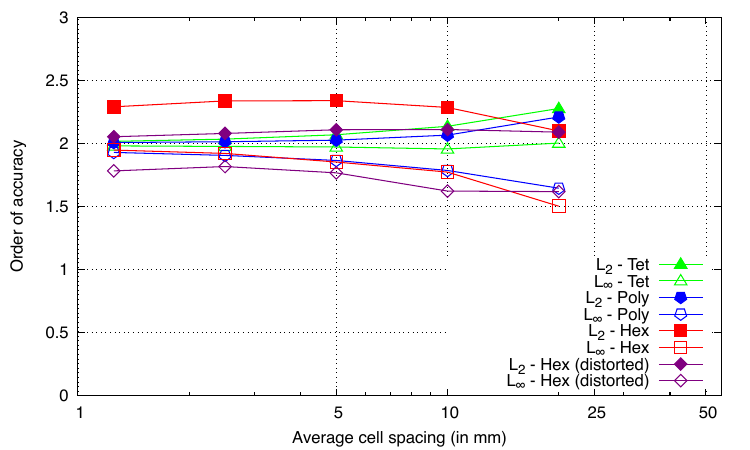}  
   	}
	\caption{Manufactured solution cube case: the accuracy and order of accuracy for displacement magnitude}
	\label{fig:mms_disp_accuracy}
\end{figure}

Figure \ref{fig:mms_disp_accuracy}(a) shows the displacement magnitude discretisation errors ($L_2$ and $L_\infty$) as a function of the average cell width for the four mesh types (tetrahedral, polyhedral, regular and distorted hexahedral), while Figure \ref{fig:mms_disp_accuracy}(b) shows the corresponding order of accuracy plots.
For ease of interpretation, the symbol shapes in the figures have been chosen to correspond to the cell shapes: a triangle for tetrahedra, a pentagon for polyhedra, a square for regular hexahedra, and a diamond for distorted hexahedra.
The maximum ($L_\infty$) and average ($L_2$) displacement discretisation errors are seen to reduce at an approximately second-order rate for all mesh types, except for the $L_2$ error on the hexahedral meshes, which is approximately 2.3 on the finest grid.
The distorted hexahedral meshes show the largest average and maximum errors, but still approach second-order accuracy as the element spacing decreases.

\begin{figure}[htbp]
	\centering
	\subfigure[Predicted $\sigma_{xx}$ stress distribution]
	{
		\label{fig:mms_sxx}
   		\includegraphics[height=0.35\textwidth]{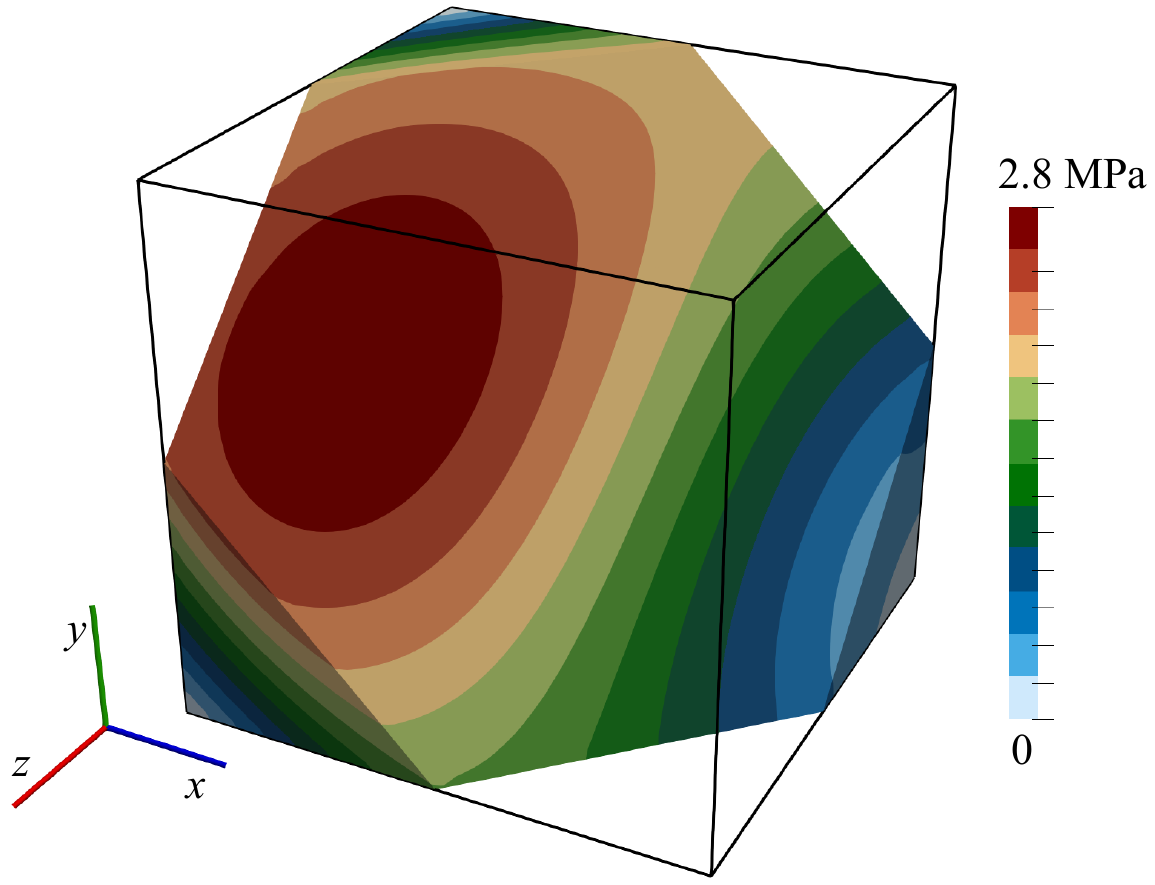} 
   	}
	\subfigure[Cell-wise $\sigma_{xx}$ error distribution]
	{
		\label{fig:mms_sxx_diff}
   		\includegraphics[height=0.35\textwidth]{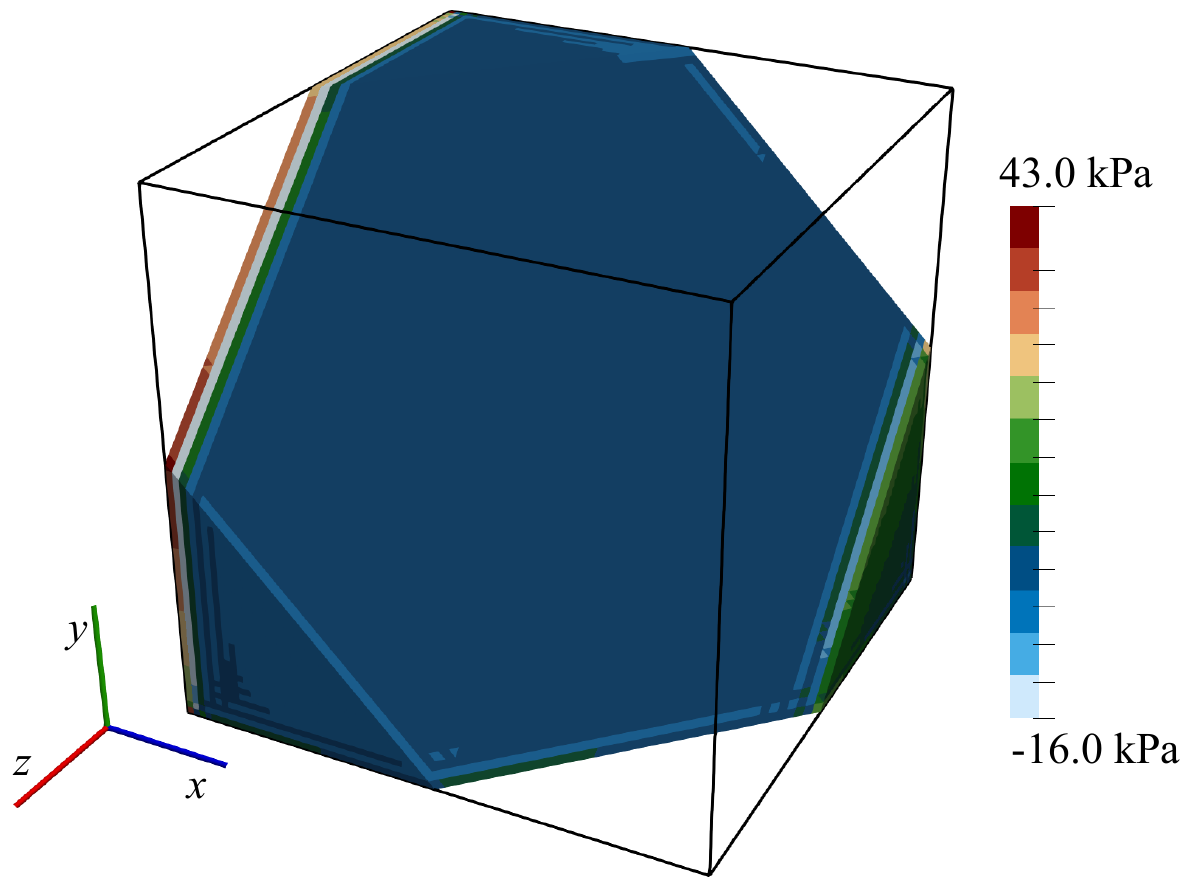}  
   	}
	\caption{Manufactured solution cube case: the predicted $\sigma_{xx}$ stress distribution on a cut-plane for the regular hexahedral mesh with $512\,000$ cells (left).}
	\label{fig:mms_stress}
\end{figure}

\begin{figure}[htbp]
	\centering
	\subfigure[Stress magnitude discretisation errors]
	{
   		\includegraphics[width=0.45\textwidth]{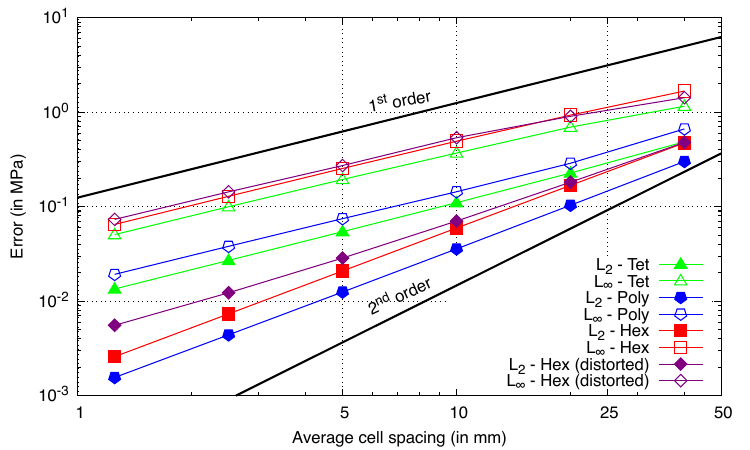} 
   	}
	\subfigure[Stress discretisation error order of accuracy]
	{
   		\includegraphics[width=0.45\textwidth]{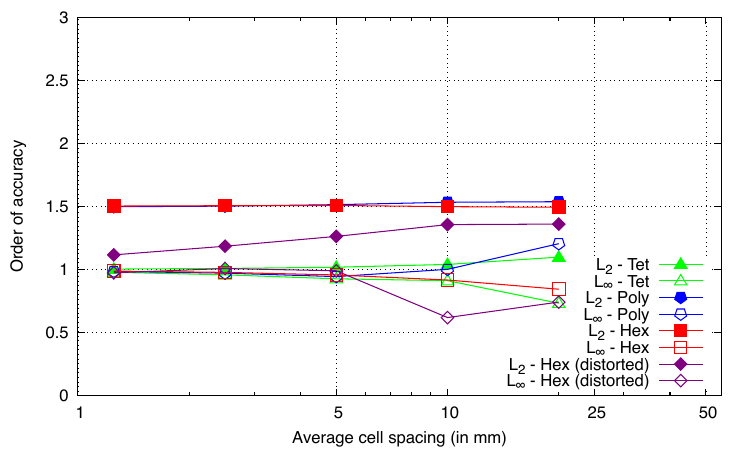}  
   	}
	\caption{Manufactured solution cube case: the accuracy and order of accuracy for stress magnitude}
	\label{fig:mms_stress_accuracy}
\end{figure}

The predicted $\sigma_{xx}$ stress distribution for a hexahedral mesh with $512\,000$ cells is shown in Figure \ref{fig:mms_stress}(a).
The corresponding cell-wise $\sigma_{xx}$ error distribution is shown in Figure \ref{fig:mms_stress}(b).
The largest errors ($43$ kPa) occur at the boundaries, where the local truncation error is highest. 
The discretisation errors in the stress magnitude as a function of cell size are shown in Figure \ref{fig:mms_stress_accuracy}(a), and the order of accuracy in Figure \ref{fig:mms_stress_accuracy}(b).
The order of accuracy for the average ($L_2$) stress error is approximately 1.5 for both hexahedral and polyhedral meshes.
In contrast, the average stress error order for the tetrahedral and distorted hexahedral meshes is 1.
Similarly, the maximum ($L_\infty$) stress order of accuracy is seen to approach 1 for all mesh types.
Unlike for the displacements, the greatest stress errors are seen in the tetrahedral rather than the distorted hexahedral meshes; however, the greatest stress maximum errors occur in the distorted hexahedral meshes.

The presented results agree with the theory: the least-squares gradient scheme should be at least first-order accurate for random grids, with higher-order accuracy on regular grids due to local truncation-error cancellations.
In addition, the displacement predictions should be at least second-order accurate on all grid types, provided the gradient scheme is at least first-order accurate.

\begin{figure}[htbp]
   \centering
	\subfigure[Polyhedral mesh with $4\,539$]
	{
	   \includegraphics[height=0.4\textwidth]{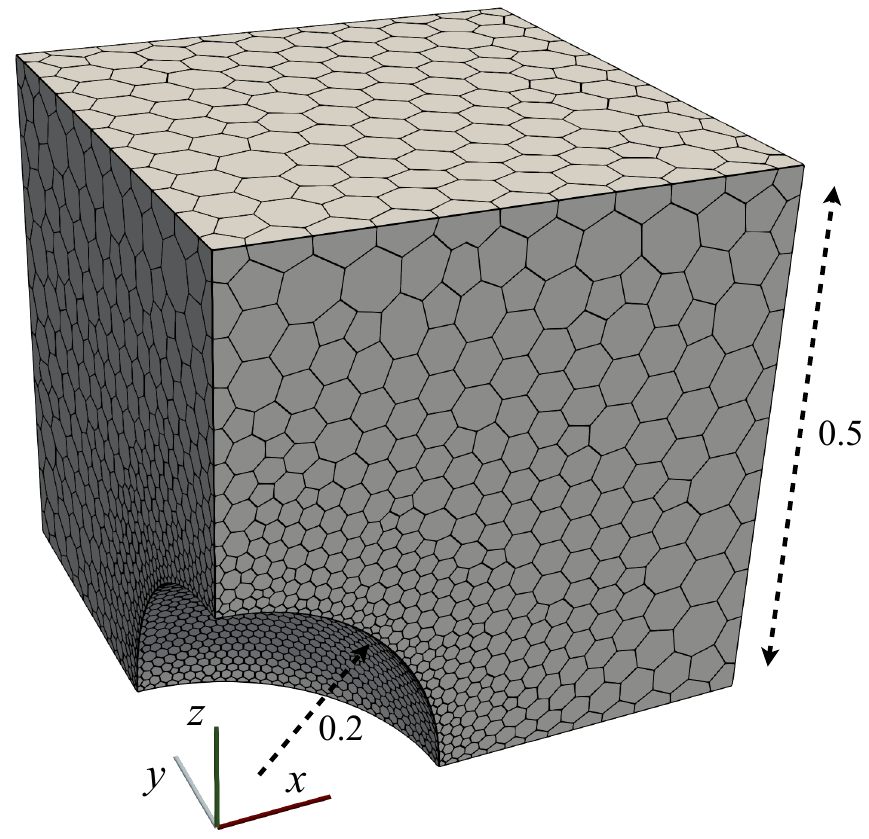} 
   	}
	\subfigure[Axial ($zz$) stress distribution on the mesh with $1\,614\,261$ cells]
	{
	   \includegraphics[height=0.4\textwidth]{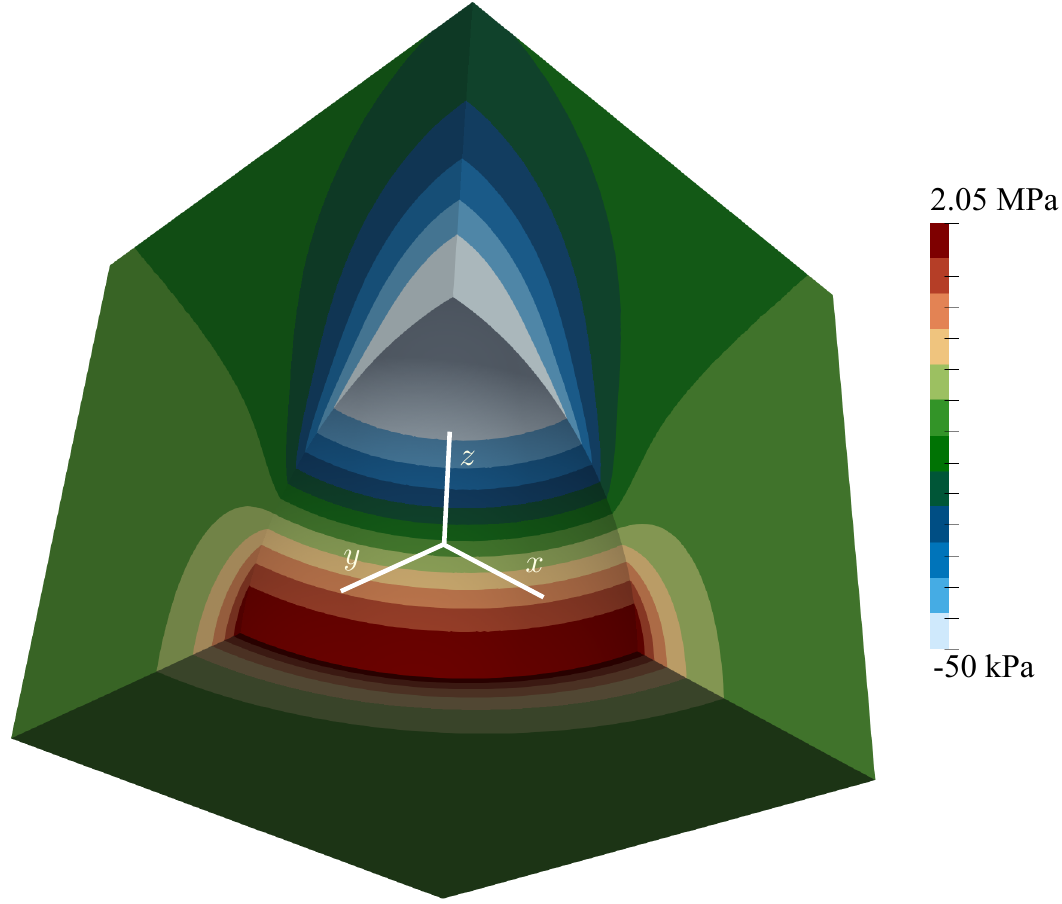} 
	   \label{fig:spherical_cavity_stress}
   	}
   \caption{Spherical cavity case}
   \label{fig:spherical_cavity}
\end{figure}

\paragraph{Case 2: Spherical Cavity in an Infinite Solid Subjected to Remote Stress}
This classic 3-D problem consists of a spherical cavity with radius $a = 0.2$ m (Figure \ref{fig:spherical_cavity}) in an infinite, isotropic linear elastic solid ($E = 200$ GPa, $\nu = 0.3$).
Far from the cavity, the solid is subjected to a tensile stress $\sigma_{zz} = T = 1$ MPa, with all other stress components zero.
The analytical expressions for the stress, derived by \citet{Southwell1926},  and displacement, derived by \citet{Goodier1933}, are given in Appendix \ref{app:sphericalCavity}.
The problem is characterised by a localised stress concentration near the cavity on the plane perpendicular to the loading, which drops off rapidly away from the cavity.

The computational solution domain is taken to be one-eighth of a $1 \times 1 \times 1$ m cube, aligned with the Cartesian axes, with one corner at the centre of the sphere.
The problem is axisymmetric but is analysed here using a 3-D domain to allow a graded unstructured polyhedral mesh to be assessed.
Analytical tractions are applied at the domain boundaries to mitigate the effects of the finite geometry.
The average cell sizes at the cavity surfaces are 100, 50, 25, 12.5, 6.25, and 3.125 mm, and the corresponding cell counts are 976, $4\,552$ (Figure \ref{fig:spherical_cavity}), $29\,611$, $213\,100$, $1\,614\,261$.
Initially, unstructured tetrahedral meshes were generated using the Gmsh meshing utility \cite{geuzaine2009gmsh}, followed by conversion to their dual polyhedral representations using the OpenFOAM \texttt{polyDualMesh} utility.

\begin{figure}[htbp]
	\centering
	\subfigure[Displacement magnitude discretisation errors]
	{
	   \includegraphics[width=0.45\textwidth]{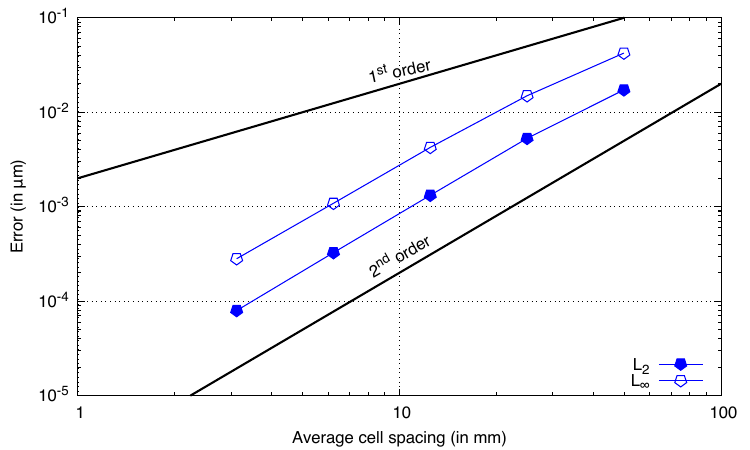} 
   	}
	\subfigure[Displacement discretisation error order of accuracy]
	{
	   \includegraphics[width=0.45\textwidth]{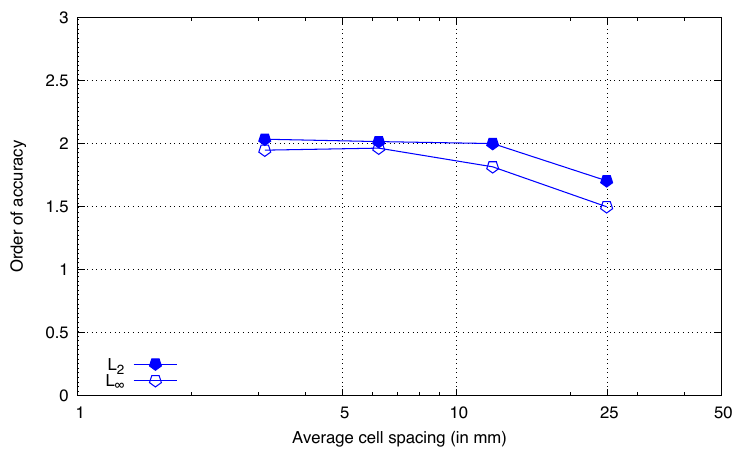} 
   	}
	\subfigure[Axial ($zz$) stress discretisation errors]
	{
	   \includegraphics[width=0.45\textwidth]{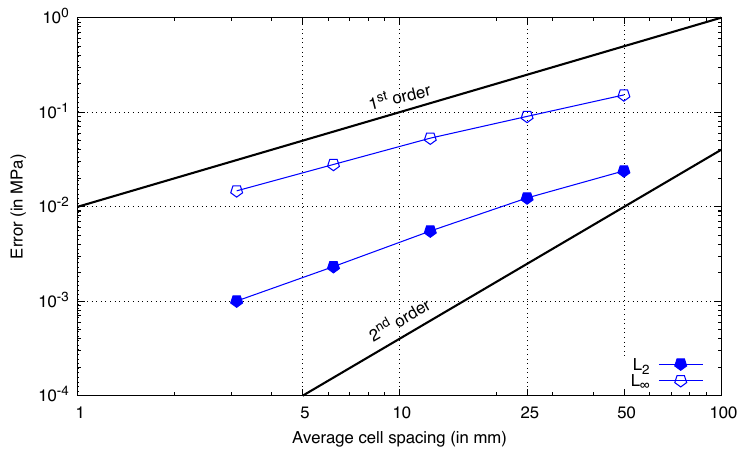} 
   	}
	\subfigure[Axial ($zz$) stress discretisation error order of accuracy]
	{
	   \includegraphics[width=0.45\textwidth]{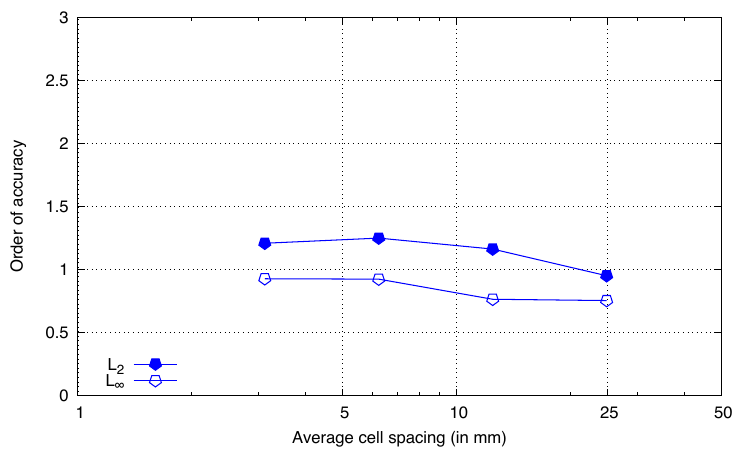} 
   	}
	\caption{Spherical cavity discretisation errors as a function of the average cell spacing}
	\label{fig:spherical_cavity_accuracy}
\end{figure}

Figure \ref{fig:spherical_cavity_stress} shows the predicted axial ($zz$) stress field on the mesh with $1\,614\,261$ cells.
The predicted mean ($L_2$) and maximum ($L_\infty$) displacement and axial stress discretisation errors as a function of the average cell width are presented in Figure \ref{fig:spherical_cavity_accuracy}.
The mean and maximum displacement errors are seen to decrease at approximately second-order rates, while the stress errors decrease at first-order rates.
As seen in the previous case, the results agree with the theoretical expectations for unstructured (irregular) meshes.

\begin{figure}[htbp]
   \centering
		\includegraphics[width=0.45\textwidth]{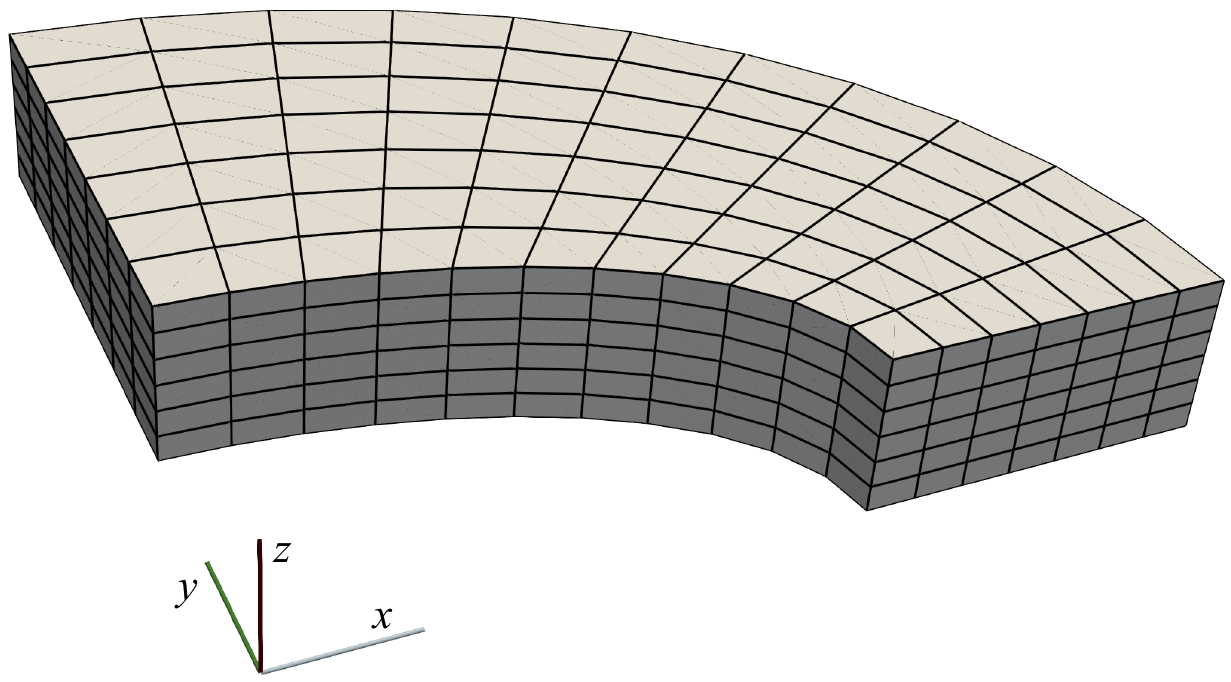} 
   \caption{Elliptic Plate geometry and mesh containing 472 hexahedral cells}
   \label{fig:elliptic_plate}
\end{figure}

\paragraph{Case 3: Out-of-plane bending of an elliptic plate}
This 3-D, static, linear elastic test case (Figure \ref{fig:elliptic_plate}) consists of a thick elliptic plate (0.6 m thick) with a centred elliptic hole, with the inner and outer ellipses given as
\begin{eqnarray}
	\left(\frac{x}{2}\right)^2 + \left(\frac{y}{1}\right)^2 = 1 & \text{inner ellipse} \\
	\left(\frac{x}{3.25}\right)^2 + \left(\frac{y}{2.75}\right)^2 = 1 & \text{outer ellipse}
\end{eqnarray}
The case has been described by the National Agency for Finite Element Methods and Standards (NAFEMS) \cite{Hitchings1987} and analysed using finite volume procedures by \citet{Demirdzic1997a} and \citet{Cardiff2016a}.
Unlike the previous cases, this case features significant bending, which is known to slow the convergence of segregated solution procedures \citep{Cardiff2016a}.
Symmetry allows one-quarter of the geometry to be simulated.
A constant pressure of 1 MPa is applied to the upper surface, and the outer surface is fully clamped.
The mechanical properties are: $E = 210$ GPa, $\nu = 0.3$.
Six successively refined hexahedral meshes are used, with cell counts of 45, 472 (Figure \ref{fig:elliptic_plate}), $4\,140$, $34\,968$, $287\,280$ and $2\,438\,242$.

\begin{figure}[htbp]
   \centering
	\subfigure[Equivalent stress distribution for the mesh with $2\,438\,242$ cells. The deformation is scaled by $1\,000$ and a translucent line indicates the outline of the undeformed geometry.]
	{
		\includegraphics[width=0.45\textwidth]{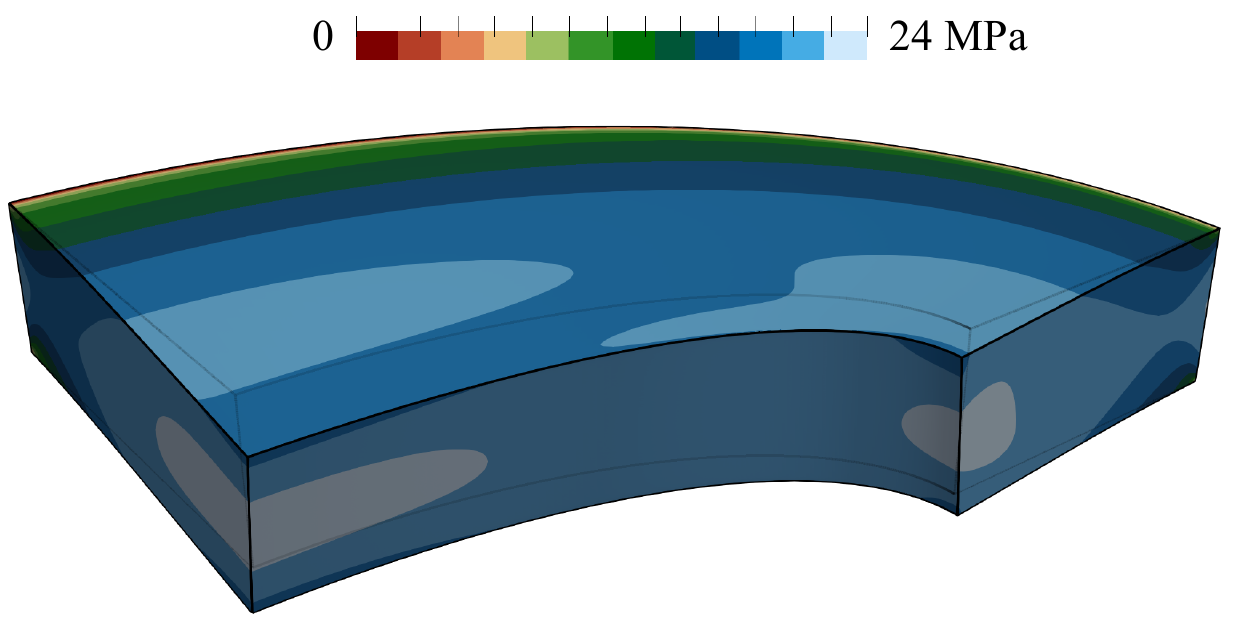}  
	}
	\subfigure[Equivalent stress along the line $r = 2.1$, $z = 0.3$ m]
	{
		\includegraphics[width=0.45\textwidth]{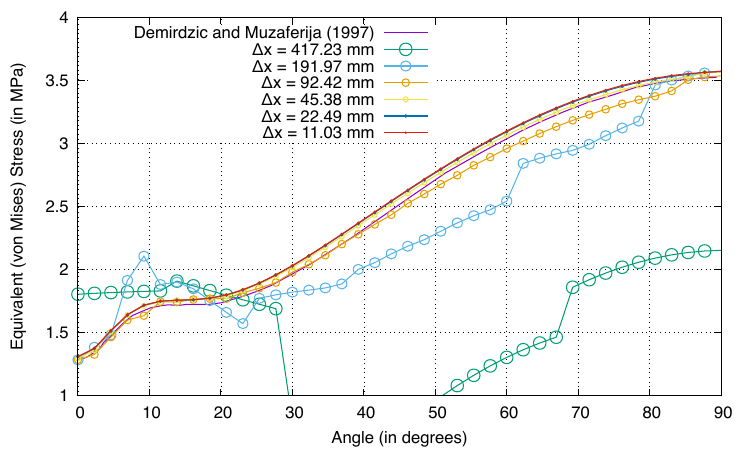} 
	}
   \caption{Equivalent (von Mises) stress distribution in the elliptic plate}
   \label{fig:elliptic_plate_sigmaEq}
\end{figure}

The prediction for the equivalent (von Mises) stress in the domain is shown in Figure \ref{fig:elliptic_plate_sigmaEq}(a), and the values along the line $r = \sqrt{x^2 + y^2} = 2.1$ m, $z = 0.3$ m are shown in Figure \ref{fig:elliptic_plate_sigmaEq}(b).
The results from the finest grid in \citet{Demirdzic1997a} are given for comparison, where good agreement is seen; the small offset between the finest mesh predictions and those from \citet{Demirdzic1997a} is likely due to errors introduced when extracting the \citet{Demirdzic1997a} results using WebPlotDigitizer \citep{WebPlotDigitizer}.
The sample values have been calculated by generating 40 uniformly spaced angles between 0 and $\nicefrac{\pi}{2}$ along the line, finding the cell of each sampled point, and extrapolating from the cell-centred value using the field gradient.

\begin{figure}[htbp]
   \centering
	\subfigure[Ventricle undeformed geometry, showing the mesh with 1,620 cells]
	{
		\includegraphics[height=0.46\textwidth]{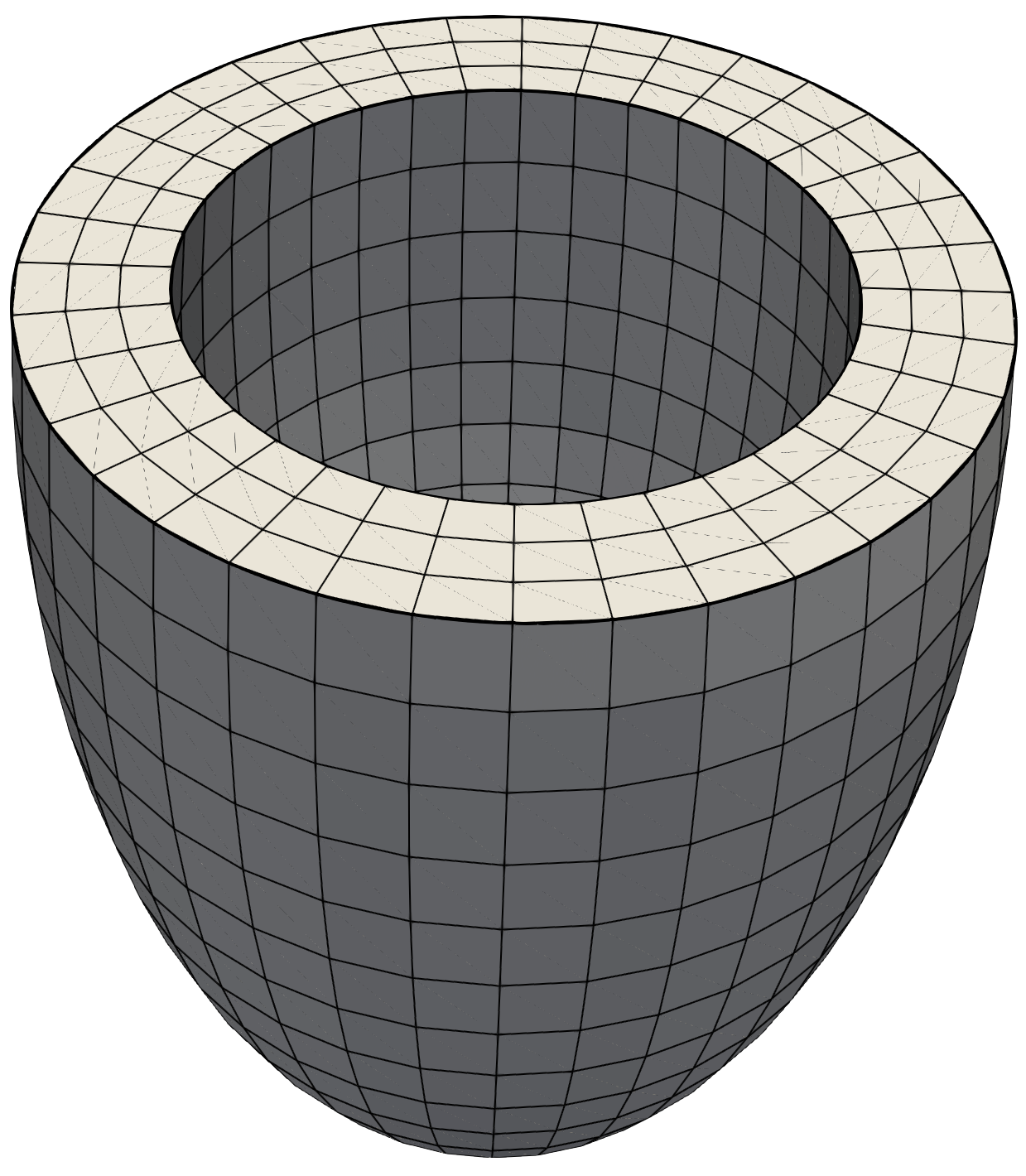} 
	}
	\subfigure[Cross-section showing the equivalent (von Mises) stress distribution for the mesh with 829,440 cells]
	{
		\includegraphics[height=0.48\textwidth]{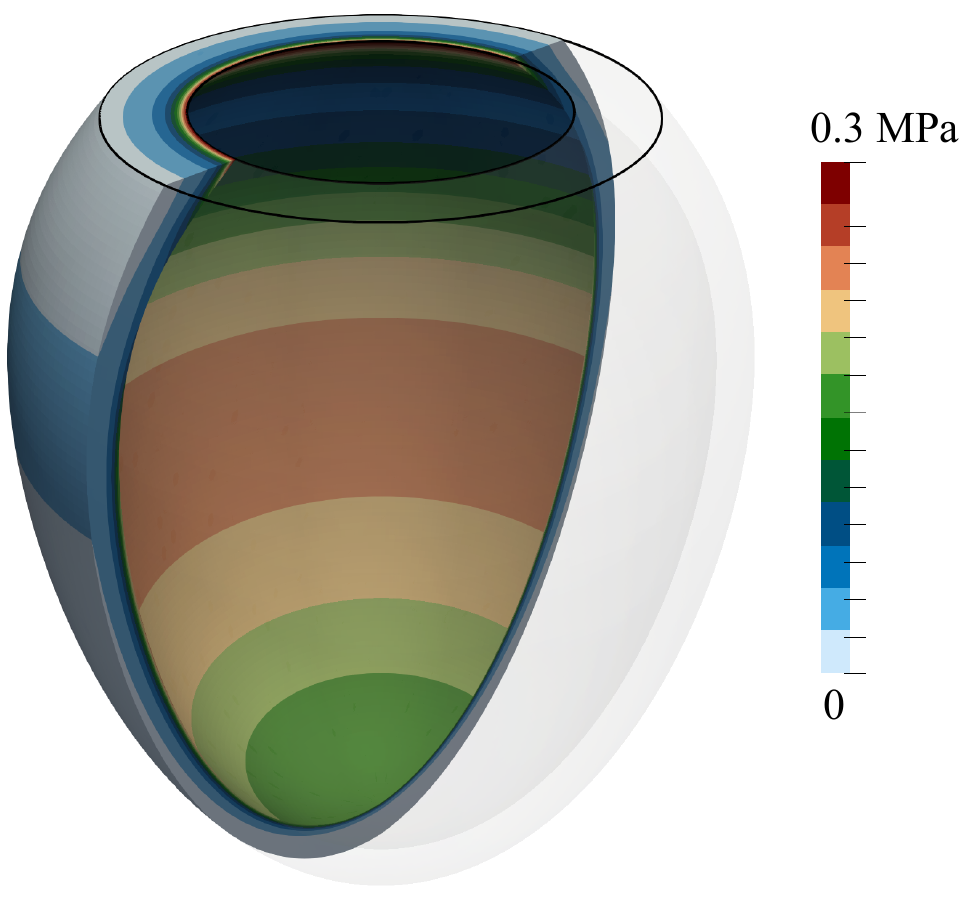}
	}
   \caption{Idealised ventricle case}
   \label{fig:ventricle}
\end{figure}

\paragraph{Case 4: Inflation of an idealised ventricle}
Inflation of an idealised ventricle (Figure \ref{fig:ventricle}) was proposed by \citet{Land2015} as a benchmark problem for cardiac mechanics software.
The case is 3-D, static, with finite hyperelastic strains.
The initial geometry is defined as a truncated ellipsoid:
\begin{eqnarray}
	x = r_s \sin(u) \cos(v), \quad
	y = r_s \sin(u) \sin(v), \quad
	z = r_l \cos(u)
\end{eqnarray}
where on the inner (endocardial) surface $r_s =7$ mm, $r_l = 17$ mm, $u \in \left[-\pi, -\arccos \left( \frac{5}{17} \right) \right]$ and $v \in \left[-\pi, \pi \right]$, while on the outer (epicardial) surface $r_s =10$ mm, $r_l = 20$ mm, $u \in \left[-\pi, -\arccos \left( \frac{5}{20} \right) \right]$ and $v \in \left[-\pi, \pi \right]$.
The base plane $z = 5$ mm is implicitly defined by the ranges for $u$.
The hyperelastic material behaviour is described by the transversely isotropic constitutive law proposed by \citet{Guccione1995} law, where the parameters are $C = 10$ kPa, $c_f = c_t = c_{fs} = 1$; the specified parameters produce isotropic behaviour.
The benchmark specifies an incompressible material; in the current work, incompressibility is enforced weakly using a penalty approach, where the bulk modulus parameter is chosen to be two orders of magnitude greater ($\kappa = 1000$ kPa) than the greatest shear modulus parameter.
A pressure of 10 kPa is applied to the inner surface over 100 equal load increments, and the base plane is fixed.
The geometry is meshed using a revolved, structured approach and is predominantly composed of hexahedra, with prism cells forming the apex.
The OpenFOAM utilities \texttt{blockMesh} and \texttt{extrudeMesh} are used to create the meshes.
Four successively refined meshes are examined: $1\,620$ (shown in Figure \ref{fig:ventricle}(a)), $12\,960$, $103\,680$, and $829\,440$ cells.
The case can be simulated as 2-D axisymmetric but is simulated here using the full 3-D geometry.

\begin{figure}[htbp]
	\centering
   		\includegraphics[width=0.6\textwidth]{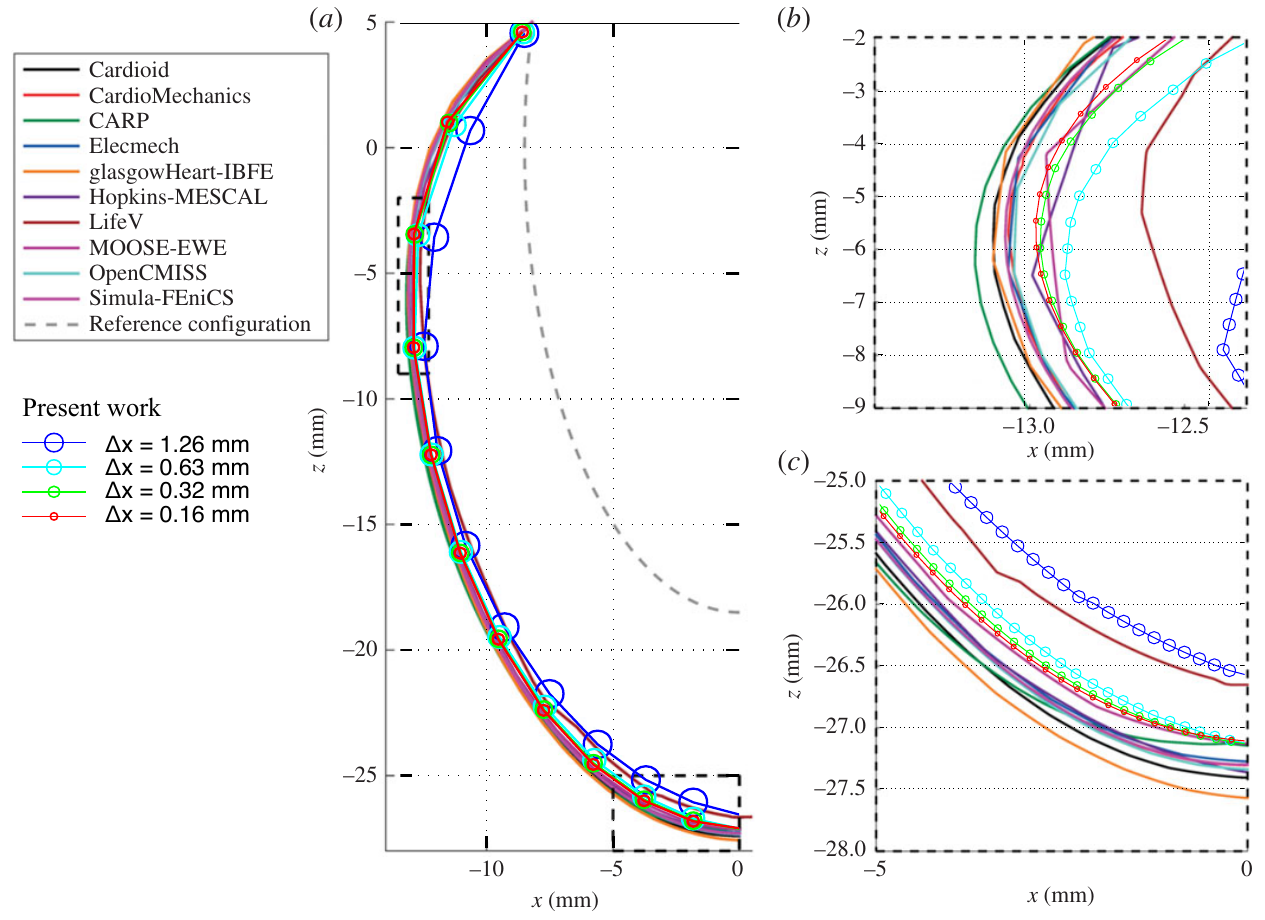} 
	\caption{Idealised ventricle: predictions for the four meshes from the present work overlaid on the round-robin predictions from \citet{Land2015}}
	\label{fig:ventricle_accuracy}
\end{figure}

The equivalent (von Mises) stress at full inflation is shown for the mesh with 829,440 cells in Figure \ref{fig:ventricle}(b), where the high stresses on the inner (endocardial) surface quickly drop off through the wall thickness towards the outer (endocardial) surface.
The predicted deformed configuration of the ventricle wall midline is shown in Figure \ref{fig:ventricle_accuracy}, where a side-by-side comparison is given with the round-robin benchmark results from \citet{Land2015}.
The predictions are seen to become quickly mesh-independent and fall within the benchmark ranges.

\begin{figure}[htbp]
	\centering
	\subfigure[Case geometry and dimensions (in mm)]
	{
		\label{fig:cooks_membrane_geometry}
   		\includegraphics[scale=1]{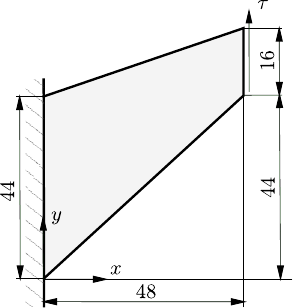} 
   	}
   	\qquad
	\subfigure[Quadrilateral mesh with 144 cells]
	{
		\label{fig:cooks_membrane_mesh}
   		\includegraphics[scale=0.14]{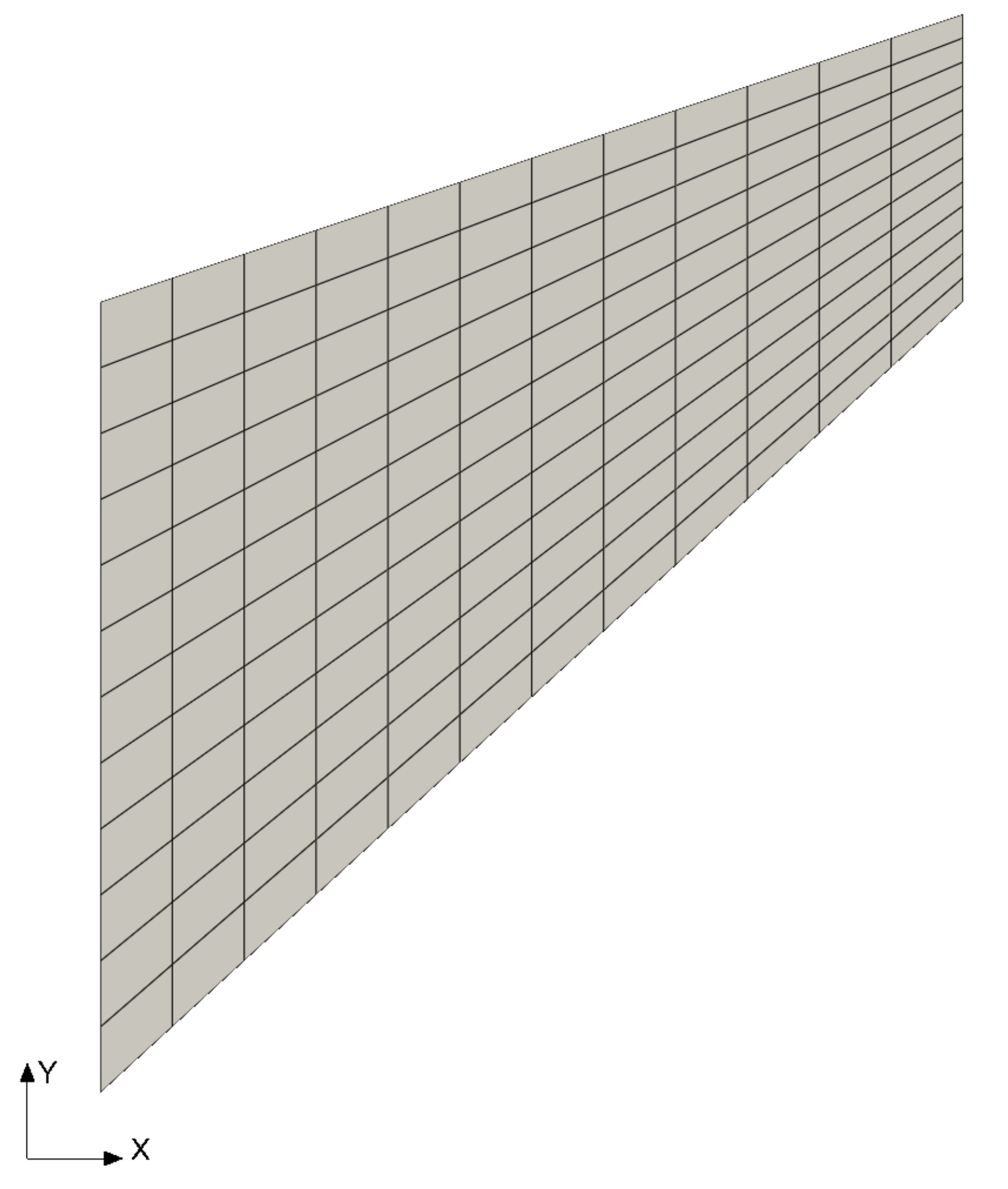}  
   	}
	\caption{Cook's membrane case geometry and mesh}
   \label{fig:cooks_membrane}
\end{figure}

\paragraph{Case 5: Cook's membrane}
Cook's membrane (Figure \ref{fig:cooks_membrane}) is a well-known bending-dominated benchmark case used in linear and non-linear analysis.
The 2-D plane-strain tapered panel (trapezoid) is fixed on one side and subjected to uniform shear traction on the opposite side.
The vertices of the trapezoid (in mm) are (0, 0), (48, 44), (48, 60),  and (0, 44).

The current work considers three forms of the problem:
\begin{enumerate}[label=\roman*.]
	\item Small strain, linear elastic \cite{Zienkiewicz2000, Simplas}: $E=70$ MPa, $\nu=1/3$, and $\tau = 6.25$ kPa.
	\item Finite strain, neo-Hookean hyperelastic \cite{Pelteret2018}: $E=1.0985$ MPa, $\nu=0.3$, and $\tau = 62.5$ kPa.
	\item Finite strain, neo-Hookean hyperelastoplastic \citep{Simo1992, Simplas, Cesar2001}: $E=206.9$ MPa, $\nu=0.29$, $\sigma_y(\bar{\varepsilon}_p) = 0.45 + 0.12924\bar{\varepsilon}_p + (0.715 - 0.45)(1- e^{-16.93\bar{\varepsilon}_p})$ MPa, and $\tau = 312.5$ kPa.
\end{enumerate}
where $E$ is the Young's modulus, $\nu$ is the Poisson's ratio, $\tau$ is the prescribed shear traction, and $\sigma_y(\bar{\varepsilon}_p)$ is the yield strength, which is a function of the equivalent plastic strain $\bar{\varepsilon}_p$.
Eight successively refined quadrilateral meshes are considered, where the cell counts are 9, 36, 144 (Figure \ref{fig:cooks_membrane}(b)), 576, $2\,304$, $9\,216$, $36\,864$, and $147\,456$.
The problem is solved quasi-statically, with no body forces.
The linear elastic case is solved in one loading increment, whereas the hyperelastic and hyperelastoplastic cases use 30 equally sized loading increments.

\begin{figure}[htbp]
   \centering
	\subfigure[Linear elastic case]
	{
		\includegraphics[width=0.28\textwidth]{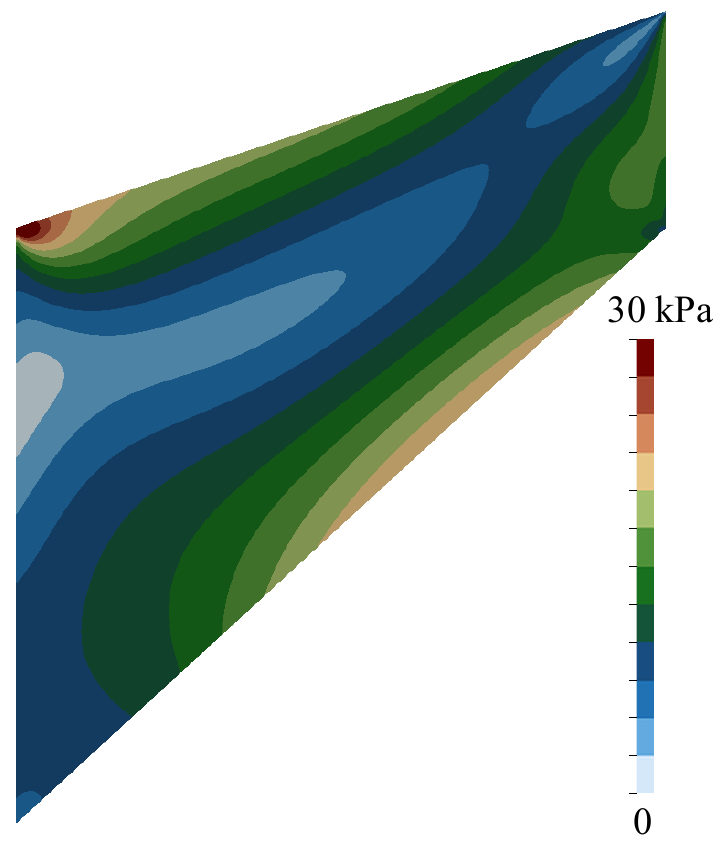}  
	}
	\subfigure[Hyperelastic case]
	{
		\includegraphics[width=0.3\textwidth]{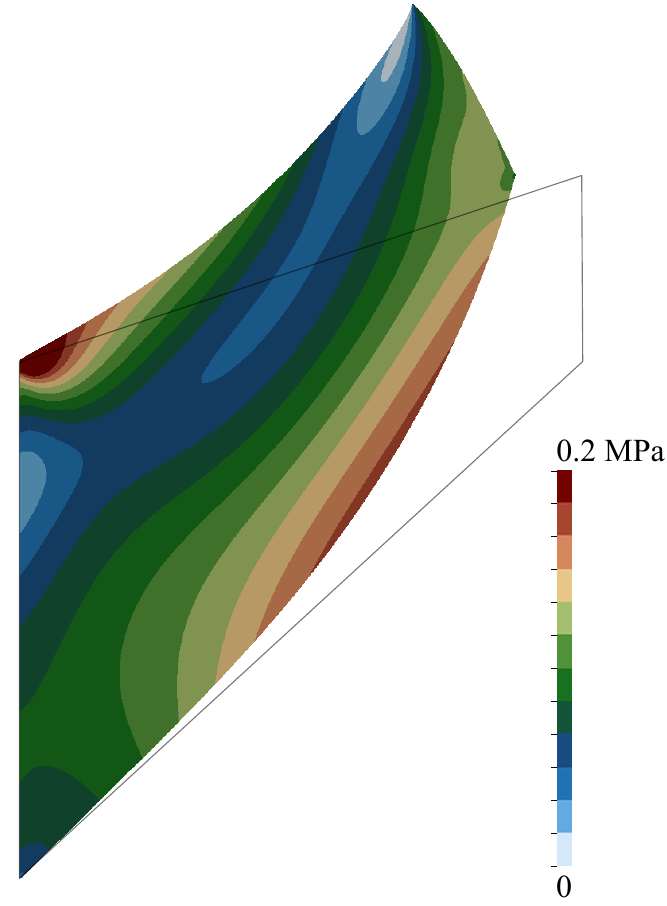}  
	}
	\subfigure[Hyperelastoplastic case]
	{
		\includegraphics[width=0.3\textwidth]{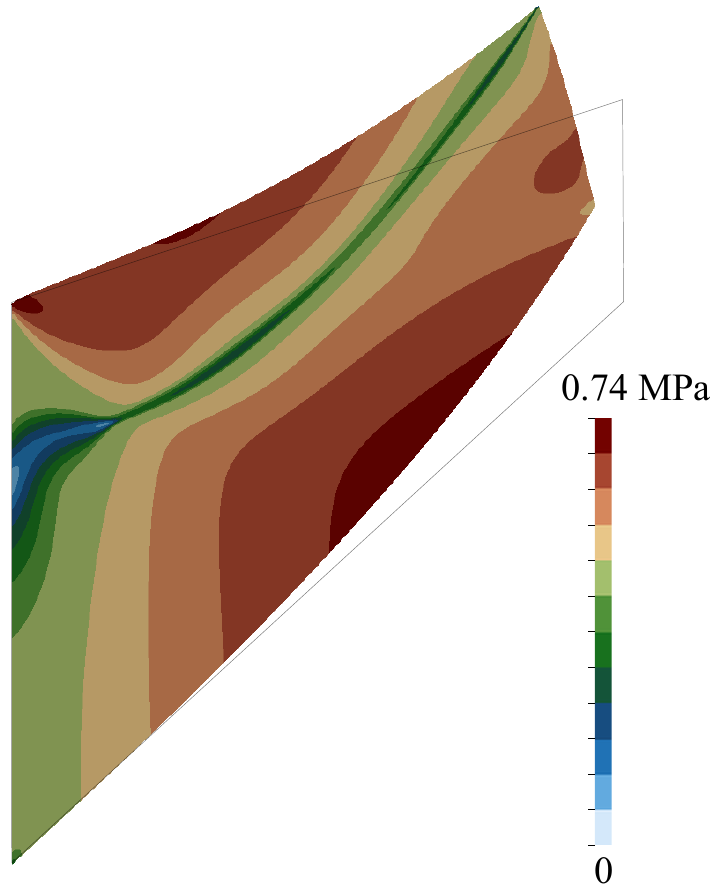}  
	}
   \caption{Equivalent (von Mises) stress distribution for the three Cook's membrane cases using the mesh with $36\,842$ cells}
   \label{fig:cooksMembrane_sigmaEq}
\end{figure}

Figure \ref{fig:cooksMembrane_sigmaEq} shows the predicted equivalent stress distribution on the mesh with $36\,842$ cells for linear elastic, hyperelastic and hyperelastoplastic cases.
The stress distribution is consistent with bending, with regions of high stress near the upper and lower surfaces and a line of relatively unstressed material in the centre.
The greatest equivalent stresses occur at the top-left corner and the lower surface in all versions of the case.
In the hyperelastoplastic case, almost the entire domain is plastically yielding with only a small thin region remaining elastic, indicated by the blue line in Figure \ref{fig:cooksMembrane_sigmaEq}(c).

\begin{figure}[htbp]
   \centering
   	\subfigure[Linear elastic case]
	{
		\includegraphics[height=0.3\textwidth]{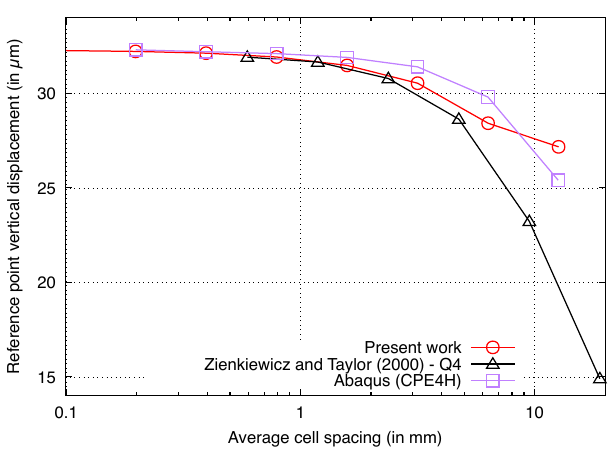}  
	}
	\subfigure[Hyperelastic case]
	{
		\includegraphics[height=0.3\textwidth]{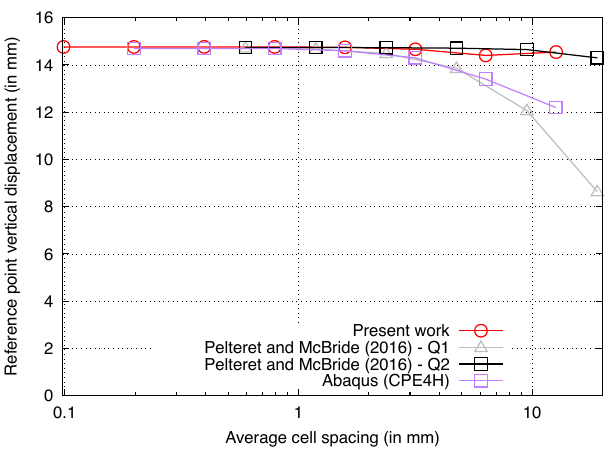}  
	}
	\subfigure[Hyperelastoplastic case]
	{
		\includegraphics[height=0.3\textwidth]{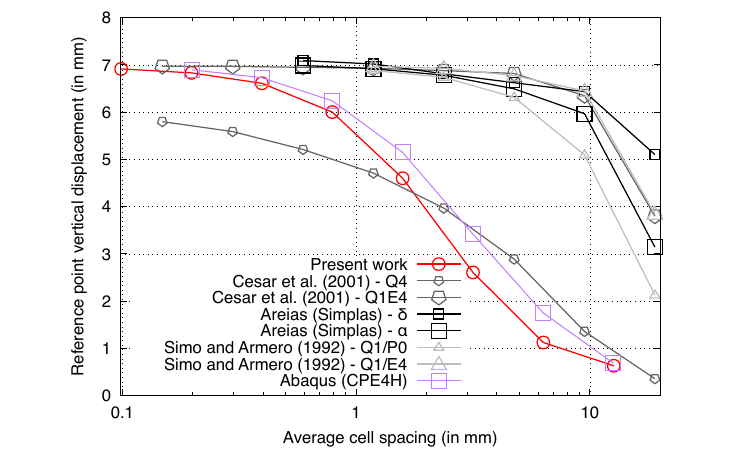}  
	}
   \caption{Cook's membrane vertical displacement predictions at the reference point as a function of the average cell width. Comparisons are given with the results from the literature \cite{Zienkiewicz2000, Pelteret2018, Simo1992, Simplas}}
   \label{fig:cooksMembrane_disp}
\end{figure}

Figure \ref{fig:cooksMembrane_disp} compares the predicted vertical displacement at the reference point as a function of the average cell widths with results from the literature and finite element software Abaqus (element code C3D8).
The reference point is taken as the top right point -- $(48,60)$ mm -- in the elastic and hyperelastoplastic cases, while it is taken as the midway point of the loading surface -- $(48,52)$ mm -- in the hyperelastic case.

The results shown for the hyperelastoplastic case were generated using the segregated solution procedure, as the Jacobian-free Newton-Krylov approach diverged for six of the eight mesh densities, as discussed further in Section \ref{sec:resource_requirements}.


\paragraph{Case 6: Vibration of a 3-D Cantilevered Beam}
This 3-D, dynamic, finite strain case geometry consists of a $2 \times 0.2 \times 0.2$ m cuboid column and was proposed in its initial form by \citet{Tukovic2007}.
A sudden, constant traction $\bb{T} = \left(50, 50, 0 \right)$ kPa is applied to the upper surface.
A neo-Hookean hyperelastic material is assumed with $E = 15.293$ MPa, $\nu = 0.3$ and density $\rho = 1000$ kg m$^{-3}$, while the area $A = 0.04$ m$^2$ and the second moment of area $I = \nicefrac{1}{7\,500}$ m$^4$. 
The geometric and material parameters were chosen to yield a first natural frequency of 1 Hz in the small-strain limit.
The magnitude of the applied traction is chosen here to test the robustness of the solution procedures under large rotations and strains.
Four succesively refined hexahedral meshes are generated using the OpenFOAM \texttt{blockMesh} utility, with cell counts of 270, $2\,160$, $17\,280$ and $138\,240$.
The time step size is 1 ms, and the total period is 1 s.

\begin{figure}[htbp]
   \centering
	\subfigure[Displacement magnitude for t = $\left\{0, \, 0.1, \, 0.2\right\}$ s (translucent) and t = $\left\{0.3, \, 0.4, \, 0.5 \right\}$ s (opaque) for the mesh with $138\,240$ cells]
	{
	   \includegraphics[width=0.5\textwidth]{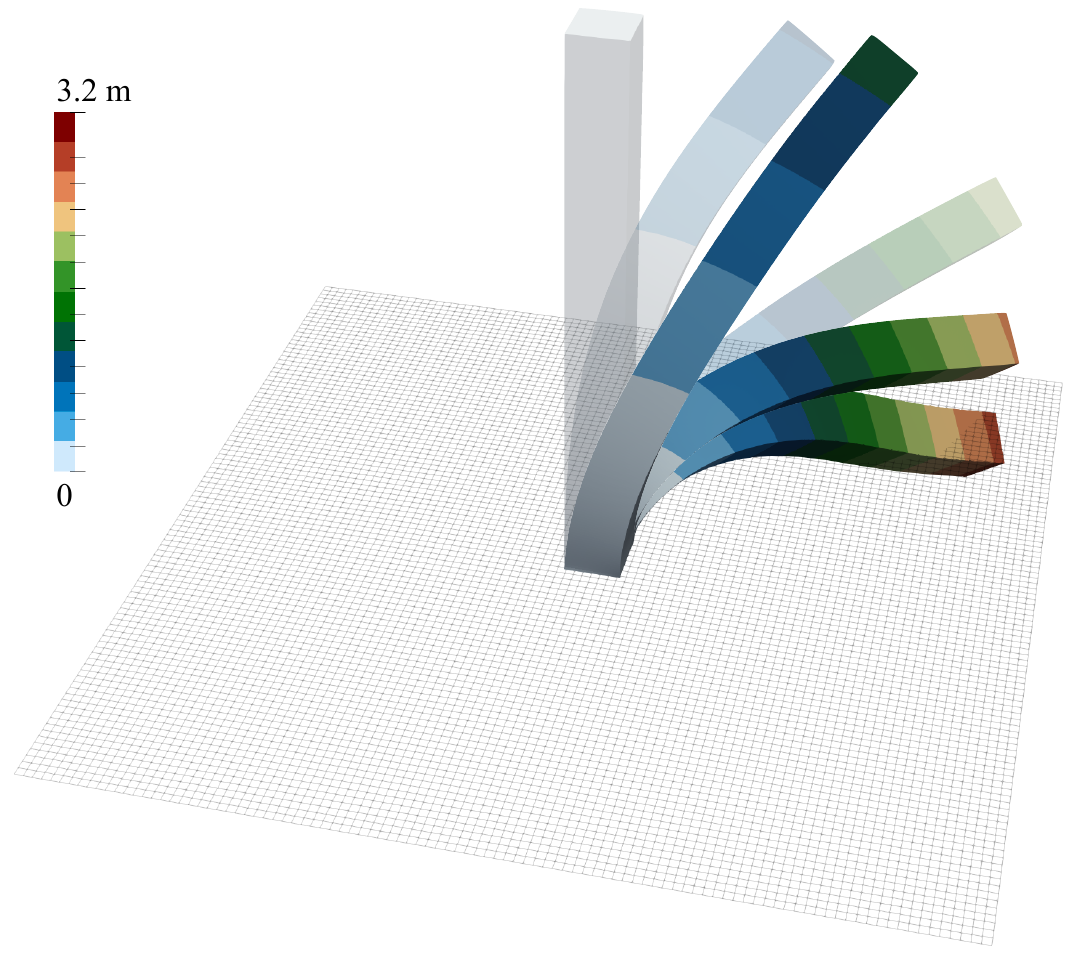}
   	} \quad
	\subfigure[Displacement magnitude of the centre of the upper surface of the beam vs time]
	{
	   \includegraphics[width=0.6\textwidth]{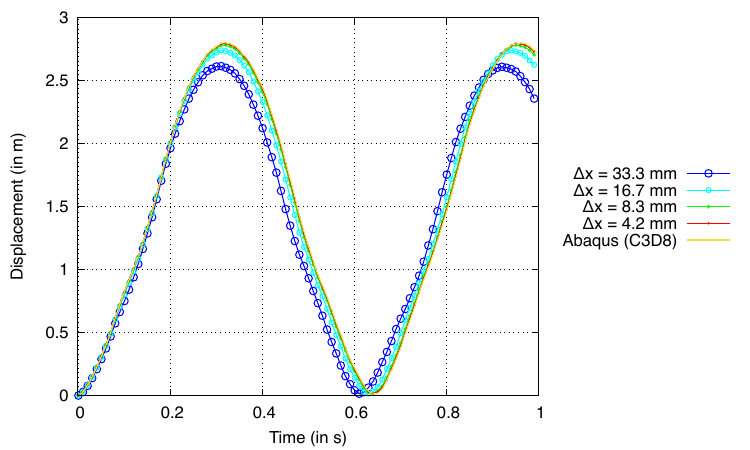} 
   	}
   \caption{Deflection of a 3-D Cantilevered Beam}
   \label{fig:dynamic_cantilever}
\end{figure}

The deformed configuration of the beam at six time steps is shown in Figure \ref{fig:dynamic_cantilever}(a) for the mesh with $2\,160$ cells, where the upper surface is seen to drop below the horizontal from approximately t = 0.2 s to 0.32 s.
The corresponding displacement magnitude of the centre of the upper surface of the beam vs time is shown in Figure \ref{fig:dynamic_cantilever}(b), where the predictions are seen to converge to the results from finite element software Abaqus (element code C3D8) using the $138\,240$ mesh.

\begin{figure}[htbp]
   \centering
	\subfigure[Mesh containing $5\,336$ cells in the fluid region and 630 cells in the solid region]
	{
   		\includegraphics[width=0.7\textwidth]{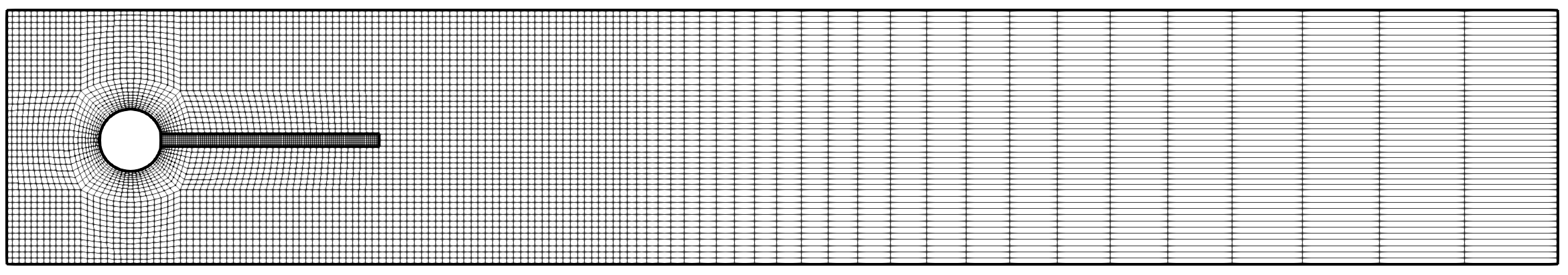}
	}
	\subfigure[Close-up of the mesh around the plate]
	{
   		\includegraphics[width=0.4\textwidth]{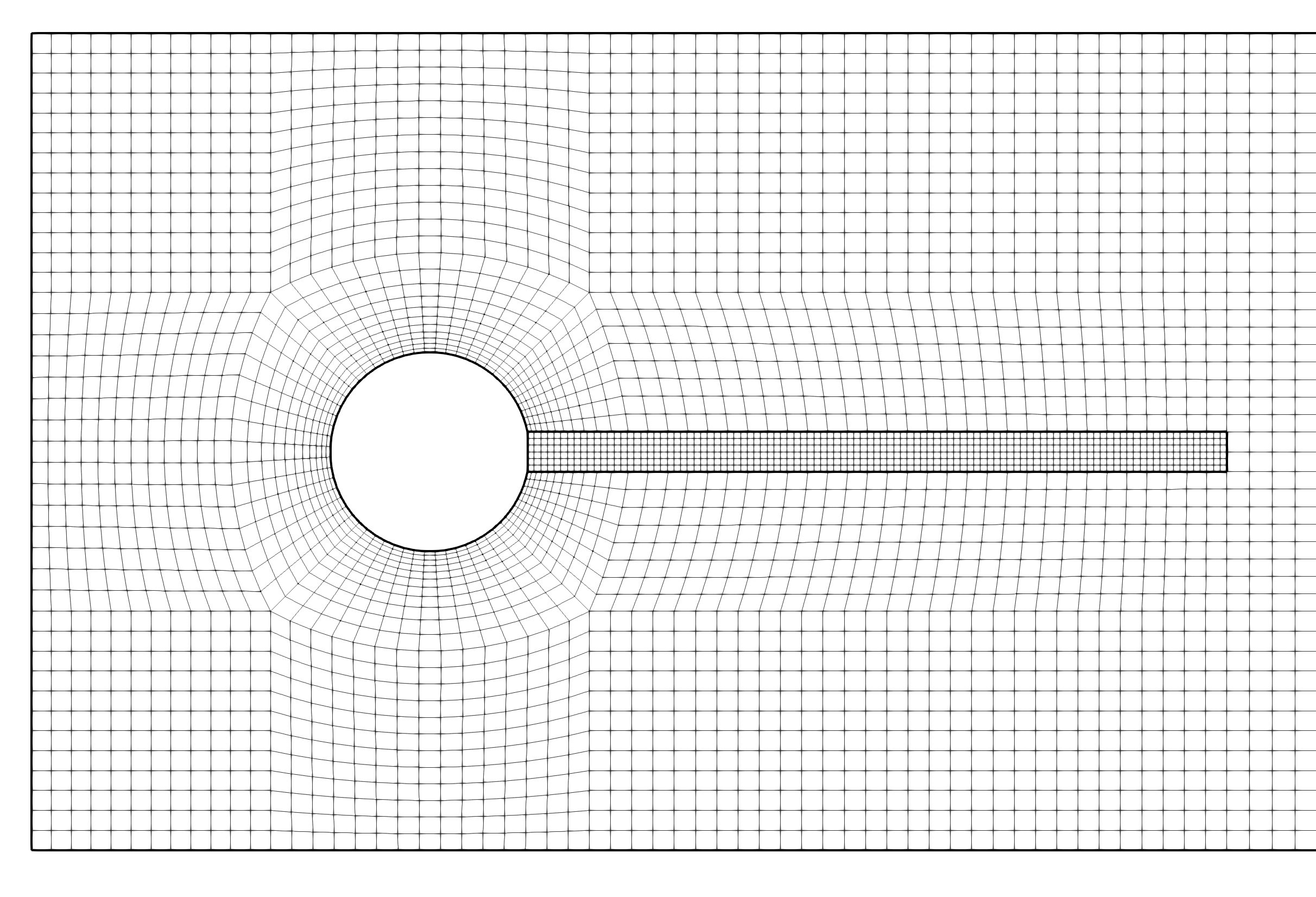}  
   	}
   \caption{Elastic plate behind a rigid cylinder case geometry and mesh}
   \label{fig:hronTurek-mesh}
\end{figure}

\paragraph{Case 7: Elastic plate behind a rigid cylinder}
The final test case, introduced by \citet{Turek2006}, extends the classic flow around a cylinder in a channel problem \citep{Ferziger2002} into a well-established fluid-solid interaction benchmark.
The \texttt{FSI3} variant of the case examined here (Figure \ref{fig:hronTurek-mesh}) consists of a horizontal channel (0.41 m in height, 2.5 m in length) with a rigid cylinder of radius 0.05 m, where the cylinder centre is 0.2 m from the bottom and inlet (left) boundaries.
A parabolic velocity is prescribed at the inlet velocity with a mean value of 2 m/s.
A St.\,Venant-Kirchhoff hyperelastic plate ($E = 5.6$ MPa, $\nu = 0.4$) of 0.35 m in length and 0.02 m in height is attached to the right-hand side of the rigid cylinder.
The fluid model is assumed to be isothermal, incompressible, and laminar and adopts the segregated PIMPLE solution algorithm.
The fluid's kinematic viscosity is 0.001 m$^2$/s and density is 1000 kg m$^{-3}$, while the solid's density is $10\,000$ kg m$^{-3}$.
The interface quasi-Newton coupling approach with inverse Jacobian from a least-squares model \cite{Degroote2009} is employed for the fluid-solid interaction coupling.
The fluid-solid interface residual is reduced by four orders of magnitude within each time step.
Further details of the fluid-solid interaction procedure are found in \citet{Tukovic2018}.
Three successively refined quadrilateral meshes are employed:
The fluid region meshes have $1\,252$, $5\,008$, and $20\,032$ cells, while the solid region meshes have 156,  624 and $2\,496$ cells.
The total simulation time is 20 s, and the time-step size is 0.5 ms. 

\begin{figure}[htbp]
   \centering
	   \includegraphics[width=\textwidth]{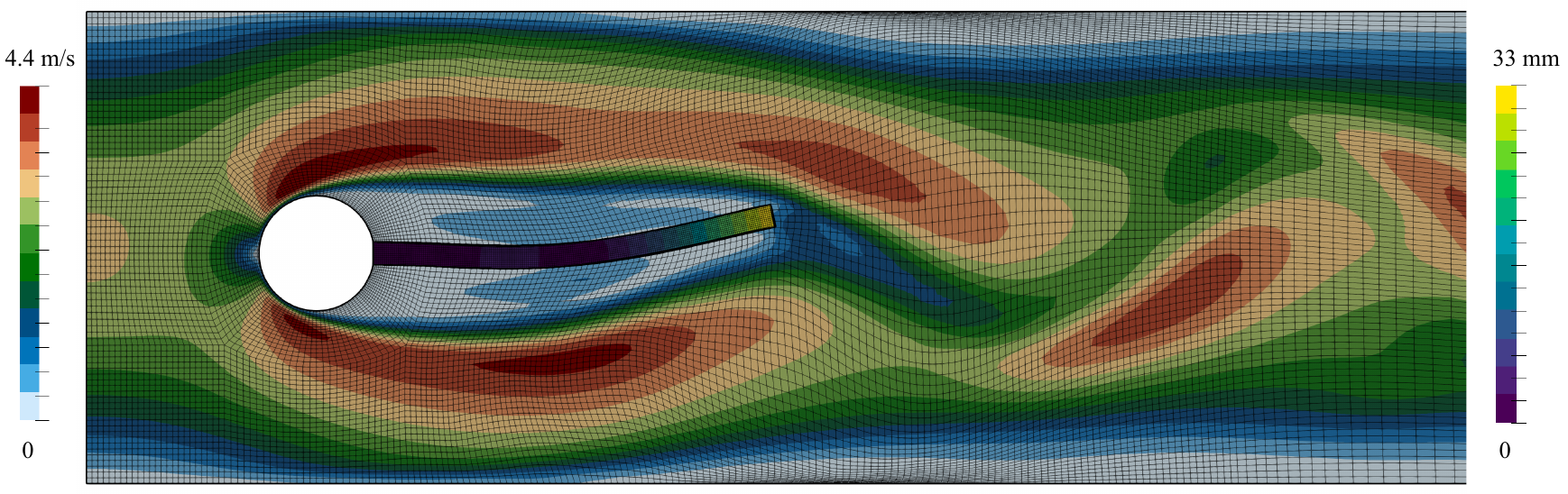} 
   \caption{Elastic plate behind a rigid cylinder case results at $t = 20$ s for the mesh with $20\,032$ cells in the fluid region and $2\,496$ in the solid region. The velocity magnitude is shown in the fluid region, while the displacement magnitude is shown in the solid region.}
   \label{fig:hronTurek-results}
\end{figure}

Figure \ref{fig:hronTurek-results} shows the velocity field in the fluid region and displacement magnitude in the solid region at $t = 4.42$ s, corresponding to a peak in the vertical displacement of the plate free-end oscillation.
The fluid accelerates as it passes the cylinder, with regular vortices being shed from the plate, causing it to oscillate.
The predicted mean, amplitude and frequency of the vertical and horizontal displacement of the end of the plate are seen to approach the results from \citet{Turek2006} in Table \ref{tab:hronTurekDisp}.
As directed in \citet{Turek2006}, the maximum and minimum values for the last period are used to calculate the mean and amplitude:
\begin{eqnarray}
	\text{mean} = \frac{1}{2}(\text{max} + \text{min}), \quad\quad
	\text{amplitude} = \frac{1}{2}(\text{max} - \text{min})
\end{eqnarray}
while the frequency is calculated as the inverse of the period.
\begin{table}[htb]
	\centering
		\begin{tabular}{lll}
			\hline
			Mesh (number of fluid and solid cells) & $u_x$ (in mm) & $u_y$ (in mm) \\
			\hline
			1 ($1252$ + $156$)  & $-0.10 \pm 0.03 \,[ 58.82 ]$ & $1.60 \pm 0.00 \,[ 5.78 ]$ \\
2 ($5008$ + $624$)  & $-2.10 \pm 2.02 \,[ 10.42 ]$ & $1.44 \pm 29.68 \,[ 5.59 ]$ \\
3 ($20032$ + $2496$)  & $-2.46 \pm 2.30 \,[ 10.75 ]$ & $1.46 \pm 32.63 \,[ 5.68 ]$ \\
\hline

			\citet{Turek2006} & $-2.69 \pm 2.53\,[10.9]$ & $1.48 \pm 34.38\,[5.3]$ \\
			\hline
		\end{tabular}
	\caption{Predicted oscillations (mean $\pm$ amplitude $[$frequency$]$) of the free end of the plate. Note that large-scale plate oscillations in mesh 1 (coarsest mesh) die out within the first 10 s.}
		\label{tab:hronTurekDisp}
\end{table}

\subsection{Resource Requirements and Robustness}
\label{sec:resource_requirements}
This section compares the resource requirements and robustness of the Jacobian-free Newton-Krylov and segregated approaches across the cases presented in Section \ref{sec:accuracy}.
Specifically, time and memory requirements are analysed, and whether the solution converged for all meshes and loading steps.
In all cases, the Jacobian-free Newton-Krylov approach used the GMRES linear solver with the direct LU preconditioner for 2-D and the Hypre BoomerAMG multigrid preconditioner for 3-D, unless stated otherwise.
A restart value of 30 is used for the GMRES solver, where the \emph{loose} variant of the GMRES solver is used, which retains a small number of extra solution vectors across restarts (2 in this case).
In contrast, the segregated approach used the conjugate gradient linear solver and the incomplete Cholesky preconditioner with no infill -- equivalent to ILU(0).
Time and memory requirements are implementation- and hardware-specific; nonetheless, it is insightful to compare the relative performance of the Jacobian-free Newton-Krylov and segregated approaches on the same hardware and within the same implementation framework.
Clock times and memory usage (measured with the GNU time utility) were generated using a Mac Studio with an M1 Ultra CPU on one core, where the code was built with the Clang compiler (version 16.0.0).
Examination of multi-CPU-core parallelisation is left to Section \ref{sec:parallelisation}.

In the AMG preconditioner, each outer GMRES iteration corresponded to a single BoomerAMG V-cycle.
On each grid level, one pre-smoothing and one post-smoothing sweep were applied.
Coarsening was performed using the HMIS strategy, which provides robust, scalable coarse grids for 3-D problems. Interpolation between levels used the extended+i scheme with truncation (\texttt{truncfactor = 0.3}) and a cap on interpolation stencil size (\texttt{P$_{max}$ = 1}), ensuring sparsity of the prolongation operator. Strong connections were defined with a threshold of 0.7, which balances coarse-grid accuracy with complexity. Aggressive coarsening was enabled (one pass with a single aggregation path), and the multigrid hierarchy was allowed up to 25 levels, though in practice far fewer were required.
Collectively, these settings aim to yield a relatively light yet robust AMG preconditioner, consistent with recommendations in the MOOSE finite element software documentation \citep{moose_hypre_doc_2022}.

Table \ref{tab:times_memory} lists the wall clock times (time according to a clock on the wall) and maximum memory usage for all cases and all meshes, where diverged cases are indicated by the symbol $\dag$.
Figure \ref{fig:times_memory}(a) plots the corresponding speedup on all cases as a function of the number of degrees of freedom, where the speedup is defined as the segregated clock time divided by the Jacobian-free Newton-Krylov clock time.
The number of degrees of freedom is $2 \times$ the cell count in 2-D cases and $3 \times$ in 3-D cases, except in the fluid-solid interaction case, where the fluid domain has three degrees of freedom (two velocity components and pressure) in 2-D.
Speedup is calculated only for clock times greater than 2 seconds to avoid comparing small numbers.
In all cases where convergence was achieved, the Jacobian-free Newton-Krylov approach was faster (speedup $> 1$), with speedups of one or two orders of magnitude in many cases.
Additionally, it can be observed from Figure \ref{fig:times_memory}(a) that the speedup increases as the number of degrees of freedom increases.
The largest speedup was 341 on the sixth mesh ($9\,216$ cells) of the membrane i (2-D, linear elastic) case, while the lowest speed was 1 and occurred in the manufactured solution coarser meshes.
\begin{table}[!htbp]
	\centering
		\begin{tabular}{ll|ll|ll}
			\hline
			\textbf{Case} & \textbf{Cell Count} & \multicolumn{2}{c|}{\textbf{JFNK}} & \multicolumn{2}{c}{\textbf{Segregated}} \\
			     &            & \textbf{Time} & \textbf{Memory} & \textbf{Time} & \textbf{Memory} \\
			     &            & (in s) & (in MB) & (in s) & (in MB) \\
			\hline
				\textbf{MMS}		& 125	        & 0 & 76 & 0 & 66 \\
	\emph{regular hex}	& $1\,000$	& 0 & 85 & 0 & 77 \\
	\emph{3-D, static}	& $8\,000$	& 1 & 156 & 1 & 118 \\
	\emph{linear elastic}	& $64\,000$	& 3 & 967 & 3 & 463 \\
				& $512\,000$	& 33 & $3\,865$ & 37 & $2\,014$ \\
				& $4\,096\,000$	& 373 & $15\,762$ & 395 & $11\,457$ \\
\hline

			\textbf{Spherical}	& 976	& 0 & 103 & 1 & 90 \\
	\textbf{Cavity}	& $4\,552$	& 1 & 204 & 1 & 143 \\
	\emph{3-D, static}	& $29\,611$	& 6 & 982 & 9 & 543 \\
	\emph{linear elastic}	& $213\,100$	& 87 & $5\,804$ & 135 & $1\,620$ \\
				& $1\,614\,261$	& $1\,265$ & $27\,173$ & $2\,367$ & $9\,017$ \\
\hline

				\textbf{Elliptic Plate}	& $45$	        & 0 & 75 & 0 & 65 \\
	\emph{3-D static}	& $472$	        & 0 & 80 & 1 & 69 \\
	\emph{linear elastic}	& $4\,140$	& 0 & 131 & 2 & 109 \\
				& $34\,968$	& 8 & 709 & 35 & 415 \\
				& $287\,280$	& 107 & $3\,852$ & 475 & $1\,497$ \\
				& $2\,438\,242$	& $1\,073$ & $22\,957$ & $7\,176$ & $7\,021$ \\
\hline

				\textbf{Ventricle}	& $1\,620$	& 52 & 150 & $\dag$ & 123 \\
	\emph{3-D, static}	& $12\,960$	& 608 & 595 & $\dag$ & 415 \\
	\emph{hyperelastic}	& $103\,680$	& $11\,893$ & 2350 & $\dag$ & 1999 \\
				& $829\,440$	& $118\,944$ & 7957 & $\dag$ & 5589 \\
\hline

				\textbf{Membrane i}	& 9	& 0 & 76 & 0 & 65 \\
	\emph{2-D static}	& 36	& 0 & 76 & 1 & 65 \\
	\emph{linear elastic}	& 144	& 0 & 77 & 1 & 66 \\
				& 576	& 0 & 81 & 1 & 73 \\
				& $2\,304$	& 0 & 99 & 1 & 81 \\
				& $9\,216$	& 0 & 174 & 9 & 143 \\
				& $36\,864$	& 2 & 606 & 67 & 422 \\
				& $147\,456$	& 9 & $1\,605$ & 532 & $1\,384$ \\
\hline

				\textbf{Membrane ii}	& 9	& 0 & 76 & 1 & 65 \\
	\emph{2-D, static,}	& 36	& 0 & 76 & 2 & 65 \\
	\emph{hyperelastic}	& 144	& 0 & 77 & 7 & 65 \\
				& 576	& 0 & 82 & 52 & 67 \\
				& $2\,304$	& 2 & 106 & 503 & 89 \\
				& $9\,216$	& 8 & 213 & $2\,725$ & 188 \\
				& $36\,864$	& 37 & 686 & $1\,795$ & 469 \\
				& $147\,456$	& 170 & $2\,118$ & $15\,573$ & 910 \\
\hline

				\textbf{Membrane iii}	& 9	& 0 & 76 & 2 & 66 \\
	\emph{2-D, static,}	& 36	& 1 & 77 & 6 & 66 \\
	\emph{hyperelastoplastic}	& 144	& $\dag$ & 79 & 29 & 67 \\
				& 576	& $\dag$ & 83 & 167 & 73 \\
				& $2\,304$	& $\dag$ & 103 & 894 & 102 \\
				& $9\,216$	& $\dag$ & 218 & $6\,936$ & 262 \\
				& $36\,864$	& $\dag$ & 813 & $57\,753$ & $1\,217$ \\
				& $147\,456$	& $\dag$ & $1\,832$ & $522\,811$ & $3\,216$ \\
\hline

				\textbf{Cantilever}	& $270$	& 12 & 85 & 11 & 73 \\
	\emph{3-D, dynamic}	& $2\,160$	& 87 & 166 & 93 & 94 \\
	\emph{hyperelastic}	& $17\,280$	& $1\,028$ & 739 & $1\,459$ & 591 \\
				& $138\,240$	& $12\,156$ & $5\,791$ & $\dag$ & $1\,838$ \\
\hline

				\textbf{FSI}	& $1\,252$ + 156			& $3\,591$ & 773 & $5\,883$ & 869 \\
	\emph{2-D, dynamic,}	& $5\,008$ + 624	& $25\,920$ & 911 & $43\,365$ & 884 \\
	\emph{hyperelastic}	& $20\,032$ + $2\,496$	& $154\,938$ & 908 & $199\,143$ & 865 \\
\hline

		\end{tabular}
	\caption{Execution times (rounded to the nearest second) and maximum memory usage for Jacobian-free Newton-Krylov (JFNK) and segregated methods (rounded to the nearest MB). $\dag$ indicates the solver diverged. }
	\label{tab:times_memory}
\end{table}

It is observed that the computational time for some 3-D cases (e.g., MMS) increases superlinearly with problem size, despite the use of AMG preconditioning.
In principle, AMG should provide optimal complexity; however, the HYPRE BoomerAMG implementation employed here is known to be sensitive to parameter choice.
In the present work, default parameters were adopted to ensure reproducibility across cases, which may not yield the best possible scalability.
Further tuning of AMG parameters is therefore expected to improve performance and will be the subject of future work.

\begin{figure}[htbp]
   \centering
	\subfigure[Speedup]
	{
   		\includegraphics[width=0.48\textwidth]{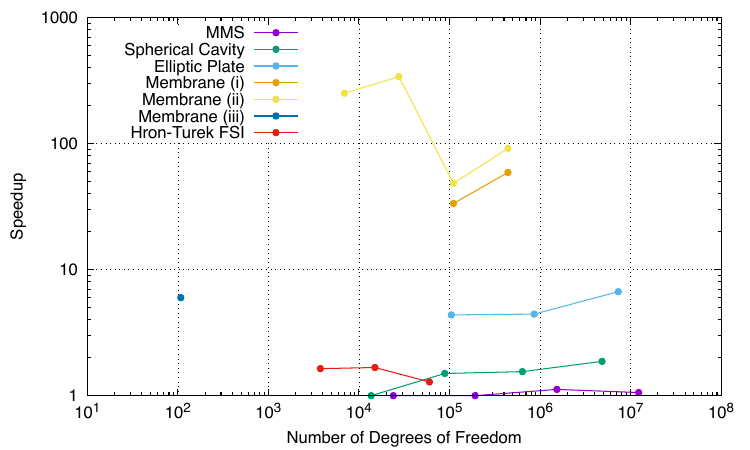}
	}
	\subfigure[Relative memory]
	{
   		\includegraphics[width=0.48\textwidth]{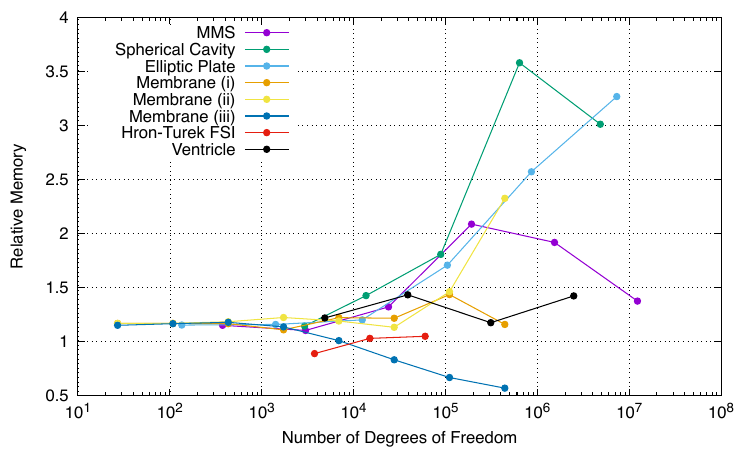}  
   	}
   \caption{The speedup (segregated clock time divided by the Jacobian-free Newton-Krylov clock time) and relative memory usage (segregated maximum memory usage divided by the Jacobian-free Newton-Krylov maximum memory usage) as a function of degrees of freedom}
   \label{fig:speedup_rel_mem}
\end{figure}

Figure \ref{fig:speedup_rel_mem}(a) shows the \emph{speedup} for each case, defined as the ratio of the time taken by the segregated solver to the time taken by the proposed Jacobian-free Newton-Krylov method.
Similarly, figure \ref{fig:speedup_rel_mem}(b) shows the \emph{relative memory} usage, defined here as the maximum memory usage of the Jacobian-free Newton-Krylov approach divided by that of the segregated approach.
The relative memory usage is close to unity for coarse meshes (less than $10^4$ degrees of freedom), but shows a general increasing trend beyond this.
The maximum relative memory value was 3.6 and occurred in the second-finest mesh for the spherical cavity case.
The variations in the relative memory trends across cases are likely due to the variable memory usage of the GMRES linear solver in the Jacobian-free Newton-Krylov, where the number of stored directions can vary depending on convergence.
In addition, the memory usage of the multigrid preconditioner can vary depending on the geometry and cell types.
The effect of the preconditioner and linear solver settings is examined further in Section \ref{sec:preconditioner}.

%

Regarding robustness, the Jacobian-free Newton-Krylov approach converged in all cases except one, requiring, on average, fewer than 5 outer Newton iterations.
The only case in which the Jacobian-free Newton-Krylov approach failed to converge was the membrane iii problem, the only elastoplastic case examined.
For the failed meshes, approximately 80\% of the total load was reached before the linear solver diverged.
The cause of this divergence is likely related to the proposed compact preconditioner matrix -- which is formulated based on small \emph{elastic} strains -- being a poor approximation of the true Jacobian for large plastic strains.
Interestingly, convergence problems were not encountered in the hyperelastic cases (membrane ii, ventricle, cantilever, fluid-solid interaction), demonstrating that the proposed compact \emph{linear elastic} preconditioner matrix is suitable for such large strain, large rotation hyperelastic cases.
In contrast, the segregated approach -- which uses the same matrix -- had no problems with the elastoplastic case (membrane iii) but failed to converge on two of the hyperelastic cases (ventricle, cantilever).
In both cases where the segregated approach failed, large elastic strains and large rotations were observed, and the segregated solver failed in the early time steps.
It is interesting to note that the same linear elastic preconditioner matrix works well with the Jacobian-free Newton-Krylov approach for hyperelastic cases but not for the segregated approach, while the opposite is true for elastoplastic cases.

Regarding geometric dimension, the same general trends are observed in 2-D and 3-D, with no major distinctions in behaviour.
The same can be said for geometric nonlinearity (small-strain vs. large-strain) with similar speed-ups in both linear elastic and hyperelastic cases.
One observation worth highlighting is the lower speedup seen for the method of manufactured solutions case (3-D, linear elastic): a possible explanation is that the segregated approach has previously been seen to be efficient on geometry with low aspect ratios (ratio of maximum to minimum dimensions) and minimal bending deformation; in this case, the aspect ratio is at its minimum (unity) and bending deformation is localised.
In contrast, in the other 3-D linear elastic cases (elliptic plate, membrane i), the aspect ratios are greater than unity and more widespread bending is present, and the segregated approach performs worse relative to the Jacobian-free Newton-Krylov approach.

\subsection{Effect of the Preconditioner Choice}
\label{sec:preconditioner}
This section examines the effect of preconditioning strategy on the performance of the Jacobian-free Newton Krylov approach.
Three choices of preconditioning procedure are compared:
\begin{itemize}
	\item \textbf{LU} - The \emph{MUltifrontal Massively Parallel sparse direct Solver} (MUMPS) \citep{MUMPS:1, MUMPS:2} LU decomposition direct solver. A direct solver is expected to be more robust but may suffer from excessive time and memory requirements for larger numbers of unknowns.
	\item \textbf{ILU($k$)} - Incomplete LU decomposition with $k$ fill-in. ILU($k$) is expected to have lower memory requirements than the LU direct solver, but at the expense of robustness.
	As the system of unknowns becomes larger, the number of ILU($k$) iterations is expected to increase. As the fill-in factor $k$ increases, ILU($k$) approaches the robustness of a LU direct solver. In the current section, $k = 5$ for all cases examined.
	\item \textbf{Algebraic multigrid} - The Hypre Boomerang \citep{hypre} parallelised multigrid preconditioner. Multigrid approaches have the potential to offer superior performance than other methods for larger problems, with near-linear scaling of time and memory requirements.
\end{itemize}

An additional consideration when selecting a preconditioner is its ability to scale in parallel as the number of CPU cores increases.
From this perspective, the iterative approaches (ILU($k$) and multigrid) are expected to show better parallel scaling than direct methods (LU), a point that is briefly examined in Section \ref{sec:parallelisation}.

\begin{figure}[htbp]
   \centering
	\subfigure[Clock time]
	{
   		\includegraphics[width=0.48\textwidth]{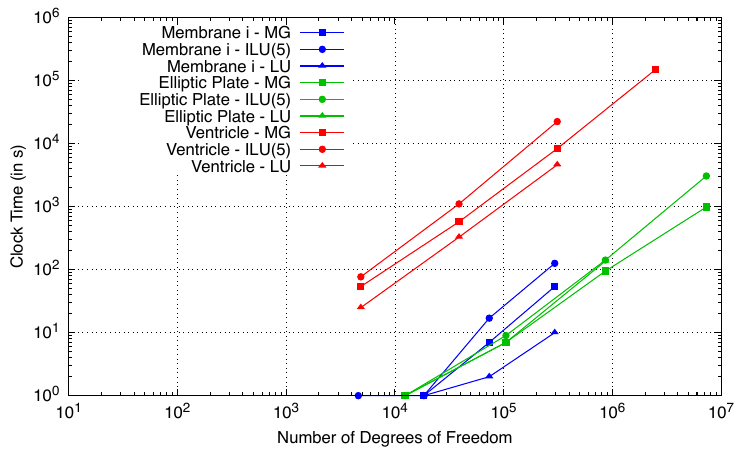}
	}
	\subfigure[Maximum memory]
	{
   		\includegraphics[width=0.48\textwidth]{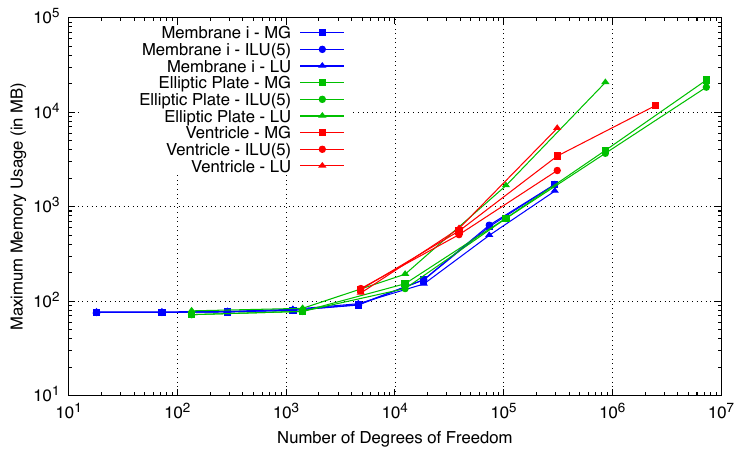}  
   	}
   \caption{The clock times and maximum memory usage for three preconditioning approaches (multigrid, ILU(5), and LU) as a function of degrees of freedom}
   \label{fig:times_memory}
\end{figure}

The preconditioning approaches are compared in three cases: membrane i (2-D, linear elastic), elliptic plate (3-D, linear elastic), and idealised ventricle (3-D, hyperelastic).
 Figure \ref{fig:times_memory} compares the clock times and maximum memory requirements for the three preconditioning approaches as a function of degrees of freedom.
 Examining the clock times, the LU preconditioner is faster for all meshes in the membrane i and idealised ventricle cases, whereas the multigrid approach is faster for all meshes in the elliptic plate cases.
The LU approach was found to be up to $5.4 \times$ faster in the membrane i case than the multigrid approach, and $1.8 \times$ faster in the idealised ventricle case.
In contrast, the multigrid was up to $1.4 \times$ faster than the LU approach in the elliptic plate case.
The ILU(5) approach is the slowest in all cases.
In addition, the ILU(5) approach diverged on the finest idealised ventricle mesh; lower values of fill-in $k$ led to earlier divergence in this case.

Regarding memory usage, below approximately $50$ k degrees of freedom, the memory usage is the same for all three approaches (approximately $80$ MB), corresponding to the solver's minimum memory overhead.
For greater numbers of degrees of freedom, the LU approach requires the most memory in both 3-D cases (idealised ventricle, elliptic plate), while the ILU(5) requires the least.
In contrast, the LU approach requires the least memory for the 2-D cases (membrane i) examined.
The rate of memory increase is seen to be much steeper for the LU approach on the 3-D cases than the multigrid and ILU(5) approaches, with the LU results for the finest meshes not shown as they required more than the maximum available memory (64 GB).
\subsection{Effect of Rhie-Chow Stabilisation} \label{sec:RhieChowResults}
This section highlights the impact of the global scaling factor $\alpha$ in the Rhie-Chow stabilisation (Equation \ref{eq:RhieChow}) on the performance of the proposed Jacobian-free Newton Krylov method.
As described in Section \ref{sec:discretisation}, the stabilisation term is introduced to quell zero-energy modes (oscillations) in the discrete solution, such as checkerboarding.
As the amount of stabilisation increases (increasing $\alpha$), these numerical modes are quelled, and the solution stabilises; however, at some point, further increases in stabilisation reduce the discretisation accuracy due to over-smoothing.
A less obvious consequence of changing the stabilisation magnitude is its effect on the convergence of the linear solver in the Jacobian-free Newton-Krylov solution procedure.
This section examines this effect.

\begin{figure}[htbp]
   \centering
	\subfigure[Clock time]
	{
   		\includegraphics[width=0.48\textwidth]{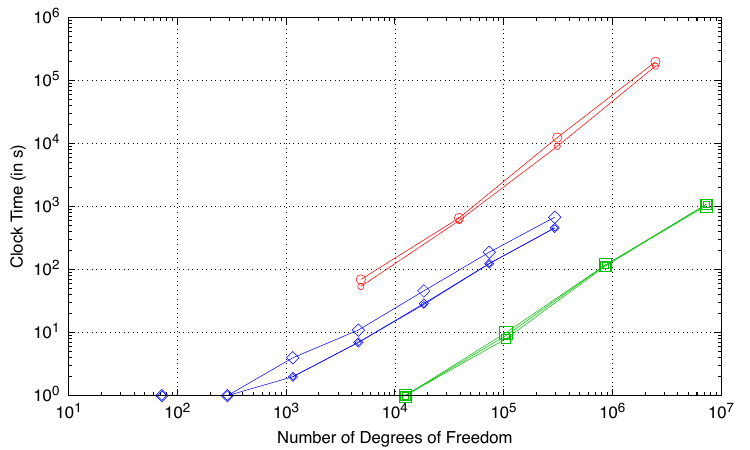}
	}
	\subfigure[Accumulated linear solver iterations]
	{
   		\includegraphics[width=0.48\textwidth]{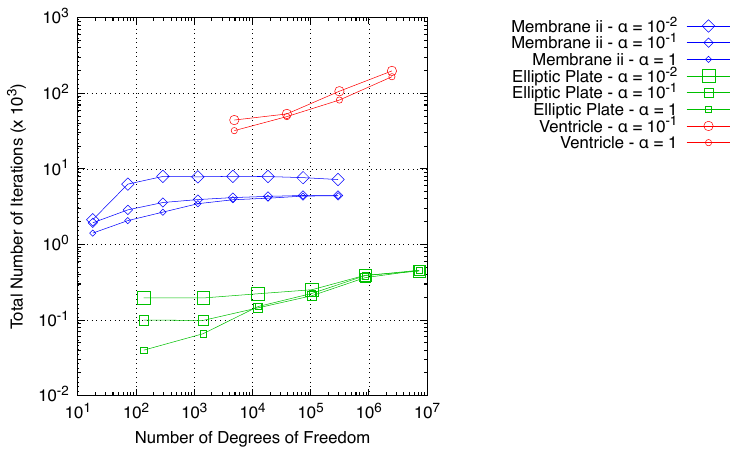}  
   	}
   \caption{The clock times and accumulated linear solver iterations for different values of Rhie-Chow global stabilisation factor $\alpha$}
   \label{fig:times_rhie_chow}
\end{figure}

As in the previous section, three cases are used to highlight the effect: membrane i (2-D, linear elastic), elliptic plate (3-D, linear elastic), and idealised ventricle (3-D, hyperelastic).
The LU preconditioning method is used for the 2-D membrane i case, while the multigrid approach is used for the two 3-D cases.
Figure \ref{fig:times_rhie_chow} presents the execution times and the number of accumulated linear solver iterations for three values of global stabilisation factor ($\alpha = \left[ 0.01, 0.1, 1 \right]$) for several mesh densities.
The memory usage is unaffected by the value of $\alpha$ and is hence not shown.
From Figure \ref{fig:times_rhie_chow}(a), the clock time is seen to increase exponentially (linearly on a log-log plot) for all cases and all values of $\alpha$.
For all meshes in all three cases, the lowest value of global stabilisation factor ($\alpha = 0.01$) -- corresponding to the least amount of stabilisation -- takes the greatest amount of time to converge, while the largest value of global stabilisation factor ($\alpha = 1$) is the fastest to converge.
In the membrane i case, the largest value of  ($\alpha$) is approximately $1.5\times$ faster than the lowest value of $\alpha$.
Nonetheless, the clock time is not a linear function of  $\alpha$, and it can be seen that $\alpha = 0.1$ requires approximately the same amount of time as $\alpha = 1$ in all cases.
In terms of robustness, the idealised ventricle cases failed with $\alpha = 0.01$, suggesting that increasing the stabilisation increases robustness.
From Figure \ref{fig:times_rhie_chow}(b), the accumulated number of linear solver (GMRES) iterations is seen to follow the same trends as the clock times, where $\alpha = 0.01$ requires the greatest number of iterations while $\alpha = 0.1$ and $\alpha = 1$ require approximately the same number of iterations.
In all cases, the number of iterations increases with the number of degrees of freedom.


The choice of the stabilisation scaling parameter $\alpha$ influences both the convergence behaviour and the pre-asymptotic accuracy of the method. Increasing $\alpha$ effectively adds diffusion to the system, lowering the condition number and making the linear systems easier to precondition, thereby reducing the number of iterations required for convergence. However, this benefit comes at the cost of reduced accuracy at coarse resolutions. As shown in Appendix B, the scheme remains formally second-order accurate asymptotically, and in our tests, displacement and stress predictions converged with mesh refinement for all $\alpha$ values.
Nonetheless, the pre-asymptotic behaviour differed between choices of $\alpha$, highlighting a trade-off between solver robustness and accuracy that must be considered in practical applications.

\subsection{Parallelisation}
\label{sec:parallelisation}
In this final analysis, the multi-CPU-core parallel scaling performance of the proposed Jacobian-free Newton-Krylov method is compared with that of the segregated approach.
A \emph{strong} scaling study is performed, where the clock time to solve a fixed-size problem is measured as the number of CPU cores is increased.
The elliptic plate (3-D, linear elastic) mesh with $2\,438\,242$ cells is chosen for the scaling analysis.
Parallelisation uses the standard OpenFOAM domain decomposition approach, where the mesh is decomposed into one sub-domain per CPU core.
In the current work, the Scotch decomposition approach \citep{Pellegrini2012} is employed.
The Jacobian-free Newton Krylov approach uses the GMRES linear solver, where the performance of the multigrid and LU preconditioners is compared.
The segregated approach uses the conjugate gradient solver with the incomplete (zero in-fill -- $k = 0$) Cholesky preconditioner.
The cases are run on the MeluXina high-performance computing system, where each standard computing node contains $2\times$ AMD EPYC Rome 7H12 64c 2.6GHz CPUs with 512 GB of memory.
In the current study, the number of CPU cores is varied from 1 to $1\,048$ ($\left[ 1, 2, 4, 8, 16, 32, 64, 128, 256, 512, 1\,048 \right]$).
In the ideal case, the speedup should double when the number of cores doubles; in reality, inter-CPU communication reduces parallel scaling efficiency below the ideal.

It is noted that the scalability of the Jacobian-free Newton–Krylov procedure is dictated by the scalability of the underlying algorithmic components (GMRES, preconditioner) rather than the Jacobian-free formulation itself. Consequently, the present implementation does not affect this property. The purpose of the scaling study presented here is therefore twofold: first, to confirm that the Jacobian-free Newton–Krylov approach inherits the expected scalability characteristics of its components; and second, to provide a practical comparison of the runtime behaviour of Jacobian-free Newton–Krylov and segregated approaches on representative problem sizes.
Larger problem sizes would be required for a rigorous scaling study; hence, these tests do not establish definitive scalability limits.

%

\begin{figure}[htbp]
	\centering
	\subfigure[Clock times vs number of CPU cores]
	{
		\label{fig:parallel_strong_times}
		\includegraphics[width=0.48\textwidth]{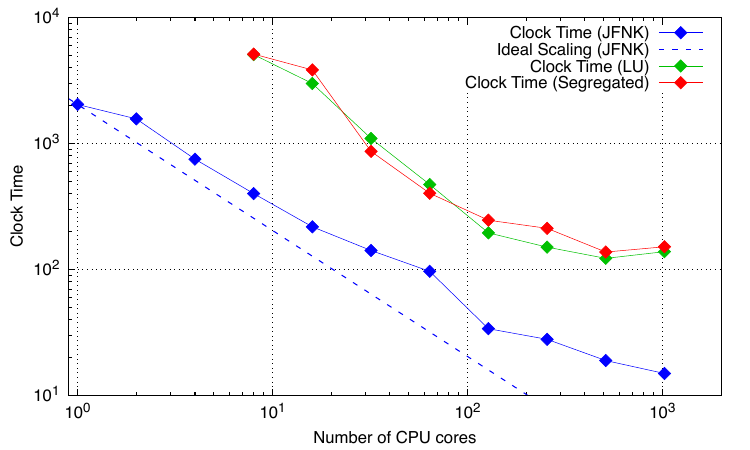}
   	}
	\subfigure[Speedup $S$ vs number of CPU cores]
	{
		\label{fig:parallel_strong_speedup}
		\includegraphics[width=0.48\textwidth]{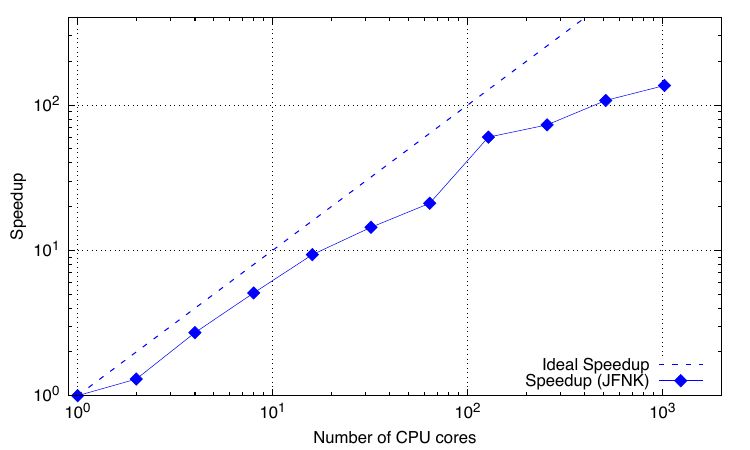}
   	}
	\caption{Strong parallel scaling study comparing the performance of the Jacobian-free Newton-Krylov approach using the multigrid and LU preconditioning strategies, and the segregated solution procedure}
	\label{fig:parallelisation_strong}
\end{figure}

The clock times from the strong scaling study are shown in Figure \ref{fig:parallelisation_strong}(a), where the \emph{ideal} scaling is shown for comparison as a dashed line.
Figure \ref{fig:parallelisation_strong}(b) shows the corresponding speedup, defined as the clock time on one CPU core divided by the clock time on a given number of CPU cores.
The Jacobian-free Newton-Krylov approach with the multigrid preconditioner is found to be the fastest by almost an order of magnitude over the segregated approach and the Jacobian-free Newton-Krylov approach with the LU preconditioner.
The times for the segregated and Jacobian-free Newton-Krylov LU preconditioner approaches using 1, 2, and 4 cores were not recorded due to the excessive times.
For all approaches, the speedup increases approximately at an ideal rate up to 128 cores, after which it continues to increase but at a lower rate.
At 128 cores (corresponding to one full computing node), the number of cells per core is approximately 19 k; increasing the number of cores beyond 128 likely shows a less than ideal scaling for two reasons: (i) the number of cells per core is becoming small relative to the amount of inter-core communication, and (ii) inter-node communication is likely slower than intra-node communication.
The segregated and Jacobian-free Newton-Krylov LU preconditioner approaches produce their fastest predictions at 512 cores, while the Jacobian-free Newton-Krylov multigrid preconditioner approach continues to speed up to the maximum number of cores tested ($1\,048$), albeit $1\,048$ cores is just $2\times$ faster than 128 cores.
At $1\,048$ cores, there are approximately 2.5 k cells per core, and the inter-core communication is expected to be significant.

As a final observation, the effect of hardware can be highlighted by comparing the time required for 1 core when using the multigrid preconditioner: the case required $2\,048$ s on the MeluXina system (EPYC Rome CPU), while it required $988$ s on a Mac Studio (M1 Ultra CPU) system as presented in Section \ref{sec:resource_requirements} -- over twice as fast.
The performance difference can primarily be attributed to the unified memory architecture of the M1 Ultra, which provides significantly greater memory bandwidth than the DDR4 memory used in the MeluXina compute nodes.
\section{Conclusions} \label{sec:conclusion}
A Jacobian-free Newton-Krylov solution algorithm has been proposed for solid-mechanics problems discretised using the cell-centred finite-volume method.
A compact-stencil discretisation of the diffusion term is proposed as the preconditioner, enabling a straightforward extension of existing segregated solution frameworks.
The key conclusions of the work are:
\begin{itemize}
	\item \textbf{Efficiency of Jacobian-free Newton-Krylov approach}: The proposed Jacobian-free Newton-Krylov solution algorithm has been shown to be faster than a conventional segregated solution algorithm for all linear and nonlinear elastic test cases examined. In particular, speedups of one order of magnitude were seen in many cases, with a maximum speedup of 341 in the 2-D linear elastic Cook's membrane case.

	\item \textbf{Additional memory overhead}: The Jacobian-free Newton-Krylov approach has been shown to have approximately the same memory requirements as a segregated solution algorithm for less than $10^4$ degrees of freedom, but generally increases beyond this.
	The greatest memory increase was 3.6 and occurred in the spherical cavity case. Nonetheless, these differences are likely due to the different linear solvers and preconditioners used in the Jacobian-free Newton-Krylov approach (GMRESS with multigrid or LU decomposition) compared to the segregated approach (conjugate gradient with zero-infill Cholesky decomposition).
	
	\item \textbf{Applicability to existing segregated frameworks}: By employing a compact-stencil diffusion approximation of the Jacobian as the preconditioner, the proposed Jacobian-free Newton-Krylov approach can be integrated into existing segregated finite volume frameworks with minimal modifications to the existing code base, in particular if existing publicly available Jacobian-free Newton-Krylov solvers (e.g., PETSc) are used.


	\item \textbf{Choice of preconditioning approach}: It has been shown that the LU direct solver preconditioning approach is faster than multigrid and ILU($k$) approaches for 2-D cases and moderately-sized 3-D problems; however, for larger 3-D problems, the multigrid approach is the fastest while also requiring less memory than the direct LU approach. In addition, the multigrid approach shows approximate ideal scaling in a parallel strong-scaling study.
	
	\item \textbf{Rhie-Chow stabilisation}: The magnitude of the Rhie-Chow stabilisation term is shown to affect the speed and robustness of the Jacobian-free Newton-Krylov approach, where $\alpha = 0.1$ and $1$ are seen to outperform $\alpha = 0.01$.

	\item \textbf{Open access and extensibility}: The implementation is made publicly available within the solids4foam toolbox for OpenFOAM to encourage implementation critique, community adoption, and comparative studies, contributing to advancing finite volume solid mechanics simulations.
\end{itemize}


The presented study highlights the potential of the Jacobian-free Newton-Krylov method in finite-volume solid mechanics simulations, yet several areas for future exploration remain:
\begin{enumerate}
	\item{Enhanced preconditioning for plasticity}: While the proposed \emph{elastic} compact-stencil preconditioning matrix has been shown to perform well in linear and nonlinear elastic scenarios, it fails in cases exhibiting plasticity, where a segregated approach succeeds. Future studies will explore modifications to the compact preconditioning matrix to extend the applicability of the proposed approach to both small- and large-strain elastoplasticity.

	\item{Incorporation of non-orthogonality in the preconditioner}: The current preconditioning matrix does not account for non-orthogonal contributions in distorted meshes. Including these effects in future work could further enhance robustness, particularly for complex geometries.

	\item{Comparison with full Jacobian approaches}: A comparison with a fully coupled Jacobian-based method would provide greater insights into the trade-offs between computational cost, memory usage, and convergence performance. Such a study could highlight the specific advantages of the Jacobian-free Newton-Krylov approach for large-scale, nonlinear problems.

	\item{Globalisation strategies}: The robustness and efficiency of the proposed Jacobian-free Newton-Krylov approach could be improved for nonlinear problems by using a segregated approach as an initialisation phase for each loading step. In its current form, the segregated approach would not be generally suitable as it is less robust on several nonlinear hyperelastic problems; however, modifications such as introducing pseudo-time-stepping (introducing a pseudo-time damping term) may enhance its robustness and suitability for these scenarios.

	\item{Higher-order discretisations}: Future work will explore the application of the proposed linear elastic compact-stencil Jacobian-free Newton-Krylov approach to higher-order finite volume discretisations. Using the same compact-stencil preconditioning matrix across these discretisations could yield efficient, accurate solutions for problems requiring higher fidelity.

	\item{Development of a monolithic fluid-solid interaction framework}: Extending the Jacobian-free Newton-Krylov method to a fully monolithic fluid-solid interaction framework presents a promising avenue. In this approach, both fluid and solid domains would be solved simultaneously using a unified Jacobian-free Newton-Krylov method, potentially improving stability and performance in highly nonlinear coupled problems.
\end{enumerate}

\backmatter

\bmhead{Data Availability}
The codes presented are publicly available at \url{https://github.com/solids4foam/solids4foam} on the \texttt{feature-petsc-snes} branch, and the cases and plotting scripts are available at \url{https://github.com/solids4foam/solid-benchmarks}.

\bmhead{Declaration of generative AI and AI-assisted technologies in the writing process}
During the preparation of this work, the authors used ChatGPT and Grammarly as writing assistants.
The authors reviewed and edited the content as needed and take full responsibility for the publication's content.

\bmhead{Acknowledgments}
Technical reviews and insightful comments from Hiroaki Nishikawa of the National Institute of Aerospace, Hampton, VA, USA, are greatly appreciated.
This project has received funding from the European Research Council (ERC) under the European Union’s Horizon 2020 research and innovation programme (Grant Agreement No. 101088740).
Financial support is gratefully acknowledged from the Irish Research Council
through the Laureate programme, grant number IRCLA/2017/45, from Bekaert through the University Technology Centre (UTC phases I and II) at UCD
(www.ucd.ie/bekaert), from I-Form, funded by Research Ireland (formerly Science Foundation Ireland)
Grant Numbers {16/RC/3872} and {21/RC/10295\_P2}, co-funded under European Regional Development Fund and by I-Form industry partners, and from NexSys, Grant Number 21/SPP/3756.
Additionally, the authors wish to acknowledge the DJEI/DES/SFI/HEA Irish Centre for High-End Computing (ICHEC) for the provision of computational facilities and support (www.ichec.ie), and part of this work has been carried out using the UCD ResearchIT Sonic cluster, which was funded by UCD IT Services and the UCD Research Office.

\newpage

\begin{appendices}

\section{Total and Updated Lagrangian Formulations of the Conservation of Linear Momentum}
\label{app:TL_UL}
More generally,  linear momentum conservation can be expressed in a nonlinear geometry form, which is suitable for finite strains.
Two equivalent nonlinear geometry forms are common: the \emph{total} Lagrangian form:
\begin{eqnarray} \label{eqn:momentum_TL}
    \int_{\Omega_o} \rho_o \frac{\partial^2 \bb{u} }{\partial t^2} d\Omega_o
    =
    \oint_{\Gamma_o} \left( J \bb{F}^{-\text{T}} \cdot \bb{n}_o \right) \cdot \bb{\sigma} \ d\Gamma_o
    + \int_{\Omega_o}  \bb{f}_b \, d\Omega_o
\end{eqnarray}
and the \emph{updated} Lagrangian form:
\begin{eqnarray} \label{eqn:momentum_UL}
    \int_{\Omega_u} \frac{\partial }{\partial t} \left( \rho_u \frac{\partial \bb{u} }{\partial t} \right) d\Omega_u
    = \oint_{\Gamma_u}(j\bb{f}^{-\text{T}}\cdot{\bb{n}_u)\cdot \bb{\sigma}}\ d\Gamma_u
    + \int_{\Omega_u}  \bb{f}_b \, d\Omega_u
\end{eqnarray}
where subscript $o$ indicates quantities in the initial reference configuration, and subscript $u$ indicates quantities in the updated configuration.
The true (Cauchy) stress tensor is indicated by $\bb{\sigma}$.
The deformation gradient is defined as $\bb{F} = \textbf{I} + (\bb{\nabla} \bb{u})^{\text{T}}$ and its determinant as $J = \text{det}(\bb{F})$.
Similarly, the \emph{relative} deformation gradient is given in terms of the displacement \emph{increment} as $\bb{f}=\textbf{I} + \left[\bb{\nabla}(\Delta \bb{u}) \right]^{\text{T}}$ and its determinant as $j = \text{det}(\bb{f})$.
The displacement increment is the change in displacement between the current and previous time steps when the time interval is discretised into a finite number of steps.

\section{Mechanical Laws}
\label{app:mechLaws}

\subsection{Linear Elasticity}
The definition of engineering stress $\bb{\sigma}_s$ for linear elasticity can be given as
\begin{eqnarray}  \label{eq:linearElastic}
	\bb{\sigma}_s
	&=& 2 \mu \bb{\varepsilon} + \lambda \, \text{tr} \left( \bb{\varepsilon} \right) \textbf{I} \notag \\
	&=& \mu \bb{\nabla} \bb{u} + \mu \left( \bb{\nabla} \bb{u}\right)^{\text{T}} + \lambda \left(\bb{\nabla} \cdot \bb{u} \right) \textbf{I}
\end{eqnarray}
where $\lambda$ is the first Lam\'{e} parameter, and $\mu$ is the second Lam\'{e} parameter, synonymous with the shear modulus.
The Lam\'{e} parameters can be expressed in terms of the Young's modulus ($E$) and Poisson's ratio $\nu$ as
\begin{eqnarray}
	\mu = \frac{E}{2(1 + \nu)}, \quad \lambda = \frac{E \nu}{(1+\nu)(1 - 2\nu)}
\end{eqnarray}

\subsection{St.\,Venant-Kirchoff Hyperelasticity}
The St.\,Venant-Kirchoff model defines the second Piola–Kirchhoff stress $\textbf{S}$ as
\begin{eqnarray}
	\bb{S} &=& 2 \mu \bb{E} + \lambda \, \text{tr} \left( \bb{E} \right) \textbf{I}
\end{eqnarray}
where, as before, $\lambda$ is the first Lam\'{e} parameter, and $\mu$ is the second Lam\'{e} parameter.
The Lagrangian Green strain $\textbf{E}$ is defined as
\begin{eqnarray}
	\bb{E} &=& \frac{1}{2} \left( \bb{\nabla} \bb{u} + \bb{\nabla} \bb{u}^{\text{T}} + \bb{\nabla} \bb{u} \cdot \bb{\nabla} \bb{u}^{\text{T}}  \right)
\end{eqnarray}

The true stress can be calculated from the second Piola–Kirchhoff stress as
\begin{eqnarray} \label{eq:S2sigma}
	\bb{\sigma} &=& \frac{1}{J} \bb{F} \cdot \bb{S} \cdot \bb{F}^{\text{T}}
\end{eqnarray}

\subsection{Neo-Hookean Hyperelasticity}
The definition of true (Cauchy) stress $\bb{\sigma}$ for neo-Hookean hyperelasticity can be given as
\begin{eqnarray} \label{eq:neoHook}
	\bb{\sigma}
	&=& \frac{\mu}{J} \, \text{dev} \left( \bar{\bb{b}} \right) + \frac{\kappa}{2} \frac{J^2 - 1}{J} \textbf{I}
\end{eqnarray}
where, once again, $\mu$ is the shear modulus, and $\kappa$ is the bulk modulus.
The bulk modulus can be expressed in terms of the Young's modulus ($E$) and Poisson's ratio $\nu$ as
\begin{eqnarray}
	\kappa = \frac{E}{3(1 - 2\nu)}
\end{eqnarray}
The volume-preserving component of the elastic left Cauchy--Green deformation tensor, $\boldsymbol{b}$, is given as
\begin{eqnarray}
	\bar{\bb{b}} = J^{-2/3} \bb{b} = J^{-2/3} \bb{F} \cdot \bb{F}^{\text{T}}
\end{eqnarray}
In the limit of small deformations $\lVert \nabla \mathbf{u} \rVert  \ll 1$, neo-Hookean hyperelasticity (Equation \ref{eq:neoHook}) reduces to linear elasticity (Equation \ref{eq:linearElastic}).

\subsection{Guccione Hyperelasticity}
The \citet{Guccione1995} hyperelastic law defines the second Piola-Kirchhoff stress as
\begin{eqnarray}
	\boldsymbol{S}
		&=& \frac{\partial Q}{\partial \boldsymbol{E} } \left( \frac{C}{2} \right) e^Q + \frac{\kappa}{2} \frac{J^2 - 1}{J} \textbf{I} 
\end{eqnarray}
where
\begin{eqnarray}
	Q(I_1, I_2, I_4, I_5) &=& c_t I_1^2 - 2 c_t I_2 + (c_f - 2c_{fs} + c_t) I_4^2   + 2(c_{fs} - c_t) I_5 \\
	\frac{\partial Q}{\partial \boldsymbol{E} }
		&=&  2 c_t \boldsymbol{E} + 2(c_f - 2 c_{fs} + c_t)I_4 ( \boldsymbol{f_0} \otimes \boldsymbol{f_0} ) \notag \\
		&&+ 2(c_{fs} - c_t)\left[\boldsymbol{E} \cdot  (\boldsymbol{f_0} \otimes \boldsymbol{f_0}) + ( \boldsymbol{f_0} \otimes \boldsymbol{f_0} ) \cdot \boldsymbol{E} \right]
\end{eqnarray}
The scalars $C$, $c_f$, $c_{fs}$, and $c_t$ are material parameters and invariants of the Green strain, $\boldsymbol{E} = \boldsymbol{F}^{\text{T}} \cdot \boldsymbol{F}$, are defined as
\begin{eqnarray}
	I_1 = \text{tr}(\boldsymbol{E}), \quad
	I_2 =  \frac{1}{2} \left[ \text{tr}^2(\boldsymbol{E}) - \text{tr}(\boldsymbol{E} \cdot \boldsymbol{E}) \right], \notag \\
	I_4 = \boldsymbol{E}  : \left( \boldsymbol{f_0} \otimes \boldsymbol{f_0} \right), \quad 
	I_5 = \left( \boldsymbol{E} \cdot \boldsymbol{E} \right): \left( \boldsymbol{f_0} \otimes \boldsymbol{f_0} \right)
\end{eqnarray}
with $\boldsymbol{f_0}$ representing the unit fibre directions in the initial configuration.

Equation \ref{eq:S2sigma} is used to convert the second Piola-Kirchhoff stress to the true stress.

\subsection{Neo-Hookean $J_2$ Hyperelastoplasticity}
For neo-Hookean $J_2$ hyperelastoplasticity, the expression for the true (Cauchy) stress $\bb{\sigma}$ takes the same form as Equation \ref{eq:neoHook}, except $\boldsymbol{b}$ is replaced by its elastic component $\boldsymbol{b}_e$.
Determination of $\boldsymbol{b}_e$ employs the definition of $J_2$ (Mises) plasticity in terms of a yield function, flow rule, Kuhn–Tucker loading/unloading conditions, and the consistency condition.
The stress calculation procedure (radial return algorithm) is described by \citet{Simo1998} (Box 9.1, page 319).

\section{Truncation Error Analysis of the Rhie-Chow Stabilisation Term}
\label{app:RhieChow}
On a 1-D uniform mesh with spacing $\Delta x$ and unity areas, the Rhie-Chow term stabilisation term (Equation \ref{eq:RhieChow}) for an internal cell $P$ (no boundary faces) becomes
\begin{eqnarray}
\mathcal{D}_P^{\text {Rhie-Chow}}
	&=&
	\sum_{f_i \in \mathcal{F}^{\text{int}}_p} \alpha \bar{K}
	\left[
	\left|\bb{\Delta}_{f_i} \right| \frac{ \bb{u}_{N_{f_i}} - \bb{u}_P}{\left|\bb{d}_f\right|}
	- \bb{\Delta}_{f_i} \cdot \left(\bb{\nabla} \bb{u} \right)_f
	\right]
	\left|\bb{\Gamma}_{{f_i}}\right| \notag \\
	&=& 
	\alpha \bar{K} \left[ \frac{ \bb{u}_{E} - \bb{u}_P}{\Delta x} - \left(\bb{\nabla} \bb{u} \right)_e \right]
	+ \alpha \bar{K} \left[ \frac{ \bb{u}_{W} - \bb{u}_P}{\Delta x} + \left(\bb{\nabla} \bb{u} \right)_w \right]
	\notag \\
\end{eqnarray}
where $E$ and $W$ indicate the east and west neighbour cell centre values, $EE$ and $WW$ are the far east and west neighbour cell centre values, and $e$ and $w$ indicate east and west face values;
$\alpha$ and $\bar{K}$ are assumed uniform, and $\left|\bb{\Gamma}_{{f_i}}\right|$ is assumed equal to unity.
The face gradients are calculated as
\begin{eqnarray}
	\left(\bb{\nabla} \bb{u} \right)_e = \frac{1}{2} \left[ \left(\bb{\nabla} \bb{u} \right)_P + \left(\bb{\nabla} \bb{u} \right)_E \right]  \notag \\
	\left(\bb{\nabla} \bb{u} \right)_w = \frac{1}{2} \left[ \left(\bb{\nabla} \bb{u} \right)_W + \left(\bb{\nabla} \bb{u} \right)_P \right]
\end{eqnarray}
where the cell centre gradients are calculated as
\begin{eqnarray}
	\left(\bb{\nabla} \bb{u} \right)_W = \frac{\bb{u}_{P} - \bb{u}_{WW}}{2\Delta x} \notag \\
	\left(\bb{\nabla} \bb{u} \right)_P = \frac{\bb{u}_{E} - \bb{u}_{W}}{2\Delta x} \notag \\
	\left(\bb{\nabla} \bb{u} \right)_E = \frac{\bb{u}_{EE} - \bb{u}_{P}}{2\Delta x}
\end{eqnarray}

The final Rhie-Chow term becomes
\begin{eqnarray} \label{eq:Rhie_Chow_1D}
	\mathcal{D}_P^{\text {Rhie-Chow}}
	&=& 
	\alpha \bar{K}
		\left\{
		\frac{ \bb{u}_{E} - \bb{u}_P}{\Delta x}
		- \frac{1}{2} \left[ \left(\bb{\nabla} \bb{u} \right)_P + \left(\bb{\nabla} \bb{u} \right)_E \right]
	+ \frac{ \bb{u}_{W} - \bb{u}_P}{\Delta x}
		+ \frac{1}{2} \left[ \left(\bb{\nabla} \bb{u} \right)_W + \left(\bb{\nabla} \bb{u} \right)_P \right]
		\right\}
	\notag \\
	&=& 
	\alpha \bar{K}
		\left\{
		\frac{ \bb{u}_{E} - \bb{u}_P}{\Delta x}
		- \frac{1}{2} \left(\bb{\nabla} \bb{u} \right)_E
	+ \frac{ \bb{u}_{W} - \bb{u}_P}{\Delta x}
		+ \frac{1}{2} \left(\bb{\nabla} \bb{u} \right)_W 
		\right\}
	\notag \\
	&=& 
	\alpha \bar{K}
		\left[
		\frac{ \bb{u}_{E} - \bb{u}_P}{\Delta x}
		- \frac{\bb{u}_{EE} - \bb{u}_{P}}{4\Delta x}
	+ \frac{ \bb{u}_{W} - \bb{u}_P}{\Delta x}
		+ \frac{\bb{u}_{P} - \bb{u}_{WW}}{4\Delta x} 
		\right]
	\notag \\
	&=& 
	\frac{\alpha \bar{K}}{4 \Delta x}	\left[- \bb{u}_{WW} + 4\bb{u}_W - 6\bb{u}_P + 4\bb{u}_{E} - \bb{u}_{EE} \right]
\end{eqnarray}

A local truncation analysis can be performed using the following truncated Taylor series about the true solution ($\bb{U}$) and its gradients ($\left( \bb{\nabla} \bb{U} \right)_P$, $\left( \bb{\nabla}^2 \bb{U} \right)_P$, $...$) at $P$:
\begin{eqnarray} \label{eq:taylor_series_1d}
	\bb{u}_{WW} &=& \bb{U}_P - 2\Delta x \left( \bb{\nabla} \bb{U} \right)_P + \frac{4\Delta x^2}{2!} \left( \bb{\nabla}^2 \bb{U} \right)_P
		- \frac{8\Delta x^3}{3!} \left( \bb{\nabla}^3 \bb{U} \right)_P + O(\Delta x^4) \notag \\
	\bb{u}_W &=& \bb{U}_P - \Delta x \left( \bb{\nabla} \bb{U} \right)_P + \frac{\Delta x^2}{2!} \left( \bb{\nabla}^2 \bb{U} \right)_P
		- \frac{\Delta x^3}{3!} \left( \bb{\nabla}^3 \bb{U} \right)_P + O(\Delta x^4) \notag \\
	\bb{u}_P &=& \bb{U}_P \notag \\
	\bb{u}_E &=& \bb{U}_P + \Delta x \left( \bb{\nabla} \bb{U} \right)_P + \frac{\Delta x^2}{2!} \left( \bb{\nabla}^2 \bb{U} \right)_P 
		+ \frac{\Delta x^3}{3!} \left( \bb{\nabla}^3 \bb{U} \right)_P + O(\Delta x^4) \notag \\
	\bb{u}_{EE} &=& \bb{U}_P + 2\Delta x \left( \bb{\nabla} \bb{U} \right)_P + \frac{4\Delta x^2}{2!} \left( \bb{\nabla}^2 \bb{U} \right)_P 
		+ \frac{8\Delta x^3}{3!} \left( \bb{\nabla}^3 \bb{U} \right)_P + O(\Delta x^4)
\end{eqnarray}
where $O(\Delta x^4)$ are higher-order terms with a leading term proportional to $\Delta x^4$.

Substituting the expressions above (Equations \ref{eq:taylor_series_1d}) for the true solution into Equation \ref{eq:Rhie_Chow_1D} results in
\begin{eqnarray}
	\mathcal{D}_P^{\text {Rhie-Chow}}
	&=& 	\frac{\alpha \bar{K}}{4 \Delta x}	\left[- \bb{u}_{WW} + 4\bb{u}_W - 6\bb{u}_P + 4\bb{u}_{E} - \bb{u}_{EE} \right] \notag \\
	&=& 	\frac{\alpha \bar{K}}{4 \Delta x}
	\Bigg[
	- \bb{U}_P + 2\Delta x \left( \bb{\nabla} \bb{U} \right)_P - 4\frac{\Delta x^2}{2!} \left( \bb{\nabla}^2 \bb{U} \right)_P 
		+ 8\frac{\Delta x^3}{3!} \left( \bb{\nabla}^3 \bb{U} \right)_P \notag \\
	&&\quad\quad + 4\bb{U}_P - 4\Delta x \left( \bb{\nabla} \bb{U} \right)_P + 4\frac{\Delta x^2}{2!} \left( \bb{\nabla}^2 \bb{U} \right)_P
		- 4\frac{\Delta x^3}{3!} \left( \bb{\nabla}^3 \bb{U} \right)_P \notag \\
	&&\quad\quad - 6\bb{U}_P \notag \\
	&&\quad\quad + 4\bb{U}_P + 4\Delta x \left( \bb{\nabla} \bb{U} \right)_P + 4\frac{\Delta x^2}{2!} \left( \bb{\nabla}^2 \bb{U} \right)_P
		+ 4 \frac{\Delta x^3}{3!} \left( \bb{\nabla}^3 \bb{U} \right)_P \notag \\
	&&\quad\quad - \bb{U}_P - 2\Delta x \left( \bb{\nabla} \bb{U} \right)_P - 4\frac{\Delta x^2}{2!} \left( \bb{\nabla}^2 \bb{U} \right)_P 
		+ 8\frac{\Delta x^3}{3!} \left( \bb{\nabla}^3 \bb{U} \right)_P + O(\Delta x^4) \Bigg] \notag \\
	&=& 	\frac{\alpha \bar{K}}{4 \Delta x}	\left[ 16\frac{\Delta x^3}{3!} \left( \bb{\nabla}^3 \bb{U} \right)_P  + O(\Delta x^4) \right] \notag \\
	&=& 	\frac{2}{3} \alpha \bar{K}  \Delta x^2 \left( \bb{\nabla}^3 \bb{U} \right)_P  + O(\Delta x^3)
\end{eqnarray}
showing the leading truncation error term to reduce at a second order rate as $\Delta x$ is reduced.

\section{Expressing the Rhie-Chow Stabilisation as a 'Jump' Term}
\label{app:RhieChowJump}
Within the square braces on the right-hand side of Equation \ref{eq:RhieChow}, the first terms represent a compact stencil (two-node) approximation of the face normal gradient, while the second terms represent a larger stencil approximation.
These two terms cancel out in the limit of mesh refinement (or if the solution varies linearly); otherwise, they produce a stabilisation effect that tends to smooth the solution fields.
The magnitude of the stabilisation is shown in Appendix \ref{app:RhieChow} to reduce at a second-order rate (not a third-order rate as stated previously \citep{Demirdzic1995}), and hence does not affect the overall scheme's second-order accuracy.
It is informative to note that the first term in the square bracket in Equation \ref{eq:RhieChow} can be expressed, assuming $(\bb{d}_{f_i} \cdot \bb{n}_{f_i}) > 0$, as
\begin{eqnarray} \label{eq:RhieChow2}
	\alpha \left[ \left|\bb{\Delta}_{f_i} \right| \frac{ \bb{u}_{N_{f_i}} - \bb{u}_P}{\left|\bb{d}_{f_i}\right|}	-  \bb{\Delta}_{f_i} \cdot \left(\bb{\nabla} \bb{u} \right)_{f_i} \right]
	&=&
	\alpha \left| \frac{\bb{d}_{f_i}}{\bb{d}_{f_i} \cdot \bb{n}_{f_i}} \right| \frac{ \bb{u}_{N_{f_i}} - \bb{u}_P}{\left|\bb{d}_{f_i}\right|}
	-  \alpha \frac{\bb{d}_{f_i}}{\bb{d}_{f_i} \cdot \bb{n}_{f_i}} \cdot \left(\bb{\nabla} \bb{u} \right)_{f_i} \notag \\
	&=&
	\frac{\alpha}{\bb{d}_{f_i} \cdot \bb{n}_{f_i}} \left[ \left( \bb{u}_{N_{f_i}} - \bb{u}_P \right) -  \bb{d}_{f_i} \cdot \left(\bb{\nabla} \bb{u} \right)_{f_i} \right]
\end{eqnarray}
Noting that $\left(\bb{\nabla} \bb{u} \right)_{f_i} = \left[ \left(\bb{\nabla} \bb{u} \right)_{N_{f_i}} + \left(\bb{\nabla} \bb{u} \right)_P \right]/2$, Equation \ref{eq:RhieChow2} can be expressed as
\begin{eqnarray}
	\alpha \left[ \left|\bb{\Delta}_{f_i} \right| \frac{ \bb{u}_{N_{f_i}} - \bb{u}_P}{\left|\bb{d}_{f_i}\right|}	-  \bb{\Delta}_{f_i} \cdot \left(\bb{\nabla} \bb{u} \right)_{f_i} \right]
	&=&
	\frac{\alpha}{\bb{d}_{f_i} \cdot \bb{n}_{f_i}} \left( \bb{u}_{N_{f_i}}^* - \bb{u}_P^* \right)
\end{eqnarray}
where
\begin{eqnarray}
	\bb{u}_P^* &=& \bb{u}_P + (\bb{d}_{f_i}/2) \cdot \left(\bb{\nabla} \bb{u} \right)_P  \\
	\bb{u}_{N_{f_i}}^* &=& \bb{u}_{N_{f_i}} - (\bb{d}_{f_i}/2) \cdot \left(\bb{\nabla} \bb{u} \right)_{N_{f_i}}
\end{eqnarray}

So, this Rhie-Chow stabilisation can take the form of a \emph{jump} term; that is, it is a function of the jump in solution between neighbouring cells.
The two states, $\bb{u}_P^*$ and $\bb{u}_{N_{f_i}}^*$, are evaluated halfway between $P$ and $N_{f_i}$, and it may not be at a face centre; however, as noted by \citet{Nishikawa2010}, this will not affect the second-order accuracy of the method; consequently, the user is free to choose the value of $\alpha$ for stabilisation and accuracy purposes.

\section{Body Force for the Method of Manufactured Solutions Case} \label{app:mms}
The body force for the manufactured solution case is \citep{Mazzanti2024}:
\begin{align}
\bb{f}_b = 
    \begin{pmatrix}
    \lambda
    \left[
        8 a_y \pi^2 \cos(4\pi x) \cos(2\pi y) \sin(\pi z) \right. \\
        \quad + 4 a_z \pi^2 \cos(4\pi x) \cos(\pi z) \sin(2\pi y) \\
        \quad \left. - 16 a_x \pi^2 \sin(4\pi x) \sin(2\pi y) \sin(\pi z)
    \right] \\
    + \mu
    \left[
        8 a_y \pi^2 \cos(4\pi x) \cos(2\pi y) \sin(\pi z) \right. \\
        \quad + 4 a_z \pi^2 \cos(4\pi x) \cos(\pi z) \sin(2\pi y) \\
        \quad \left. - 5 a_x \pi^2 \sin(4\pi x) \sin(2\pi y) \sin(\pi z)
    \right] \\
    - 32 a_x \mu_ \pi^2 \sin(4\pi x) \sin(2\pi y) \sin(\pi z) \\
    \\
    \lambda
    \left[
        8 a_x \pi^2 \cos(4\pi x) \cos(2\pi y) \sin(\pi z) \right. \\
        \quad + 2 a_z \pi^2 \cos(2\pi y) \cos(\pi z) \sin(4\pi x) \\
        \quad \left. - 4 a_y \pi^2 \sin(4\pi x) \sin(2\pi y) \sin(\pi z)
    \right] \\
    + \mu
    \left[
        8 a_x \pi^2 \cos(4\pi x) \cos(2\pi y) \sin(\pi z) \right. \\
        \quad + 2 a_z \pi^2 \cos(2\pi y) \cos(\pi z) \sin(4\pi x) \\
        \quad \left. - 17 a_y \pi^2 \sin(4\pi x) \sin(2\pi y) \sin(\pi z)
    \right] \\
    - 8 a_y \mu_ \pi^2 \sin(4\pi x) \sin(2\pi y) \sin(\pi z) \\
    \\
    \lambda
    \left[
        4 a_x \pi^2 \cos(4\pi x) \cos(\pi z) \sin(2\pi y) \right. \\
        \quad + 2 a_y \pi^2 \cos(2\pi y) \cos(\pi z) \sin(4\pi x) \\
        \quad \left. - a_z \pi^2 \sin(4\pi x) \sin(2\pi y) \sin(\pi z)
    \right] \\
    + \mu
    \left[
        4 a_x \pi^2 \cos(4\pi x) \cos(\pi z) \sin(2\pi y) \right. \\
        \quad + 2 a_y \pi^2 \cos(2\pi y) \cos(\pi z) \sin(4\pi x) \\
        \quad \left. - 20 a_z \pi^2 \sin(4\pi x) \sin(2\pi y) \sin(\pi z)
    \right] \\
    - 2 a_z \mu_ \pi^2 \sin(4\pi x) \sin(2\pi y) \sin(\pi z)
    \end{pmatrix}
\end{align}
\section{Meshes Used with the Method of Manufactured Solutions}
\label{app:meshes}
\begin{figure}[htbp]
	\centering
	\subfigure[Regular tetrahedral mesh with $4\,374$ cells]
	{
   		\includegraphics[width=0.45\textwidth]{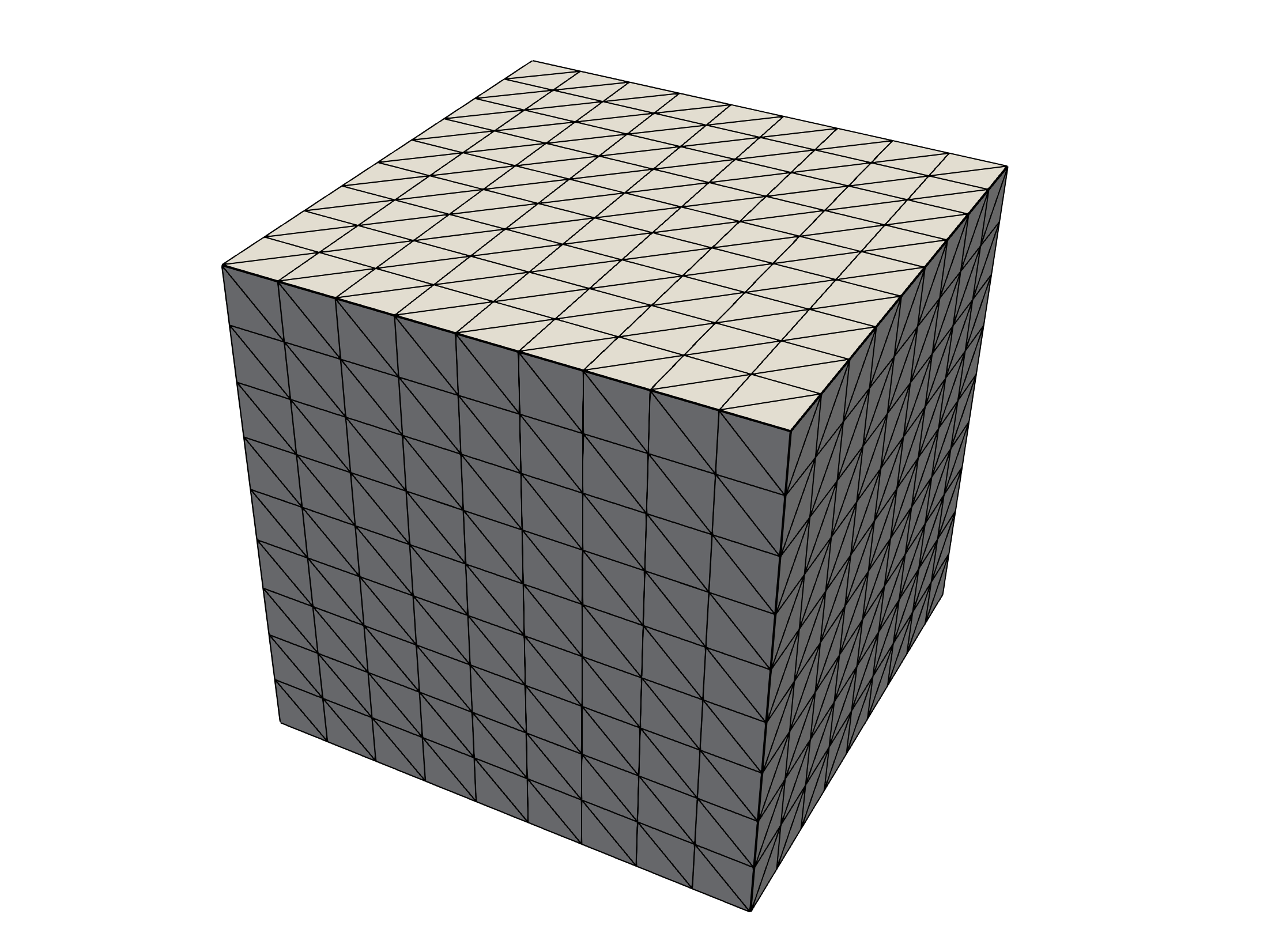}  
   	}
	\subfigure[Regular polyhedral mesh with $1\,000$ cells]
	{
   		\includegraphics[width=0.45\textwidth]{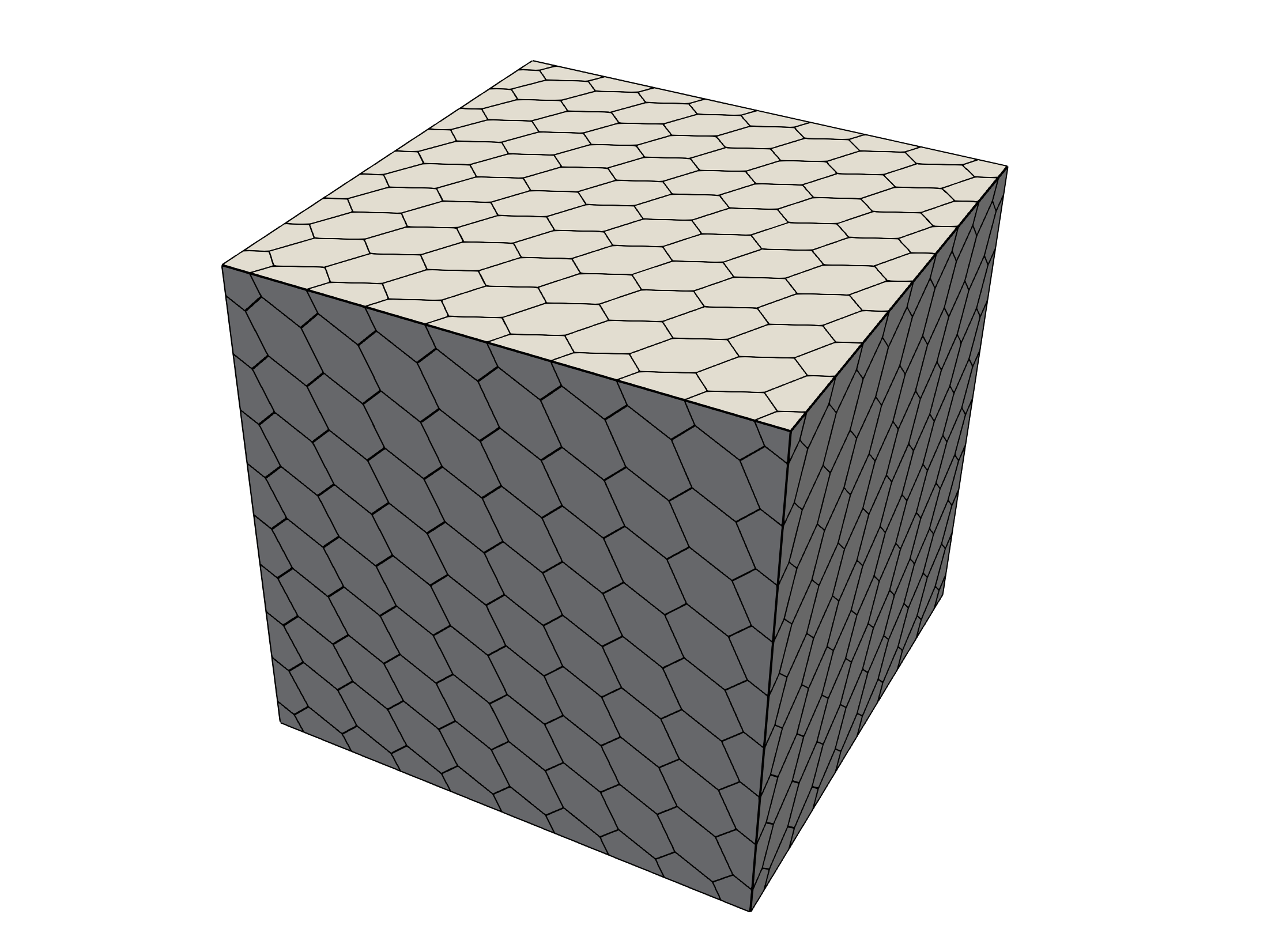}  
   	}
	\subfigure[Regular hexahedral mesh with $1\,000$ cells]
	{
   		\includegraphics[width=0.45\textwidth]{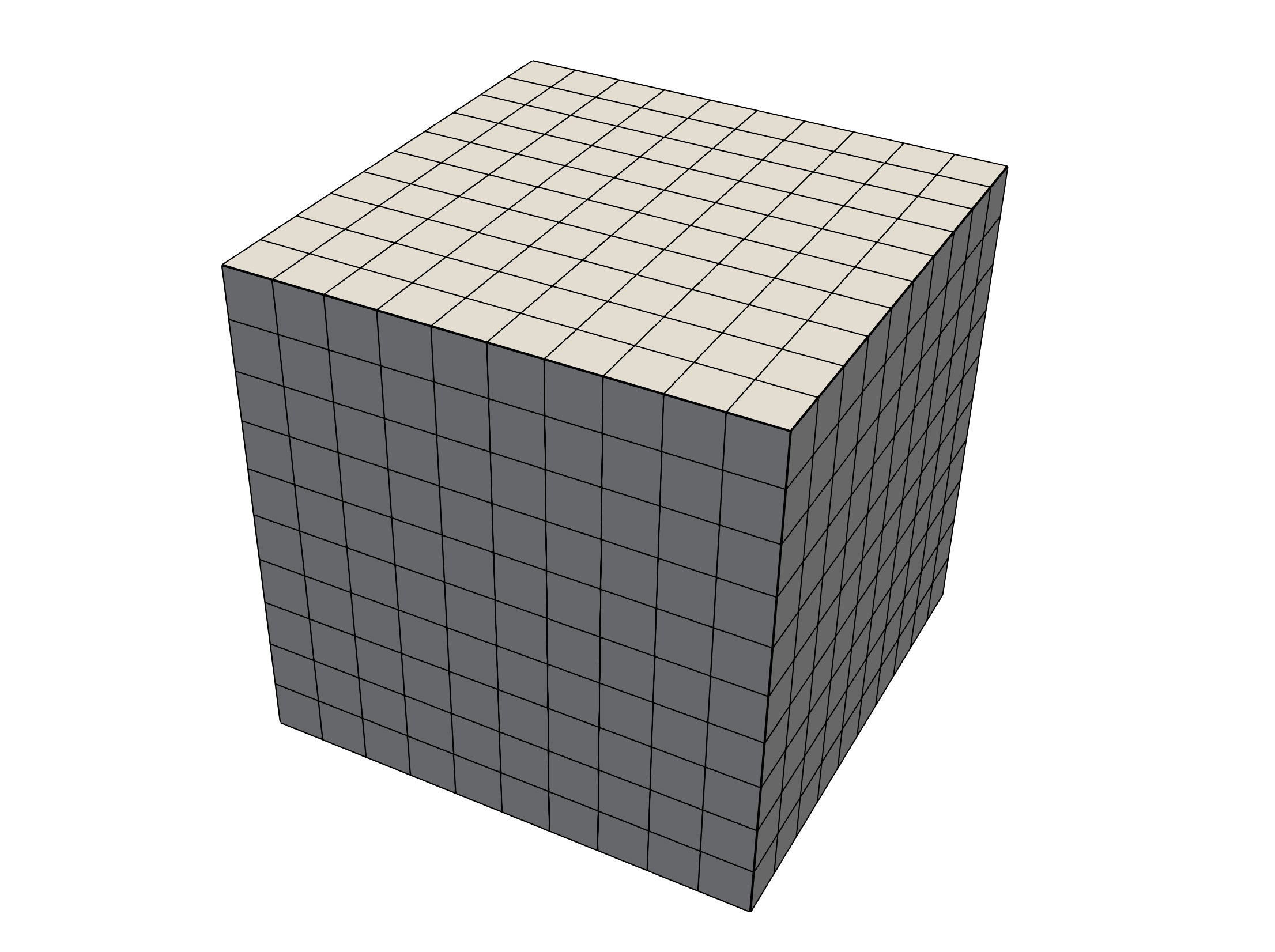}  
   	}
	\subfigure[Distorted hexahedral mesh with $1\,000$ cells]
	{
   		\includegraphics[width=0.45\textwidth]{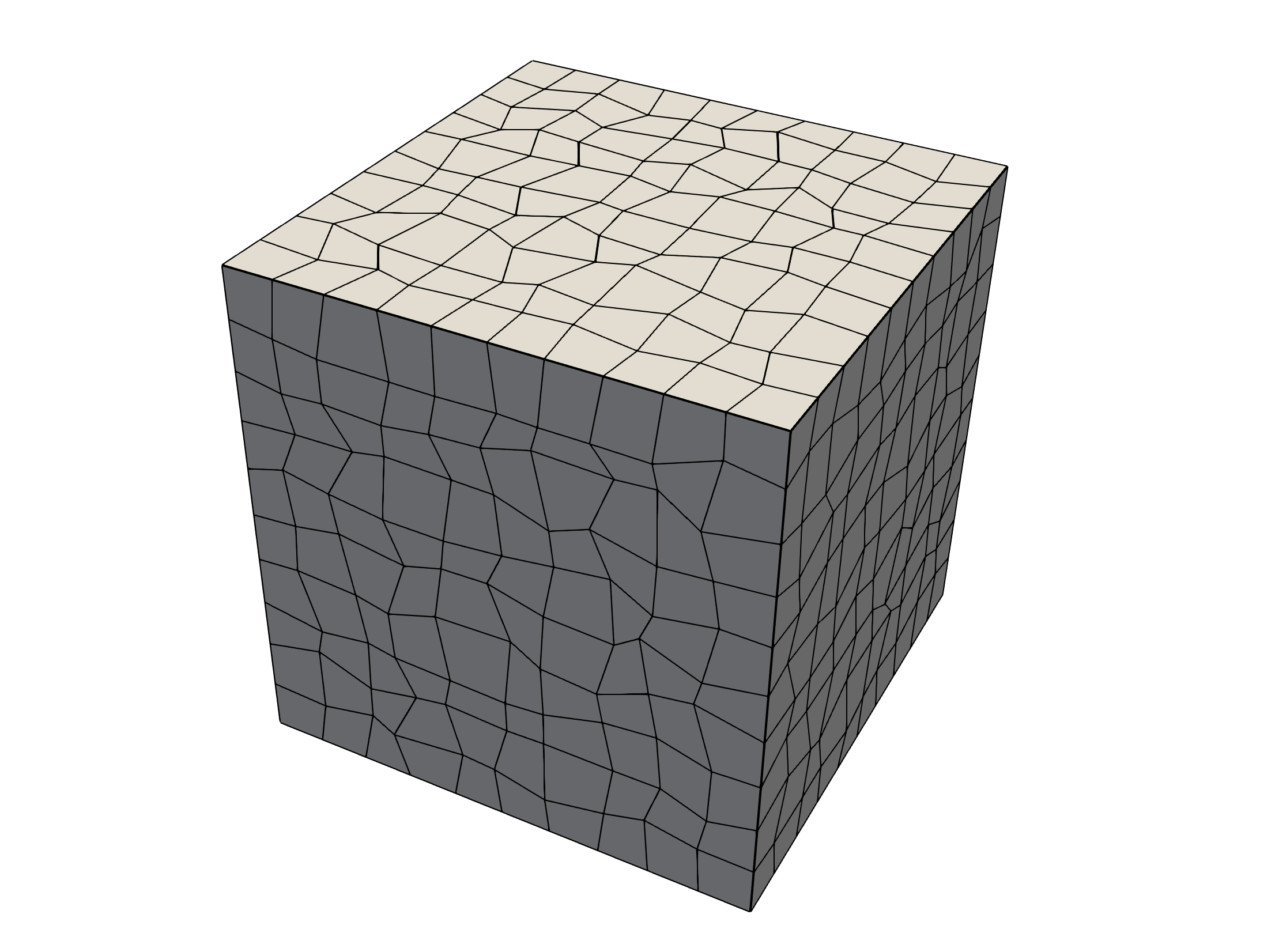}  
   	}
	\caption{Meshes used in the method of manufactured solutions case}
	\label{fig:mms_all}
\end{figure}

\section{Analytical Expressions for the Stress and Displacement Distributions around a Spherical Cavity under Uniaxial Tensions}
\label{app:sphericalCavity}
The analytical expressions for the stress distributions around the cavity, first derived by \citet{Southwell1926}, are given in cylindrical coordinates ($r$, $\theta$, $z$) as
\begin{eqnarray}
	\sigma_{rr} &=&
		\frac{T}{14 - 10\nu} \frac{a^3}{R^3}
		\left[ 9 - 15\nu - 12 \frac{a^2}{R^2}  - \frac{r^2}{R^2} \left( 72 - 15\nu - 105 \frac{a^2}{R^2} \right) + 15 \frac{r^4}{R^4} \left( 5 - 7 \frac{a^2}{R^2} \right) \right], \\
	\sigma_{\theta\theta} &=&
		\frac{T}{14 - 10\nu} \frac{a^3}{R^3}
		\left[ 9 - 15\nu - 12 \frac{a^2}{R^2}  - 15 \frac{r^2}{R^2} \left( 1 - 2\nu - \frac{a^2}{R^2} \right) \right], \\
	\sigma_{zz} &=&
		T \left[ 1 - \frac{1}{14 - 10\nu} \frac{a^3}{R^3} \left\{ 38 - 10\nu - 24 \frac{a^2}{R^2} 
		- \frac{r^2}{R^2} \left( 117 - 15\nu - 120 \frac{a^2}{R^2} \right)
		+ 15 \frac{r^4}{R^4} \left( 5 - 7 \frac{a^2}{R^2} \right) \right\} \right], \notag \\
		\\
	\sigma_{zr} &=&
	\frac{T}{14 - 10\nu} \frac{a^3 z r}{R^5}
	\left[ -3(19 - 5\nu) + 60 \frac{a^2}{R^2} + 15 \frac{r^2}{R^2} \left( 5 - 7 \frac{a^2}{R^2} \right)  \right].
\end{eqnarray}
where $a$ is the hole radius, $T$ is the distant stress applied in the $z$ direction, $\nu$ is the Poisson's ratio, $r^2 = x^2 + y^2$ is the cylinderical radial coordinate, $R^2 = r^2 + z^2$ is the spherical radial coodinate, and $x$, $y$, $z$ are the Cartesian coordinates.

To obtain Cartesian stresses \(\boldsymbol{\sigma}^{\mathrm{cart}}\) we rotate the cylindrical tensor
\(\boldsymbol{\sigma}^{\mathrm{cyl}}\) with the in-plane rotation
\[
Q(\theta)=
\begin{bmatrix}
\cos\theta & -\sin\theta & 0\\
\sin\theta & \phantom{-}\cos\theta & 0\\
0 & 0 & 1
\end{bmatrix},\qquad
\boldsymbol{\sigma}^{\mathrm{cart}} = Q\,\boldsymbol{\sigma}^{\mathrm{cyl}}\,Q^{\text{T}},\quad
\theta=\operatorname{atan2}(y,x).
\]
Because \(\sigma_{r\theta}=\sigma_{\theta z}=0\), the Cartesian components reduce to
\[
\begin{aligned}
\sigma_{xx} &= \sigma_{rr}\cos^2\theta + \sigma_{\theta\theta}\sin^2\theta,\\
\sigma_{yy} &= \sigma_{rr}\sin^2\theta + \sigma_{\theta\theta}\cos^2\theta,\\
\sigma_{xy} &= (\sigma_{rr}-\sigma_{\theta\theta})\sin\theta\cos\theta,\\
\sigma_{xz} &= \sigma_{rz}\cos\theta,\qquad
\sigma_{yz} = \sigma_{rz}\sin\theta,\\
\sigma_{zz} &= \sigma_{zz}.
\end{aligned}
\]

The displacement distributions were later derived by \citet{Goodier1933}:
\begin{eqnarray}
u_r &=& -\frac{A}{r^2} - \frac{3B}{r^4} + \left[ \frac{5-4\nu}{1-2\nu} \frac{C}{r^2}-9\frac{B}{r^4} \right]\cos (2\theta),\\
u_{\theta} &=& - \left[ \frac{2C}{r^2} + 6\frac{B}{r^4}  \right]\sin(2\theta),
\end{eqnarray}
where the constants $A$, $B$ and $C$ are defined as follows:
\begin{equation}
\frac{A}{a^3} = -\frac{T}{8\mu}\frac{13-10\nu}{7-5\nu}, \qquad
\frac{B}{a^5} = \frac{T}{8\mu}\frac{1}{7-5\nu}, \qquad
\frac{C}{a^3} = \frac{T}{8\mu}\frac{5(1-2\nu)}{7-5\nu}.
\end{equation}
and $\mu$ is the shear modulus.

Transforming the spherical-polar displacement \((u_r,u_\theta,0)\) to Cartesian
(using \(x=r\sin\theta\cos\phi\), \(y=r\sin\theta\sin\phi\), \(z=r\cos\theta\)):
\begin{align}
u_x^{\text{cav}} &= u_r\sin\theta\cos\phi + u_\theta\cos\theta\cos\phi,\\
u_y^{\text{cav}} &= u_r\sin\theta\sin\phi + u_\theta\cos\theta\sin\phi,\\
u_z^{\text{cav}} &= u_r\cos\theta - u_\theta\sin\theta.
\end{align}

The uniform uniaxial contribution (no cavity) is
\begin{equation}
\mathbf u_{\infty}(\mathbf x)
= \left(-\nu\,\frac{T}{E}\,x,\; -\nu\,\frac{T}{E}\,y,\; \frac{T}{E}\,z\right),
\end{equation}
with \(E\) the Young’s modulus.

Hence, the total Cartesian displacement field is
\begin{equation}
\mathbf u(\mathbf x)
= \big(u_x^{\text{cav}},\,u_y^{\text{cav}},\,u_z^{\text{cav}}\big)
  + \mathbf u_{\infty}(\mathbf x).
\end{equation}

\end{appendices}

\bibliography{bibliography}

@misc{moose_hypre_doc_2022,
  author       = {MOOSE Finite ELement Software Documentation},
  title        = {Hypre / BoomerAMG — Application Usage Documentation, MOOSE Framework},
  howpublished = {\url{https://mooseframework.inl.gov/releases/moose/2022-06-10/application_development/hypre.html}},
  year         = {2022},
  note         = {Accessed: 2025-09-18}
}

@article{Whelan2025,
author = {Whelan, A. and Tang, T. and Pakrashi, V. and Cardiff, P.},
year = {2025},
title = {A Finite Volume Framework for Damage and Fracture Prediction in Wire Drawing},
journal = {Int J Numer Methods Eng},
volume = {126},
number = {e7640},
doi = {https://doi.org/10.1002/nme.7640}
}

@article{Giudicelli2024,
   title = {3.0 - {MOOSE}: Enabling massively parallel multiphysics simulations},
   author = {Guillaume Giudicelli and Alexander Lindsay and Logan Harbour and Casey Icenhour and
             Mengnan Li and Joshua E. Hansel and Peter German and Patrick Behne and Oana Marin and
             Roy H. Stogner and Jason M. Miller and Daniel Schwen and Yaqi Wang and Lynn Munday and
             Sebastian Schunert and Benjamin W. Spencer and Dewen Yushu and Antonio Recuero and
             Zachary M. Prince and Max Nezdyur and Tianchen Hu and Yinbin Miao and
             Yeon Sang Jung and Christopher Matthews and April Novak and Brandon Langley and
             Timothy Truster and Nuno Nobre and Brian Alger and David Andr{\v{s}} and
             Fande Kong and Robert Carlsen and Andrew E. Slaughter and John W. Peterson and
             Derek Gaston and Cody Permann},
    year = {2024},
 journal = {{SoftwareX}},
  volume = {26},
   pages = {101690},
    issn = {2352-7110},
     doi = {https://doi.org/10.1016/j.softx.2024.101690},
     url = {https://www.sciencedirect.com/science/article/pii/S235271102400061X},
}

@article{Zhu2025,
title = {Development and validation of fully-coupled subchannel analysis code for LMR within the MOOSE framework},
journal = {Annals of Nuclear Energy},
volume = {214},
pages = {111199},
year = {2025},
issn = {0306-4549},
doi = {https://doi.org/10.1016/j.anucene.2025.111199},
url = {https://www.sciencedirect.com/science/article/pii/S0306454925000167},
author = {Xinyang Zhu and Jinshun Wang and Ronghua Chen and Yuhang Niu and Yingwei Wu and Dalin Zhang and Wenxi Tian}
}

@article{Wu2024,
  author    = {Wu, Y. and Zhang, H. and Liu, L. and Tang, H. and Dou, Q. and Guo, J. and Li, F.},
  title     = {An Efficient and Robust ILU(k) Preconditioner for Steady-State Neutron Diffusion Problem Based on MOOSE},
  journal   = {Energies},
  year      = {2024},
  volume    = {17},
  pages     = {1499},
  doi       = {10.3390/en17061499}
}

@article{Williamson2021,
author = {Richard L. Williamson and Jason D. Hales and Stephen R. Novascone and Giovanni Pastore and Kyle A. Gamble and Benjamin W. Spencer and Wen Jiang and Stephanie A. Pitts and Albert Casagranda and Daniel Schwen and Adam X. Zabriskie and Aysenur Toptan and Russell Gardner and Christoper Matthews and Wenfeng Liu and Hailong Chen},
title = {BISON: A Flexible Code for Advanced Simulation of the Performance of Multiple Nuclear Fuel Forms},
journal = {Nuclear Technology},
volume = {207},
number = {7},
pages = {954--980},
year = {2021},
publisher = {Taylor \& Francis},
doi = {10.1080/00295450.2020.1836940},
URL = { https://doi.org/10.1080/00295450.2020.1836940},
eprint = { https://doi.org/10.1080/00295450.2020.1836940}
}

@article{Africa2024,
  author  = {Pasquale C. Africa and Daniel Arndt and Wolfgang Bangerth and Bruno Blais and
             Marc Fehling and Rene Gassm{\"o}ller and Timo Heister and Luca Heltai and
             Sebastian Kinnewig and Martin Kronbichler and Matthias Maier and Peter Munch and
             Magdalena Schreter-Fleischhacker and Jan P. Thiele and Bruno Turcksin and
             David Wells and Vladimir Yushutin},
  title   = {The deal.II library, Version 9.6},
  journal = {Journal of Numerical Mathematics},
  year    = 2024,
  volume  = 32,
  number  = 4,
  pages   = {369--380},
  doi     = {10.1515/jnma-2024-0137}
}

@article{Chacon2025,
title = {A scalable multidimensional fully implicit solver for Hall magnetohydrodynamics},
journal = {Journal of Computational Physics},
volume = {526},
pages = {113789},
year = {2025},
issn = {0021-9991},
doi = {https://doi.org/10.1016/j.jcp.2025.113789},
url = {https://www.sciencedirect.com/science/article/pii/S0021999125000725},
author = {L. Chacón}
}

@article{Sukas2025,
title = {A robust monolithic nonlinear Newton method for the compressible Reynolds averaged Navier–Stokes Equations},
journal = {Computers \& Fluids},
volume = {289},
pages = {106549},
year = {2025},
issn = {0045-7930},
doi = {https://doi.org/10.1016/j.compfluid.2025.106549},
url = {https://www.sciencedirect.com/science/article/pii/S0045793025000106},
author = {Hulya Sukas and Mehmet Sahin}
}

@article{Ahmed2025,
title = {Modified preconditioned Newton-Krylov approaches for Navier-Stokes equations using nodal integral method},
journal = {Computers \& Mathematics with Applications},
volume = {181},
pages = {163-192},
year = {2025},
issn = {0898-1221},
doi = {https://doi.org/10.1016/j.camwa.2024.12.027},
url = {https://www.sciencedirect.com/science/article/pii/S0898122124005777},
author = {Nadeem Ahmed and Suneet Singh and Ram Prakash Bharti}
}

@article{Liu2025,
author = {Lixun Liu and Han Zhang and Xinru Peng and Qinrong Dou and Yingjie Wu and Jiong Guo and Fu Li},
title = {A Parallel Coloring Newton-Krylov Method for Multiphysics Coupling System in Nuclear Reactors},
journal = {Nuclear Science and Engineering},
volume = {199},
number = {1},
pages = {61--81},
year = {2025},
publisher = {Taylor \& Francis},
doi = {10.1080/00295639.2024.2344956}
}

@article{Eilmer2023,
title = {Eilmer: An open-source multi-physics hypersonic flow solver},
journal = {Computer Physics Communications},
volume = {282},
pages = {108551},
year = {2023},
issn = {0010-4655},
doi = {https://doi.org/10.1016/j.cpc.2022.108551},
url = {https://www.sciencedirect.com/science/article/pii/S0010465522002703},
author = {Nicholas N. Gibbons and Kyle A. Damm and Peter A. Jacobs and Rowan J. Gollan}
}

@article{Amir2021,
  author    = {Laila Amir and Michel Kern},
  title     = {Jacobian Free Methods for Coupling Transport with Chemistry in Heterogeneous Porous Media},
  journal   = {Water},
  year      = {2021},
  volume    = {13},
  number    = {3},
  pages     = {370},
  doi       = {10.3390/w13030370},
  url       = {https://doi.org/10.3390/w13030370}
}

@inproceedings{Johannes2021,
  author    = {K. Johannes and L. Versbach and P. Birken and G. Gassner and R. Klöfkorn},
  title     = {A Finite Volume Based Multigrid Preconditioner for DG-SEM for Convection-Diffusion},
  booktitle = {WCCM-ECCOMAS 2020},
  year      = {2021},
  url       = {https://www.scipedia.com/public/Kasimir_et_al_2021}
}

@article{Nguyen2022,
author = {Nguyen, T. Binh and De  Sterck, Hans and Freret, Lucie and Groth, Clinton P. T.},
title = {High-order implicit time-stepping with high-order central essentially-non-oscillatory methods for unsteady three-dimensional computational fluid dynamics simulations},
journal = {International Journal for Numerical Methods in Fluids},
volume = {94},
number = {2},
pages = {121-151},
doi = {https://doi.org/10.1002/fld.5049},
url = {https://onlinelibrary.wiley.com/doi/abs/10.1002/fld.5049},
eprint = {https://onlinelibrary.wiley.com/doi/pdf/10.1002/fld.5049},
year = {2022}
}

@article{Zhang2024twofluid,
title = {Highly robust Jacobian-free Newton–Krylov method for solving fully implicit two-fluid equations},
journal = {International Communications in Heat and Mass Transfer},
volume = {159},
pages = {108260},
year = {2024},
issn = {0735-1933},
doi = {https://doi.org/10.1016/j.icheatmasstransfer.2024.108260},
url = {https://www.sciencedirect.com/science/article/pii/S0735193324010224},
author = {Yuhang Zhang and Zhaofei Tian and Lei Li and Guangliang Chen and Dabin Sun and Rui Li and Hao Qian and Lixuan Zhang and Jinchao Li}
}

@article{Ma2024,
    author = {Ma, Pengfei and Cai, Li and Wang, Xuan and Wang, Yongheng and Luo, Xiaoyu and Gao, Hao},
    title = {An unconditionally stable scheme for the immersed boundary method with application in cardiac mechanics},
    journal = {Physics of Fluids},
    volume = {36},
    number = {8},
    pages = {081914},
    year = {2024},
    month = {08},
    issn = {1070-6631},
    doi = {10.1063/5.0225605},
    url = {https://doi.org/10.1063/5.0225605},
    eprint = {https://pubs.aip.org/aip/pof/article-pdf/doi/10.1063/5.0225605/20131633/081914\_1\_5.0225605.pdf},
}

@article{Munch2024,
title = {On the construction of an efficient finite-element solver for phase-field simulations of many-particle solid-state-sintering processes},
journal = {Computational Materials Science},
volume = {231},
pages = {112589},
year = {2024},
issn = {0927-0256},
doi = {https://doi.org/10.1016/j.commatsci.2023.112589},
url = {https://www.sciencedirect.com/science/article/pii/S0927025623005839},
author = {Peter Munch and Vladimir Ivannikov and Christian Cyron and Martin Kronbichler},
}

@article{Zhang2024,
    author = {Zhang, Jianshe and Zhang, Ziqing and Dong, Xu and Yuan, Hang and Zhang, Yanfeng and Lu, Xingen},
    title = {A robust Jacobian-free Newton–Krylov method for turbomachinery simulations},
    journal = {Physics of Fluids},
    volume = {36},
    number = {12},
    pages = {126116},
    year = {2024},
    month = {12},
    issn = {1070-6631},
    doi = {10.1063/5.0243628},
    url = {https://doi.org/10.1063/5.0243628},
    eprint = {https://pubs.aip.org/aip/pof/article-pdf/doi/10.1063/5.0243628/20283515/126116\_1\_5.0243628.pdf},
}

@article{Syrakos2023,
  title={A unification of least-squares and Green--Gauss gradients under a common projection-based gradient reconstruction framework},
  author={Syrakos, Alexandros and Oxtoby, Oliver and de Villiers, Eugene and Varchanis, Stylianos and Dimakopoulos, Yannis and Tsamopoulos, John},
  journal={Mathematics and Computers in Simulation},
  volume={205},
  pages={108--141},
  year={2023},
  month={March},
  doi={10.1016/j.matcom.2022.09.008},
  publisher={Elsevier}
}

@incollection{Pellegrini2012,
  author       = {François Pellegrini},
  title        = {Scotch and {PT-Scotch} Graph Partitioning Software: An Overview},
  booktitle    = {Combinatorial Scientific Computing},
  editor       = {Uwe Naumann and Olaf Schenk},
  publisher    = {Chapman and Hall/CRC},
  year         = {2012},
  chapter      = {14},
  pages        = {373--406},
  doi          = {10.1201/b11644-15},
  url          = {https://hal.inria.fr/hal-00770422},
  address      = {London, United Kingdom}
}

@article{Demirdzic2022,
  author = {I. Demirdžić and P. Cardiff},
  title = {Symmetry plane boundary conditions for cell-centered finite-volume continuum mechanics},
  journal = {Numerical Heat Transfer, Part B: Fundamentals},
  year = {2022},
  doi = {10.1080/10407790.2022.2105073}
}

@inproceedings{Nishikawa2010,
  author    = {Hiroaki Nishikawa},
  title     = {Beyond Interface Gradient: A General Principle for Constructing Diffusion Schemes},
  booktitle = {40th Fluid Dynamics Conference and Exhibit},
  year      = {2010},
  address   = {Chicago, Illinois, USA},
  month     = {June 28 -- July 1},
  doi       = {10.2514/6.2010-5093},
  publisher = {AIAA},
  note      = {Published online: 13 Nov 2012},
}

@article{Nishikawa2017,
  author = {Hiroaki Nishikawa and Yoshitaka Nakashima and Norihiko Watanabe},
  title = {Effects of high-frequency damping on iterative convergence of implicit viscous solver},
  journal = {Journal of Computational Physics},
  volume = {348},
  pages = {66--81},
  year = {2017},
  publisher = {Elsevier},
  doi = {10.1016/j.jcp.2017.07.021}
}

@article{Cardiff2018,
  author =       {P. Cardiff and A. Kara\v{c} and P. De Jaeger and H. Jasak and J. Nagy and A. Ivankovi\'{c} and {\v{Z}}. Tukovi\'{c}},
  title =        {An open-source finite volume toolbox for solid mechanics and fluid-solid interaction simulations},
  year =         {2018},
  doi =          {10.48550/arXiv.1808.10736},
  note =         {\url{https://arxiv.org/abs/1808.10736}}
}

@article{Badcock1996,
  author = {K. J. Badcock and A. L. Gaitonde},
  title = {An unfactored implicit moving mesh method for the two-dimensional unsteady {N-S} equations},
  journal = {International Journal for Numerical Methods in Fluids},
  volume = {23},
  number = {6},
  pages = {607--631},
  year = {1996},
  publisher = {Wiley},
  doi = {10.1002/(SICI)1097-0363(19960930)23:6<607::AID-FLD433>3.0.CO;2-E},
  url = {https://doi.org/10.1002/(SICI)1097-0363(19960930)23:6<607::AID-FLD433>3.0.CO;2-E},
  note = {First published: 30 September 1996}
}

@article{Pernice2001,
  author  = {Michael Pernice and Michael D. Tocci},
  title   = {A Multigrid-Preconditioned {Newton}--{Krylov} Method for the Incompressible {Navier}--{Stokes} Equations},
  journal = {SIAM Journal on Scientific Computing},
  volume  = {23},
  number  = {2},
  pages   = {398--418},
  year    = {2001},
  doi     = {10.1137/S1064827500372250}
}

@article{Geuzaine2001,
author = {Geuzaine, Philippe},
title = {{Newton}-{Krylov} Strategy for Compressible Turbulent Flows on Unstructured Meshes},
journal = {AIAA Journal},
volume = {39},
number = {3},
pages = {528-531},
year = {2001},
doi = {10.2514/2.1339}
}

@article{Qin2000,
  author    = {Ning Qin and David K. Ludlow and Scott T. Shaw},
  title     = {A matrix-free preconditioned {Newton}/{GMRES} method for unsteady {Navier}--{Stokes} solutions},
  journal   = {International Journal for Numerical Methods in Fluids},
  volume    = {33},
  number    = {2},
  pages     = {223--248},
  year      = {2000},
  month     = {May},
  doi       = {10.1002/(SICI)1097-0363(20000530)33:2<223::AID-FLD10>3.0.CO;2-V}
}

@article{Brown1990,
  author  = {P. N. Brown and Y. Saad},
  title   = {Hybrid {Krylov} methods for nonlinear systems of equations},
  journal = {SIAM Journal on Scientific and Statistical Computing},
  volume  = {11},
  number  = {3},
  pages   = {450--481},
  year    = {1990},
  doi     = {10.1137/0911026}
}

@article{Chan1984,
  author  = {T. F. Chan and K. R. Jackson},
  title   = {Nonlinearly preconditioned {Krylov} subspace methods for discrete {Newton} algorithms},
  journal = {SIAM Journal on Scientific and Statistical Computing},
  volume  = {5},
  number  = {3},
  pages   = {533--542},
  year    = {1984},
  doi     = {10.1137/0905039}
}

@article{Brown1986,
  author  = {P. N. Brown and A. C. Hindmarsh},
  title   = {Matrix-free methods for stiff systems of {ODE}'s},
  journal = {SIAM Journal on Numerical Analysis},
  volume  = {23},
  number  = {3},
  pages   = {610--638},
  year    = {1986},
  doi     = {10.1137/0723034}
}

@article{Gear1983,
  author  = {C. W. Gear and Y. Saad},
  title   = {Iterative solution of linear equations in {ODE} codes},
  journal = {SIAM Journal on Scientific and Statistical Computing},
  volume  = {4},
  number  = {4},
  pages   = {583--601},
  year    = {1983},
  doi     = {10.1137/0904040}
}

@article{Mchugh1994,
  author    = {Paul R. McHugh and Dana A. Knoll},
  title     = {Comparison of standard and matrix-free implementations of several {Newton}--{Krylov} solvers},
  journal   = {AIAA Journal},
  volume    = {32},
  number    = {12},
  pages     = {2394--2400},
  year      = {1994},
  doi       = {10.2514/3.12305} 
}

@article{Lucas2010,
  author    = {Peter Lucas and Alexander H. van Zuijlen and Hester Bijl},
  title     = {Fast unsteady flow computations with a Jacobian-free {Newton}--{Krylov} algorithm},
  journal   = {Journal of Computational Physics},
  volume    = {229},
  number    = {24},
  pages     = {9201--9215},
  year      = {2010},
  month     = {December},
  publisher = {Elsevier},
  doi       = {10.1016/j.jcp.2010.08.033}
}

@article{Vaassen2008,
  author    = {J.-M. Vaassen and D. Vigneron and J.-A. Essers},
  title     = {An implicit high order finite volume scheme for the solution of 3D {Navier}--{Stokes} equations with new discretization of diffusive terms},
  journal   = {Journal of Computational and Applied Mathematics},
  volume    = {215},
  number    = {2},
  pages     = {595--601},
  year      = {2008},
  month     = {June},
  publisher = {Elsevier},
  doi       = {10.1016/j.cam.2006.04.066}
}

@article{Nejat2011,
  author  = {Amir Nejat and Alireza Jalali and Mahkame Sharbatdar},
  title   = {A {Newton}--{Krylov} finite volume algorithm for the power-law non-{Newtonian} fluid flow using pseudo-compressibility technique},
  journal = {Journal of Non-Newtonian Fluid Mechanics},
  volume  = {166},
  number  = {19--20},
  pages   = {1158--1172},
  year    = {2011},
  month   = {October},
  doi     = {10.1016/j.jnnfm.2011.07.003},
  publisher = {Elsevier}
}

@article{Nejat2008,
  author  = {Amir Nejat and Carl Ollivier-Gooch},
  title   = {A high-order accurate unstructured finite volume {Newton}--{Krylov} algorithm for inviscid compressible flows},
  journal = {Journal of Computational Physics},
  volume  = {227},
  number  = {4},
  pages   = {2582--2609},
  year    = {2008},
  month   = {February},
  doi     = {10.1016/j.jcp.2007.11.011},
  publisher = {Elsevier}
}

@InProceedings{hypre,
author="Falgout, Robert D.
and Yang, Ulrike Meier",
editor="Sloot, Peter M. A.
and Hoekstra, Alfons G.
and Tan, C. J. Kenneth
and Dongarra, Jack J.",
title="hypre: A Library of High Performance Preconditioners",
booktitle="Computational Science --- ICCS 2002",
year="2002",
publisher="Springer Berlin Heidelberg",
address="Berlin, Heidelberg",
pages="632--641",
}

@article{MUMPS:1,
   title   = {A Fully Asynchronous Multifrontal Solver Using Distributed Dynamic Scheduling},
   author  = {P.R. Amestoy and I. S. Duff and J. Koster and J.-Y. L'Excellent},
   journal = {SIAM Journal on Matrix Analysis and Applications},
   volume  = {23},
   number  = {1},
   year    = {2001},
   pages   = {15-41}
 }

@article{MUMPS:2,
  title = {{Performance and Scalability of the Block Low-Rank Multifrontal
  Factorization on Multicore Architectures}},
  author = {P.R. Amestoy and A. Buttari and J.-Y. L'Excellent and T. Mary},
  journal = {ACM Transactions on Mathematical Software},
  volume = 45,
  issue = 1,
  pages = {2:1--2:26},
  year={2019},
}

@article{Degroote2009,
  title={Performance of a new partitioned procedure versus a monolithic procedure in fluid-structure interaction},
  author={Degroote, Joris and Bathe, Klaus-J{\"u}rgen and Vierendeels, Jan},
  journal={Computers \& Structures},
  volume={87},
  number={11-12},
  pages={793--801},
  year={2009},
  publisher={Elsevier},
  doi={10.1016/j.compstruc.2009.01.006}
}

@misc{Pelteret2018,
  author       = {Pelteret, Jean-Paul and McBride, Andrew},
  title        = {{The deal.II code gallery: Quasi-Static Finite-Strain Compressible Elasticity}},
  month        = apr,
  year         = 2018,
  publisher    = {Zenodo},
  doi          = {10.5281/zenodo.1228964},
  url          = {https://doi.org/10.5281/zenodo.1228964}
}

@misc{Simplas,
  author       = {P. Areias},
  title        = {{Simplas}},
  publisher    = {Portuguese Software Association (ASSOFT)},
  url          = {http://www.simplassoftware.com.}
}

@misc{WebPlotDigitizer,
    author = {Ankit Rohatgi},
    title = {WebPlotDigitizer},
    url = {https://automeris.io},
    version = {5.2},
}

@article{Cesar2001,
    author = {César de Sá, José M. A. and Areias, Pedro M. A. and Natal Jorge, Renato M.},
	title = {Quadrilateral elements for the solution of elasto-plastic finite strain problems},
	journal = {International Journal for Numerical Methods in Engineering},
	volume = {51},
	number = {8},
	pages = {883-917},
	doi = {https://doi.org/10.1002/nme.183},
	url = {https://onlinelibrary.wiley.com/doi/abs/10.1002/nme.183},
	eprint = {https://onlinelibrary.wiley.com/doi/pdf/10.1002/nme.183},
	year = {2001}
}

@incollection{Turek2006,
  title={Proposal for numerical benchmarking of fluid-structure interaction between an elastic object and laminar incompressible flow},
  author={Stefan Turek and Jaroslav Hron},
  booktitle={Fluid-Structure Interaction: Modelling, Simulation, Optimisation},
  pages={371--385},
  year={2006},
  publisher={Springer},
  series={Lecture Notes in Computational Science and Engineering},
  volume={53},
  doi={10.1007/3-540-34596-5_20},
  editor={Hans-Joachim Bungartz and Michael Schäfer},
  address={Berlin, Heidelberg},
}

@article{Goodier1933,
  title={Concentration of stress around spherical and cylindrical inclusions and flaws},
  author={Goodier, J. N.},
  journal={Journal of Applied Mechanics},
  volume={1},
  number={2},
  pages={39--44},
  year={1933}
}

@book{Zienkiewicz2000,
title={The finite element method},
author={Zienkiewicz, Olgierd Cecil and Taylor, Robert Leroy},
volume={2},
year={2000},
publisher={Butterworth-heinemann},
address={Berlin, Germany}
}

@article{Cardiff2016,
  author = {Cardiff, P. and Tukovi\'{c}, {\v{Z}}. and Jasak, H. and Ivankovi\'{c}, A.},
  title = {A block{-}coupled Finite Volume methodology for linear elasticity and unstructured meshes},
  journal = {{C}omputers \& {S}tructures},
  volume = {175},
  pages = {100 - 122},
  year = {2016},
  issn = {0045-7949},
  doi = {https://doi.org/10.1016/j.compstruc.2016.07.004},
  url = {http://www.sciencedirect.com/science/article/pii/S0045794916306046}
}

@article{Das2011,
  author = { Shankhadeep   Das  and  Sanjay R.   Mathur  and  Jayathi Y.   Murthy },
  title = {An {U}nstructured {F}inite-{V}olume {M}ethod for {S}tructure–{E}lectrostatics {I}nteractions in {MEMS}},
  journal = {Numerical Heat Transfer, Part B: Fundamentals},
  volume = {60},
  number = {6},
  pages = {425-451},
  year  = {2011},
  publisher = {Taylor & Francis},
  doi = {10.1080/10407790.2011.628252},
  URL = {https://doi.org/10.1080/10407790.2011.628252},
  eprint = {https://doi.org/10.1080/10407790.2011.628252}
}

@article{Castrillo2024,
  author = {Pablo Castrillo and Eugenio Schillaci and Joaquim Rigola},
  title = {High-order cell-centered finite volume method for solid dynamics on unstructured meshes},
  journal = {Computers \& Structures},
  volume = {295},
  pages = {107288},
  year = {2024},
  issn = {0045-7949},
  doi = {https://doi.org/10.1016/j.compstruc.2024.107288},
  url = {https://www.sciencedirect.com/science/article/pii/S0045794924000178}
}

@article{Nishikawa2020,
  title = {A hyperbolic Poisson solver for tetrahedral grids},
  journal = {Journal of Computational Physics},
  volume = {409},
  pages = {109358},
  year = {2020},
  issn = {0021-9991},
  doi = {https://doi.org/10.1016/j.jcp.2020.109358},
  url = {https://www.sciencedirect.com/science/article/pii/S0021999120301327},
  author = {Hiroaki Nishikawa}
}

@article{Jacobs1986,
  author = {Jsvobs, D. A. H.},
  title = "{A Generalization of the Conjugate-Gradient Method to Solve Complex Systems}",
  journal = {IMA Journal of Numerical Analysis},
  volume = {6},
  number = {4},
  pages = {447-452},
  year = {1986},
  month = {10},
  issn = {0272-4979},
  doi = {10.1093/imanum/6.4.447},
  url = {https://doi.org/10.1093/imanum/6.4.447},
  eprint = {https://academic.oup.com/imajna/article-pdf/6/4/447/2181350/6-4-447.pdf},
}

@article{Knoll2004,
  author = {D.A. Knoll and D.E. Keyes},
  title = {Jacobian-free {Newton}–{Krylov} methods: a survey of approaches and applications},
  journal = {Journal of Computational Physics},
  volume = {193},
  number = {2},
  pages = {357-397},
  year = {2004},
  issn = {0021-9991},
  doi = {https://doi.org/10.1016/j.jcp.2003.08.010},
  url = {https://www.sciencedirect.com/science/article/pii/S0021999103004340}
}

@article{Cardiff2017,
  author = {Cardiff, P. and {\v{Z}}. Tukovi{\'c} and Jaeger, P. De and Clancy, M. and Ivankovi{\'c}, A.},
  title = {A {L}agrangian cell-centred finite volume method for metal forming simulation},
  journal = {International {J}ournal for {N}umerical {M}ethods in {E}ngineering},
  volume = {109},
  number = {13},
  pages = {1777-1803},
  year = {2017},
  doi = {https://doi.org/10.1002/nme.5345},
  url = {https://onlinelibrary.wiley.com/doi/abs/10.1002/nme.5345},
  eprint = {https://onlinelibrary.wiley.com/doi/pdf/10.1002/nme.5345}
}

@article{Jasak2000,
  author = {Jasak, H. and Weller, H. G.},
  title = {Application of the finite volume method and unstructured meshes to linear elasticity},
  journal = {International {J}ournal for {N}umerical {M}ethods in {E}ngineering},
  volume = {48},
  number = {2},
  pages = {267-287},
  year = {2000},
  doi = {https://doi.org/10.1002/(SICI)1097-0207(20000520)48:2<267::AID-NME884>3.0.CO;2-Q},
  url = {https://onlinelibrary.wiley.com/doi/abs/10.1002/\%28SICI\%291097-0207\%2820000520\%2948\%3A2\%3C267\%3A\%3AAID-NME884\%3E3.0.CO\%3B2-Q},
  eprint = {https://onlinelibrary.wiley.com/doi/pdf/10.1002/\%28SICI\%291097-0207\%2820000520\%2948\%3A2\%3C267\%3A\%3AAID-NME884\%3E3.0.CO\%3B2-Q}
}

@misc{Jasak2011,
  author = {Hrvoje Jasak},
  title = {Finite Volume Discretisation with Polyhedral Cell Support: What is Discretisation?},
  year = {2011},
  note = {Presentation at NUMAP-FOAM Summer School, Zagreb, Croatia, 2-15 September 2009},
  url = {https://www.slideshare.net/slideshow/57969246-finitevolume/16128767}
}

@article{Rhie1983,
  author = {Rhie, C. M. and Chow, W. L.},
  title = {Numerical study of the turbulent flow past an airfoil with trailing edge separation},
  journal = {AIAA Journal},
  volume = {21},
  number = {11},
  pages = {1525-1532},
  year = {1983},
  doi = {10.2514/3.8284},
  URL = {https://doi.org/10.2514/3.8284 },
  eprint = {https://doi.org/10.2514/3.8284}
}

@article{Tukovic2013,
  author = {Tukovi\'{c}, {\v{Z}} and Ivankovi\'{c}, A. and Kara\v{c}, A.},
  title = {Finite{-}volume stress analysis in multi{-}material linear elastic body},
  journal = {International Journal for Numerical Methods in Engineering},
  volume = {93},
  number = {4},
  pages = {400-419},
  year = {2013},
  doi = {https://doi.org/10.1002/nme.4390},
  url = {https://onlinelibrary.wiley.com/doi/abs/10.1002/nme.4390},
  eprint = {https://onlinelibrary.wiley.com/doi/pdf/10.1002/nme.4390}
}

@article{Batistic2022,
  author = {Ivan Batisti\'{c} and Philip Cardiff and {\v{Z}}eljko Tukovi\'{c}},
  title = {A finite volume penalty based segment-to-segment method for frictional contact problems},
  journal = {Applied Mathematical Modelling},
  volume = {101},
  pages = {673-693},
  year = {2022},
  issn = {0307-904X},
  doi = {https://doi.org/10.1016/j.apm.2021.09.009},
  url = {https://www.sciencedirect.com/science/article/pii/S0307904X21004248}
}

@article{Demirdzic1995,
  author = {I. Demirdžić and S. Muzaferija},
  title = {Numerical method for coupled fluid flow, heat transfer and stress analysis using unstructured moving meshes with cells of arbitrary topology},
  journal = {Computer Methods in Applied Mechanics and Engineering},
  volume = {125},
  number = {1},
  pages = {235-255},
  year = {1995},
  issn = {0045-7825},
  doi = {https://doi.org/10.1016/0045-7825(95)00800-G},
  url = {https://www.sciencedirect.com/science/article/pii/004578259500800G}
}

@phdthesis{Jasak1996,
  title={Error analysis and estimation for the finite volume method with applications to fluid flows},
  author={Jasak, Hrvoje},
  year={1996},
  school={Imperial College London (University of London)}
}

@phdthesis{Mazzanti2024,
  title={Coupled Vertex-Centred Finite Volume Methods for Large-Strain Elastoplasticity},
  author={Mazzanti, Federico},
  year={2024},
  school={University College Dublin},
  note = {For examination}
}

@Misc{PETSc,
  author = {Satish Balay and Shrirang Abhyankar and Mark~F. Adams and Jed Brown and Peter Brune
            and Kris Buschelman and Lisandro Dalcin and Victor Eijkhout and William~D. Gropp
            and Dinesh Kaushik and Matthew~G. Knepley
            and Lois Curfman McInnes and Karl Rupp and Barry~F. Smith
            and Stefano Zampini and Hong Zhang},
  title =  {{PETS}c {W}eb page},
  url =    {http://www.mcs.anl.gov/petsc},
  howpublished = {\url{http://www.mcs.anl.gov/petsc}},
  year = {2015}
}

@article{Southwell1926,
  author = {R.V. Southwell and H.J. Gough},
  title = {VI. On the concentration of stress in the neighbourhood of a small spherical flaw; and on the propagation of fatigue fractures in “Statistically Isotropic” materials },
  journal = {The London, Edinburgh, and Dublin Philosophical Magazine and Journal of Science},
  volume = {1},
  number = {1},
  pages = {71--97},
  year = {1926},
  publisher = {Taylor \& Francis},
  doi = {10.1080/14786442608633614},
  URL = {https://doi.org/10.1080/14786442608633614},
  eprint = {https://doi.org/10.1080/14786442608633614}
}

@article{Guccione1995,
  author = {Julius M. Guccione and Kevin D. Costa and Andrew D. McCulloch},
  title = {Finite element stress analysis of left ventricular mechanics in the beating dog heart},
  journal = {Journal of Biomechanics},
  volume = {28},
  number = {10},
  pages = {1167-1177},
  year = {1995},
  issn = {0021-9290},
  doi = {https://doi.org/10.1016/0021-9290(94)00174-3},
  url = {https://www.sciencedirect.com/science/article/pii/0021929094001743}
}

@article{Simo1992,
  author = {Simo, J. C. and Armero, F.},
  title = {Geometrically non-linear enhanced strain mixed methods and the method of incompatible modes},
  journal = {International Journal for Numerical Methods in Engineering},
  volume = {33},
  number = {7},
  pages = {1413-1449},
  doi = {https://doi.org/10.1002/nme.1620330705},
  url = {https://onlinelibrary.wiley.com/doi/abs/10.1002/nme.1620330705},
  eprint = {https://onlinelibrary.wiley.com/doi/pdf/10.1002/nme.1620330705},
  year = {1992}
}

@book{Simo1998,
	author={Simo, J. C. and Hughes, T. J. R.},
	year={1998},
	title={Computational Inelasticity},
	publisher={Springer-Verlag},
	volume={7},
	address={New York}
}

@article {Land2015,
  author = {Land, Sander and Gurev, Viatcheslav and Arens, Sander and
              Augustin, Christoph M. and Baron, Lukas and Blake, Robert and
              Bradley, Chris and Castro, Sebatian and Crozier, Andrew and
              Favino, Marco and Fastl, Thomas E. and Fritz, Thomas and Gao, Hao and
              Gizzi, Alessio and Griffith, Boyce E. and Hurtado, Daniel E. and
              Krause, Rolf and Luo, Xiaoyu and Nash, Martyn P. and
              Pezzuto, Simone and Plank, Gernot and Rossi, Simone and
              Ruprecht, Daniel and Seemann, Gunnar and Smith, Nicolas P. and
              Sundnes, Joakim and Rice, J. Jeremy and Trayanova, Natalia and
              Wang, Dafand and Wang, Zhinuo Jenny and Niederer, Steven A.},
  title = {Verification of cardiac mechanics software: benchmark problems and
              solutions for testing active and passive material behaviour},
  jounral = {Proc. R. Soc. A},
  fjounral = {Proceedings of the Royal Society A: Mathematical, Physical and
              Engineering Sciences},
  volume = {471},
  year = {2015},
  month = {dec},
  number = {2184},
  pages = {20150641},
  issn = {1364-5021},
  doi = {10.1098/rspa.2015.0641},
  url = {https://doi.org/10.1098/rspa.2015.0641}
}

@article{geuzaine2009gmsh,
  title={Gmsh: A 3-D finite element mesh generator with built-in pre- and post-processing facilities},
  author={Geuzaine, Christophe and Remacle, Jean-Fran{\c{c}}ois},
  journal={International journal for numerical methods in engineering},
  volume={79},
  number={11},
  pages={1309--1331},
  year={2009},
  publisher={Wiley Online Library},
  doi={10.1002/nme.2579}
}

@article{Tukovic2007,
  author={Tukovi{\'c}, {\v{Z}}eljko and Jasak, Hrvoje},
  title = {Updated {L}agrangian finite volume solver for large deformation dynamic response of elastic body},
  pages = {55-70},
  journal = {Transactions of {FAMENA}},
  volume = {31},
  year = {2007},
  number = {1},
  issn_ = {1333-1124}
}

@article{Cardiff2016a,
  author =       {Cardiff, P. and Tukovi\'{c} and Jasak, H. and
                  Ivankovi\'{c}, A.},
  title =        {A Block-Coupled Finite Volume Methodology for Linear
                  Elasticity and Unstructured Meshes},
  journal =      {Computers \& Structures},
  year =         {2016},
  volume =       {175},
  pages =        {100-122},
  doi =          "10.1016/j.compstruc.2016.07.004"
}

@article{Demirdzic1988,
  author =       {I. Demird{\v{z}}i\'{c} and D. Martinovi\'{c} and
                  A. Ivankovi\'{c}},
  title =        {Numerical simulation of thermal deformation in
                  welded workpiece},
  journal =      {Zavarivanje},
  year =         {1988},
  volume =       {31},
  number =       {},
  pages =        {209-219},
  note =         {In Croatian. English translation available at \url{https://www.researchgate.net/profile/Alojz_Ivankovic/publication/296148474_Numerical_simulation_of_thermal_deformation_in_welded_workpiece/links/5d07642ba6fdcc39f12219eb/Numerical-simulation-of-thermal-deformation-in-welded-workpiece.pdf}},
}

@incollection{Demirdzic1997a,
  author =       {I. Demird{\v{z}}i\'{c} and S. Muzaferija and
                  M. Peri\'{c}},
  title =        {Advances in computation of heat transfer, fluid
                  flow, and solid body deformation using finite volume
                  approaches},
  booktitle =    {Advances in Numerical Heat Transfer, Chapter 2},
  year =         {1997},
  publisher =    {Taylor \& Francis},
  editor =       {W. J. Minkowycz, E. M. Sparrow},
  pages =        {59-96},
  address =      {London, United Kingdom}
}

@book{Ferziger2002,
  author =       {Ferziger, J. H. and Peric, M.},
  year =         {2002},
  title =        {Computational methods for fluid dynamics},
  edition =      {3$^{rd}$},
  publisher =    {Springer},
  address =      {Berlin, Germany}
}

@article{Fryer1991,
  author =       {Y. D. Fryer and C. Bailey and M. Cross and
                  C.-H. Lai},
  title =        {A control volume procedure for solving the elastic
                  stress-strain equations on an unstructured mesh},
  journal =      {Applied Mathematical Modelling},
  year =         {1991},
  volume =       {15},
  pages =        {639-645},
}

@article{Haider2017,
  author =       {J. Haider and C. H. Lee and A. J. Gil and J. Bonet},
  title =        {A first order hyperbolic framework for large strain
                  computational solid dynamics: An upwind cell centred
                  Total {Lagrangian} scheme},
  journal =      {International Journal for Numerical Methods in Engineering},
  year =         {2017},
  volume =       {109},
  pages =        {407-456},
}

@book{Hitchings1987,
  author =       {Hitchings, D. and Davies, G. A. O. and Kamoulakos,
                  A.},
  publisher =    {{International Association for the Engineering
                  Analysis Community \& National Agency for Finite
                  Element Methods \& Standards (NAFEMS)}},
  year =         {1987},
  title =        {Linear static benchmarks},
  address =      {Glasgow, UK}
}

@article{Kluth2010,
  author =       {G. Kluth and B. Despr\'{e}s},
  title =        {Discretization of hyperelasticity on unstructured
                  mesh with a cell-centered {Lagrangian} scheme},
  journal =      {Journal of Computational Physics},
  year =         {2010},
  volume =       {229},
  pages =        {9092-9118},
}

@article{Lee2013,
  author =       {C. H. Lee and A. J. Gil and J. Bonet},
  title =        {Development of a cell centred upwind finite volume
                  algorithm for a new conservation law formulation in
                  structural dynamics},
  journal =      {Computers \& Structures},
  year =         {2013},
  volume =       {118},
  pages =        {13-38},
}

@article{Trangenstein1991,
  author =       {J. A. Trangenstein and P. Colella},
  title =        {A higher-order {Godunov} method for modeling finite
                  deformation in elastic-plastic solids},
  journal =      {Communications on Pure and Applied Mathematics},
  year =         {1991},
  volume =       {44},
  pages =        {41-100},
}

@article{Tukovic2018,
  author =       {{\v{Z}}. Tukovi\'{c} and A. Kara\v{c} and P. Cardiff
                  and H. Jasak and A. Ivankovi\'{c}},
  title =        {{OpenFOAM} finite volume solver for fluid-solid interaction},
  year =         {2018},
  volume =       {42},
  number =        {3},
  pages =        {1-31},
  journal =      {Transactions of {FAMENA}},
  doi =         {10.21278/TOF.42301},
}

@article{Weller1998,
  author =       {H. G. Weller and G. Tabor and H. Jasak and
                  C. Fureby},
  title =        {A tensorial approach to computational continuum
                  mechanics using object orientated techniques},
  journal =      {Computers in Physics},
  year =         {1998},
  volume =       {12},
  pages =        {620-631},
}

\end{document}